\documentclass[aos,preprint]{imsart} 

\usepackage{amsthm,amsmath,amsfonts,amssymb,mathtools}
\usepackage{bbm}
\usepackage{dsfont}
\usepackage{tikz}
\usetikzlibrary{calc}
\let\le\leqslant
\let\ge\geqslant

\def\argmin{\mathop{\rm arg\, min}}

\def\B{{\mathcal{B}}}
\def\F{{\mathcal{F}}}
\def\mJ{{\mathcal{J}}}

\def\N{{\mathbb{N}}}
\def\mN{{\mathcal{N}}}

\def\PP{{\mathbb{P}}}
\def\hp{{\hat{p}}}
\def\R{{\mathbb{R}}}
\def\hR{{\widehat{R}}}
\def\hC{{\widehat{C}}}

\def\S{{\mathbb{S}}}
\def\Z{{\mathbb{Z}}}

\makeatletter
\newcommand*\wt[1]{\mathpalette\wthelper{#1}}
\newcommand*\wthelper[2]{%
        \hbox{\dimen@\accentfontxheight#1%
                \accentfontxheight#11.3\dimen@
                $\m@th#1\widetilde{#2}$%
                \accentfontxheight#1\dimen@
        }%
}

\newcommand*\accentfontxheight[1]{%
        \fontdimen5\ifx#1\displaystyle
                \textfont
        \else\ifx#1\textstyle
                \textfont
        \else\ifx#1\scriptstyle
                \scriptfont
        \else
                \scriptscriptfont
        \fi\fi\fi3
}
\makeatother

\makeatletter
\newcommand*\wh[1]{\mathpalette\whhelper{#1}}
\newcommand*\whhelper[2]{%
        \hbox{\dimen@\accentfontxheight#1%
                \accentfontxheight#11.3\dimen@
                $\m@th#1\widehat{#2}$%
                \accentfontxheight#1\dimen@
        }%
}

\makeatother

\newcommand{\setd}[1]{\setminus\{#1\}}
\newcommand{\indic}[1]{\mathds 1_{#1}}

\newcommand*{\ie}{i.e.\@\xspace, }
\newcommand*{\wlo}{w.l.o.g.\@\xspace}
\newcommand*{\eg}{e.g.\@\xspace, }
\newcommand*{\as}{a.s.\@\xspace}
\newcommand*{\aev}{a.e.\@\xspace}
\newcommand*{\wrt}{w.r.t.\@\xspace}
\newcommand*{\iid}{i.i.d.\@\xspace}

\newcommand*{\hmu}{\widehat \mu}

\newcommand*{\hm}{\widehat m}
\newcommand*{\tphi}{\wt{\phi}}
\newcommand*{\hphi}{\wh{\phi}}
\newcommand*{\halpha}{\wh{\alpha}}
\newcommand*{\talpha}{\wt{\alpha}}
\newcommand*{\alphas}{{\alpha_\star}}
\newcommand*{\rhos}{{\rho_\star}}

\newcommand*{\us}{{\underline s}}
\newcommand*{\ut}{{\underline t}}

\newcommand*{\hPsi}{\widehat \Psi}

\newcommand*{\push}{\mathsmaller{\#}}

\DeclareMathOperator{\Bin}{Bin}
\DeclareMathOperator{\conv}{conv}
\DeclareMathOperator{\diam}{diam}
\DeclareMathOperator{\dom}{dom}
\DeclareMathOperator{\dir}{Dir}
\DeclareMathOperator{\Id}{Id}

\DeclareMathOperator{\sgn}{sgn}
\DeclareMathOperator{\supp}{supp}
\DeclareMathOperator{\tv}{TV}

\DeclarePairedDelimiterX{\inner}[2]{\langle}{\rangle}{#1, #2}
\DeclarePairedDelimiterX{\norm}[1]{\lVert}{\rVert}{#1}
\DeclarePairedDelimiterX{\floor}[1]{\lfloor}{\rfloor}{#1}
\DeclarePairedDelimiterX{\ceil}[1]{\lceil}{\rceil}{#1}
\DeclarePairedDelimiterX{\set}[1]{\{}{\}}{#1}
\DeclarePairedDelimiterXPP{\E}[1]{\mathbb{E}}[]{}{#1}
\DeclarePairedDelimiterXPP{\V}[1]{\mathbb{V}}[]{}{#1}
\let\P\relax
\DeclarePairedDelimiterXPP{\P}[1]{\mathbb{P}}(){}{#1}
\DeclarePairedDelimiterXPP{\Pin}[1]{\mathbb{P}_*}(){}{#1}
\DeclarePairedDelimiterXPP{\Pout}[1]{\mathbb{P}^*}(){}{#1}

\usepackage[numbers]{natbib}
\usepackage[colorlinks,citecolor=blue,urlcolor=blue]{hyperref}
\usepackage{cleveref}
\usepackage{xspace}
\usepackage{xcolor}
\usepackage{float}
\usepackage{graphicx}
\graphicspath{{figures/}} 
\usepackage{caption}
\usepackage{subcaption}
\usepackage[textwidth=2.7cm, backgroundcolor=white]{todonotes}
\usepackage{enumitem}
\usepackage{dirtytalk}
\usepackage{upgreek}
\usepackage{relsize}

\startlocaldefs

\theoremstyle{plain}
\newtheorem{theorem}{Theorem}[section]
\newtheorem{lemma}[theorem]{Lemma}
\newtheorem{proposition}[theorem]{Proposition}
\newtheorem{corollary}[theorem]{Corollary}
\newtheorem{assumption}{Assumption}
\theoremstyle{definition}
\newtheorem{definition}[theorem]{Definition}
\newtheorem{example}[theorem]{Example}

\theoremstyle{definition}
\newtheorem{remark}{Remark}[section]

\newcommand*\diff{\mathop{}\!\mathrm{d}}

\endlocaldefs

\begin{document}

\begin{frontmatter}

\title{Convex generalized Fr\'echet means in a metric tree}
\runtitle{Convex Fr\'echet \texorpdfstring{$\ell$}{ℓ}-means in a metric tree}

\begin{aug}
\author[A]{\fnms{Gabriel}~\snm{Romon}\ead[label=e1]{gabriel.romon@uni.lu}}
\author[B]{\fnms{Victor-Emmanuel}~\snm{Brunel}\ead[label=e2]{victor.emmanuel.brunel@ensae.fr}}
\address[A]{Department of Mathematics, University of Luxembourg\\ 6, Avenue de la Fonte, Esch-sur-Alzette, 4364, Luxembourg\\\printead[presep={}]{e1}}
\address[B]{CREST, ENSAE, IP Paris\\ 5 Avenue Henry Le Chatelier, 91120 Palaiseau, France\\\printead[presep={}]{e2}}

\end{aug}

\begin{abstract}
    We are interested in measures of central tendency for a population $\mu$ on a network, which is modeled by a metric tree.
    The location parameters that we study are generalized Fr\'echet means defined as minimizers of the objective function $\alpha \mapsto \E{\ell(d(\alpha,X)) - \ell(d(o,X))}$,
    where $\ell$ is a convex, strictly increasing loss, $X\sim \mu$ and $o$ is an arbitrary origin.
     
    We leverage the geometry of the tree and the geodesic convexity of the objective to develop a notion of directional derivative in the tree, which helps us locate and characterize the minimizers.
    We then extend to a metric tree the concept of stickiness defined by Hotz et al.\ \cite{hotz2013sticky}.
    When $\mu$ is nondegenerate,
    any Fr\'echet $\ell$-mean of $\mu$ is either sticky, one-sided partly sticky or two-sided partly sticky,
    depending on the number of vanishing directional derivatives.

    Estimation is performed using a sample analog. 
    We establish laws of large numbers and central limit theorems 
    that reveal distinct asymptotic behaviors of empirical means across these three regimes: 
    collapse onto the population mean in the sticky case, 
    partial collapse accompanied by one-sided fluctuations in the one-sided case, 
    and two-sided fluctuations without collapse in the two-sided case.
    We further develop two consistent estimators of the stickiness regime, 
    based respectively on pairwise collisions among empirical means 
    and on empirical difference quotients.
    For the particular case of the Fr\'echet median, 
    we develop distribution-free non-asymptotic confidence regions
    for the entire set of minimizers.
\end{abstract}

\end{frontmatter}

\section{Introduction}

\subsection{Context}

Across many domains of application, it is common to model observations as points on a network.
These points may lie not only at the vertices of the network, but also in the interior of the edges.
The network is usually represented as a linear network: 
its edges are line segments of prescribed lengths, 
these segments are joined at their endpoints, and the resulting space is equipped with the shortest path metric.

A common example in spatial statistics is the analysis of traffic accidents along a road network.
The accident locations are modeled as a point pattern on the network.
The analysis can then be used to identify stretches of road where accident risk is elevated \cite{baddeley2021analysing}. 
Another example arises in the study of neurons. Spines are small protrusions distributed along the dendritic network of a neuron, 
and their locations are analyzed using point-process methods. 
Questions of interest include whether spine density varies among branches and whether spines exhibit clustering 
\cite{jammalamadaka2013statistical,baddeley2014multitype}.

In metagenomics, an environmental sample can be studied by sequencing a large number of short DNA fragments, called reads, 
from the collection of microbes present in the sample. 
A fixed reference phylogenetic tree represents the evolutionary relationships among previously characterized DNA sequences, 
with edge lengths encoding amounts of evolutionary change. 
Phylogenetic placement methods assign each read to a position on this reference tree,
according to its inferred evolutionary relationship with the reference sequences. 
Thus, the resulting placements form a collection of points $x_1,\ldots,x_n$ on the reference tree, 
and the environmental sample may be represented by the empirical measure $\frac 1n \sum_{i=1}^n \delta_{x_i}$. 
Given two samples from different environments, the $1$-Wasserstein distance between their empirical measures on a common reference tree 
quantifies the dissimilarity between the two microbial communities \cite{evans2012phylogenetic}.

In operations research, points on a transportation or distribution network may represent users requiring service or, more generally, sources of demand. 
A probability distribution on the network then describes where a future service request is likely to occur. 
A common objective in the literature is to place one or more facilities 
in order to minimize the expected distance from a random demand point to its nearest facility. 
Facility location models often assume that demand is concentrated at the vertices. 
This assumption can be unrealistic, for instance when positioning railway maintenance facilities, 
since a track failure may occur at any point along the railway network.
Continuous-demand models therefore allow demand to arise along edges, 
possibly with additional mass at the vertices \cite{nkansah1986network,brandeau1986locating,blanquero2013solving}.

Dendritic, phylogenetic, and river networks are naturally without cycles. 
In other settings, cycles may be avoided for economic reasons, as in sparse transportation or communication networks.
When the network is connected and contains no cycles, there is a unique path between any pair of points, 
and the resulting linear network is called a metric tree. 
Basic examples include the real line 
(viewed as infinitely many line segments chained together)
and star-shaped networks organized around a single hub.

These examples motivate a population-level model 
in which the ambient metric space is a metric tree $T$ equipped with the shortest path metric $d$,
and the distribution of locations is represented by a probability measure $\mu$.
It is then natural to seek a point of $T$ that summarizes the central tendency of $\mu$.
For a convex and strictly increasing loss $\ell:[0,+\infty)\to\R$
and an origin $o\in T$ such that $\int_T \ell_+'(d(o,x)) \diff\mu(x) < +\infty$,
we consider minimizers of the objective function
$\alpha \mapsto \int_T  \ell(d(\alpha,x)) - \ell(d(o,x))  \diff\mu(x)$, which we call Fr\'echet $\ell$-means.
Two special cases of this formulation have already appeared in the applications above.
For $\ell:z\mapsto z$, Fr\'echet medians coincide with solutions of the continuous-demand one-median problem studied in operations research
\cite{nkansah1986network,blanquero2013solving}.
For $\ell:z\mapsto z^2$, the resulting parameter is the Fr\'echet mean, or barycenter.
Evans and Matsen \cite[Section~5.2]{evans2012phylogenetic} 
note that it provides a one-point summary of a probability distribution on a phylogenetic tree.
In the statistical model studied here, $\mu$ is unknown and one observes an \iid sample $X_1,\ldots,X_n$ from $\mu$.
The corresponding population parameter is naturally estimated by replacing $\mu$ with the empirical measure
$\hmu_n=\frac 1n \sum_{i=1}^n \delta_{X_i}$,
or, equivalently, by minimizing the empirical objective
$\alpha \mapsto \sum_{i=1}^n \ell(d(\alpha,X_i))$.

Unlike the real line, a metric tree with a branching point has no natural linear structure. 
The estimation problem considered above thus belongs to the broader field of non-Euclidean statistics.
The notion of a mean on a metric space dates back to Fr\'echet \cite{frechet1948elements}. 
Since then, an extensive literature has developed on the statistical theory of empirical Fr\'echet means. 
By contrast, results covering losses beyond $\ell:z\mapsto z^2$ remain comparatively scarce. 
Laws of large numbers have been established for broad classes of generalized Fr\'echet means 
\cite{huckemann2011intrinsic,schotz2022strong,aveni2026uniform};
convergence rates \cite{schotz2019convergence,schotz2026transformed}
and variance inequalities \cite{schotz2025variance} have also been obtained. 
On Riemannian manifolds, asymptotic normality has been established for geodesically convex $M$-estimators \cite{brunel2023geodesically}.
A metric tree is a particular instance of a Hadamard space \cite{bridson1999metric}, 
hence results on Fr\'echet means
\cite{sturm2003probability,mccormack2022stein,brunel2024concentration,escande2024concentration,lammers2023sticky}
and generalized Fr\'echet means 
\cite{schotz2025variance,schotz2026transformed}
in general Hadamard spaces
apply also to metric trees. 

There is little statistical literature on Fr\'echet $\ell$-means in the specific setting of a metric tree. 
Funano \cite{funano2010rate} focuses on the Fr\'echet mean in an $\R$-tree and
obtains a Gaussian tail bound for the inductive mean, a different estimator from the empirical Fr\'echet mean.
Basrak \cite{basrak2010limit} establishes a central limit theorem for the inductive mean in a binary metric tree. 
Hotz et al.\ \cite{hotz2013sticky} develop laws of large numbers and central limit theorems for the Fr\'echet mean on an open book, 
\ie finitely many Euclidean halfspaces glued along their common boundary. 
In the special case where the pages are one-dimensional, hence Euclidean rays, the resulting open book is itself a metric tree.

\subsection{Contributions and outline}

The goal of this work is to 
develop a theory of Fr\'echet $\ell$-means on metric trees,
covering both the population and the statistical standpoint.
We describe below how the paper is organized and we give a brief summary of our contributions.

\begin{itemize}
    \item In \Cref{sec:frechet},
    by leveraging the geodesic convexity of the objective function (\Cref{prop:phi-trees}) 
    and the geometry of the tree (\Cref{def:directions}),
    we develop a notion of directional derivative (\Cref{def:directional-derivative}).
    This allows us to describe the shape of the set of minimizers (\Cref{thm:shape}), 
    as well as to locate (\Cref{prop:sign-derivative}) and characterize (\Cref{prop:optimality}) 
    Fr\'echet $\ell$-means according to the signs of these directional derivatives.

    \item In \Cref{sec:stickiness}, we extend to the metric tree 
    the notion of stickiness introduced by Hotz et al.\ \cite{hotz2013sticky}.
    An arbitrary point $\alpha \in T$ is either sticky, partly sticky or nonsticky 
    according to the signs of directional derivatives at $\alpha$
    (\Cref{def:sticky}).
    Depending on the stickiness regime of $\alpha$, 
    we refine the localization constraints on minimizers (\Cref{thm:localization-sticky}).
    According to the number of vanishing derivatives, we further divide the partly sticky category into two subcategories: 
    one-sided partly sticky and two-sided partly sticky (\Cref{def:partly-sticky-sided}).
    We then draw a connection between the stickiness regime of $\alpha$ 
    and the geometry of the tree (\Cref{prop:sticky-geometry}).
    Finally, we provide an equivalent definition of stickiness, 
    in terms of robustness to small perturbations 
    of the data-generating distribution (\Cref{thm:robust-sticky}).

    \item In \Cref{sec:stats}, we turn to the estimation of Fréchet $\ell$-means using a sample analog.
    We first establish a set-valued law of large numbers 
    (\Cref{thm:lln-sets}).
    We then show that, in the sticky case, 
    the empirical Fr\'echet $\ell$-means eventually collapse 
    exactly onto the population parameter (\Cref{thm:sticky-lln}). 
    We complement this asymptotic result with non-asymptotic bounds 
    (\Cref{thm:sticky-nonasymp-subgaussian}).

    Finally, assuming that the Fr\'echet $\ell$-mean is unique and partly sticky, 
    we establish central limit theorems in the two-sided and one-sided regimes.
    In the two-sided case, there is no collapse at the population mean: 
    empirical means fluctuate along the two zero-derivative directions 
    (\Cref{thm:two-sided-clt}).
    In the one-sided case, empirical Fr\'echet $\ell$-means instead exhibit partial collapse, 
    being equal to the population parameter with asymptotic probability $1/2$ 
    and otherwise fluctuating into the unique zero-derivative direction 
    (\Cref{prop:one-sided-collapse} and \Cref{thm:one-sided-clt}).

    \item In \Cref{sec:estimate-stickiness}, 
    we consider the problem of estimating the stickiness regime of a Fr\'echet $\ell$-mean.
    We introduce a first estimator based on 
    the proportion of pairwise collisions among empirical means 
    (\Cref{def:estimator-collisions}), 
    and a second estimator based on empirical difference quotients 
    (\Cref{def:estimator-difference-quotients}).
    We then establish consistency of both estimators
    (\Cref{corol:estimator-stickiness},\Cref{thm:estimator-stickiness-2}).

    \item In \Cref{sec:medians}, we specialize our results to Fr\'echet medians.
    We further study the set of population medians,
    deriving constraints on its endpoints (\Cref{prop:endpoints-median})
    and characterizing when the median is unique (\Cref{prop:unique-median}).
    Finally, we construct distribution-free, non-asymptotic confidence regions
    for the whole set of Fr\'echet medians (\Cref{thm:confidence-region})
    and show that, when the population median is unique,
    they shrink to it almost surely (\Cref{prop:region-shrink}).
\end{itemize}

\subsection{Notations and terminology}
$\subset$ means inclusion, possibly with equality.
We let $\N = \set{0,1,\ldots}$ and $[n] = \set{1,\ldots,n}$ for each integer $n\ge 1$.
The cardinality of a set $A$ is denoted by $|A|$.
We say that $A$ is countable if it is finite or countably infinite.
Throughout, $T$ denotes the metric tree
and $d$ is the shortest path metric. 
Given $x\in T$ and $r>0$, 
$B(x,r)$ (resp., $\overline B(x,r)$) is the open (resp., closed) ball with center $x$ and radius $r$.
The interior of $A$ is denoted by $\operatorname{int}(A)$. 
We write $\lambda$ for the Lebesgue measure on $T$ or $\R$,
and $\delta_x$ for the Dirac measure centered at $x\in T$.
$N(\mu, \sigma^2)$ is the normal distribution with mean $\mu\in \R$ and variance $\sigma^2>0$.
$\Bin(n,p)$ is the binomial distribution with parameters $n\in \N$ and $p\in [0,1]$.
Given a function $f:\R \to \R$ and $x\in \R$,
$f(x+)$ (resp., $f(x-)$) denotes the right-hand (resp., left-hand) limit of $f$ at $x$,
while
$f_+'(x)$ (resp., $f_-'(x)$) is the right (resp., left) derivative of $f$ at $x$.
The terms \say{increasing} and \say{decreasing} do not imply strictness; when strictness is intended, we write \say{strictly increasing} and \say{strictly decreasing}.
Given $x\in \R$, we let $x_+=\max(x,0)$.

\section{Fr\'echet \texorpdfstring{$\ell$}{ℓ}-means in a metric tree}
\label{sec:frechet}

\subsection{Terminology and setting}

Let us make precise what we mean by a metric tree and introduce further terminology.

\begin{definition}
    \label{def:tree}
    \begin{enumerate}
        \item 
        Let $T$ denote an undirected, connected (meaning that any two distinct vertices are linked by a finite path) and acyclic graph with weighted edges.
        The weight of an edge is understood as the length of this edge, i.e., as the distance between the corresponding adjacent
        vertices.
        We make the following assumptions on $T$:
        \begin{enumerate}
            \item There are at least two vertices.
            \item Each vertex has finitely many neighbors.
            \item Each edge length belongs to $(0,+\infty)$.
        \end{enumerate}
        Moreover, each edge is identified with a compact interval of $\R$ that has the same length.
        The tree $T$ is then equipped with the shortest path metric $d$: the distance between two points of $T$ (not necessarily vertices) is the length of the shortest path between them.
        Then, $(T,d)$ is a metric space, which is referred to as a \textit{metric tree}. 
        In the sequel, we denote by $\mathcal B(T)$ its Borel $\sigma$-algebra.

        \item The metric tree with vertex set $\N$, edge set $\set{\set{n,n+1}:n\in \N}$, and weight $1$ on each edge is called the \textit{ray}.
        The metric tree with vertex set $\Z$, edge set $\set{\set{n,n+1}:n\in \Z}$, and weight $1$ on each edge is called the \textit{double ray}.

        \item Let $m\ge 2$. $T$ is an \textit{$m$-spider} if the underlying graph-theoretic tree has exactly $m$ leaves and there is a single vertex adjacent to all of them.
        $T$ is an \textit{open book} if it consists of $m$ rays glued together at their origin.
    \end{enumerate}
\end{definition}

\begin{remark}
    In order to shorten the exposition, we glossed over some technicalities when defining the concept of metric tree.
    A comprehensive construction can be found in \cite[Section~I.1.9]{bridson1999metric} or \cite[Examples 2.1.3 (iv)]{papadopoulos2014metric}.
\end{remark}

\begin{remark}
    The real line $\R$ is isometric to the double ray. Throughout this work, we recover the classical theory
    for one-dimensional parameters of central tendency
    by considering the double ray as the ambient tree $T$. 
\end{remark}

As seen in \Cref{lemma:countable} of the Appendix, $T$ has countably many vertices, hence also countably many edges. We are then able to define an analog of the Lebesgue measure on $(T,d)$.

\begin{definition} 
    \label{def:Lebesgue}
    Let $\mathcal E$ denote the set of edges of $T$. 
    Each edge $e\in\mathcal E$ being identified with a compact interval $I_e\subset \R$ of identical length, 
    $e$ naturally inherits the Lebesgue measure from $I_e$, denoted by $\lambda_e$. 
    For any $A\in \mathcal B(T)$, we set $\lambda(A)=\sum_{e\in\mathcal E} \lambda_e(A\cap e)$, where $A\cap e$ is identified isometrically with a subset of $I_e$. 
\end{definition}

Next, we introduce some relevant concepts from metric geometry.
Given $x,y\in T$, a constant speed geodesic from $x$ to $y$ is a map $\gamma$ from some interval $[a,b]\subset \R$  to $T$ such that 
$\gamma(a)=x$, $\gamma(b)=y$ and $d(\gamma(t_1), \gamma(t_2)) = v |t_1-t_2|$ for some $v\in [0,+\infty)$ and every $t_1,t_2\in [a,b]$.
If $a\neq b$, $v=\frac{d(x,y)}{b-a}$ is called the speed of the geodesic $\gamma$.
For the sake of legibility, we will often write $\gamma_t$ in lieu of $\gamma(t)$.
The space $(T,d)$ is uniquely geodesic, meaning that between any two points $x,y\in T$, there exists a geodesic from $x$ to $y$ that is unique up to reparametrization. 
Its image is denoted by $[x,y]$ and it is referred to as the geodesic segment joining $x$ and $y$.
We also define open and half-open geodesic intervals $(x,y)$, $[x,y)$, $(x,y]$. For instance, $[x,y) = \gamma([a,b))$
for some $a\le b$ and a geodesic $\gamma$ from $x$ to $y$ defined on $[a,b]$.

In addition to the assumptions in \Cref{def:tree}, we assume in the rest of the paper that the metric tree satisfies the following condition.
\begin{enumerate}[label=\alph*), ref=\alph*)]
    \setcounter{enumi}{3}
    \item \label{assum:lower-bound}
    The set of all edge lengths is bounded below by some positive constant.
\end{enumerate}

This assumption on edge length is key in establishing the following topological and geometric properties of $(T,d)$.
\begin{proposition}
    \label{prop:topology}
    \begin{enumerate}
        \item $(T,d)$ is a proper metric space, meaning that closed and bounded subsets of $T$ are compact. 
        As a consequence, $(T,d)$ is complete, separable and locally compact.
        \item $(T,d)$ is a Hadamard space (see \cite[Chapter II.1]{bridson1999metric} for the definition).
    \end{enumerate}
\end{proposition}
In Hadamard spaces it is possible to develop a theory of convex analysis, convex optimization and probability
that generalizes to nonlinear settings the classical results known in Hilbert spaces \cite{bacak2014convex}. 
Let us recall the definition of geodesic convexity in $(T,d)$.

\begin{definition} 
    \begin{enumerate}
        \item A subset $G\subset T$ is called geodesically convex (convex, for short) if for all $x,y\in G$, the geodesic segment $[x,y]$ is a subset of $G$.
        \item A function $f:T\to(-\infty,+\infty]$ is called geodesically convex (convex, for short) if 
        for every geodesic $\gamma:[0,1]\to T$ and $t\in (0,1)$ we have the inequality $f(\gamma_t)\le (1-t)f(\gamma_0)+tf(\gamma_1)$.
        The set $\dom f \coloneqq \set{\alpha\in T:f(\alpha)< +\infty}$ is called the domain of $f$.
        When $\dom f$ is nonempty, we say that $f$ is proper.
        We call $f$ geodesically strictly convex (strictly convex, for short) if the previous inequality is strict, provided that 
        $\gamma_0\neq \gamma_1$ and 
        $\set{\gamma_0, \gamma_1}\subset \dom f$. %
    \end{enumerate}
\end{definition}

We next introduce the concept of loss function and we define parameters of central tendency.

\begin{definition}
    \label{def:objective}
    A \textit{loss function} is a mapping $\ell:[0,+\infty)\to \R$ that is convex and strictly increasing.
    Let $\mu$ be a probability measure on $(T,\mathcal B(T))$. 
    We define the subset $F=\set{\alpha\in T: \int_T \ell_+'(d(\alpha,x)) \diff \mu(x) < \infty}$
    and we pick an arbitrary $o\in T$, which we call the \textit{origin}.
    
    We assume from now on that $o\in F$, and we define the \textit{objective function} $\phi$ 
    \begin{align}
        \label{eq:objective}
        \phi \colon T &\to (-\infty,+\infty] \\
        \alpha &\mapsto \int_T \ell(d(\alpha,x))-\ell(d(o,x)) \diff\mu(x). \nonumber
    \end{align}
    Minimizers of $\phi$ are called \textit{Fr\'echet $\ell$-means of $\mu$}, and 
    we denote by $M(\mu)$ the set of all minimizers.
\end{definition}

\begin{remark}
    The notation $F$ is mnemonic for \say{finite}.
    By the convexity of $\ell$ and the triangle inequality, the following estimate holds:
    $-\ell_+'(d(o,x)) d(\alpha, o) \le \ell(d(\alpha,x)) - \ell(d(o,x))$ for every $\alpha, x \in T$.
    For the integral in \eqref{eq:objective} to be well-defined in $(-\infty,+\infty]$, it suffices therefore that $o$ belongs to $F$,
    which is often less stringent in practice than requiring instead integrability of $x\mapsto \ell(d(\alpha,x))$.
    This motivates subtracting $\ell(d(o,x))$ in the definition of $\phi$.
\end{remark}

\begin{remark}
    The choice of origin $o\in F$ does not affect the set of minimizers $M(\mu)$.
    Indeed, choosing a different origin in $F$ only shifts the objective function $\phi$ by some additive constant.
\end{remark}

In the next proposition, we provide foundational properties of $\phi$ and $M(\mu)$.

\begin{proposition}
    \label{prop:phi-trees}
    \begin{enumerate}
        \item $\phi$ is well-defined, convex, proper and $F\subset \dom \phi$.
        \item $\phi$ is lower semicontinuous on $T$ and continuous on $\operatorname{int}(\dom \phi)$.
        \item If $\ell_+'$ is bounded from above, then $$\frac{\phi(\alpha)}{\ell(d(\alpha,o))}\xrightarrow[d(\alpha,o)\to +\infty]{} 1.$$
        Otherwise, $\ell_+'(z)\xrightarrow[z\to +\infty]{} +\infty$ and  $$\frac{\phi(\alpha)}{d(\alpha,o)}\xrightarrow[d(\alpha,o)\to +\infty]{} +\infty.$$
        In both cases, $\phi$ is coercive, meaning that $\phi(\alpha) \to +\infty$ as $d(\alpha,o)\to +\infty$.
        \item $M(\mu)$ is a nonempty, compact and convex subset of $T$.
        \item $\ell$ is strictly convex if and only if $\ell_+'$ is strictly increasing. 
        In that case, $\phi$ is strictly convex and $M(\mu)$ is a singleton.
    \end{enumerate}
\end{proposition}

\begin{remark}
    Convex losses which are increasing but not strictly increasing lead to location parameters that do not denote centrality.
    Indeed, if $\ell$ is increasing but not in the strict sense, it is constant on an interval of the form $[0,r]$ where $r>0$.
    Then, letting $\mu=\delta_{o}$ be the Dirac measure at $o$, we have the inclusion
    $\set{\alpha\in T: d(\alpha,o)\le r}\subset M(\mu)$ and we obtain location parameters which are arguably not representative of 
    the central tendency of $\delta_{o}$. It is therefore not restrictive in \Cref{def:objective} to only consider losses that are strictly increasing.
\end{remark}

Recall that we silently assume throughout that $F$ contains at least the origin, \ie $o\in F$. 
Since the condition $\alpha \in F$ appears repeatedly in later results, 
it is useful to know when it holds automatically.
The following assumption provides a favorable setting where $F=T$.
We let $\supp \mu$ denote the support of $\mu$.
Since $T$ is a separable metric space, $\supp \mu$ is well-defined as 
the smallest closed subset of $T$ that has $\mu$-probability $1$ \cite[Proposition~7.2.9]{bogachev2007measure}.

\begin{assumption}
    \label{assum:bounded}
    The subset $\supp \mu$ is bounded (\eg $T$ is bounded), or the function $\ell_+'$ is bounded.
\end{assumption}

To deal with the unbounded case, we formulate an assumption on the growth of $\ell_+'$.

\begin{assumption}
    \label{assum:growth}
    $\exists C>0, \exists R>0, \forall z\ge R, \, \ell_+'(z+1)\le C\ell_+'(z)$.
\end{assumption}

For losses that satisfy \ref{assum:bounded} or \ref{assum:growth}, the integrability of  $x\mapsto \ell_+'(d(o,x))$ implies
integrability of $x\mapsto \ell_+'(d(\alpha,x))$ for any $\alpha\in T$, as seen in the following proposition.
Moreover, we provide sufficient conditions which ensure that \Cref{assum:growth} holds.

\begin{proposition}
    \label{prop:growth}
    \begin{enumerate}
        \item Under Assumption~\ref{assum:bounded} or \ref{assum:growth}, $F=\dom \phi = T$.
        \item If $\ell_+'$ is concave, then it satisfies \Cref{assum:growth}.
        \item Suppose that $\ell_+'(z)\asymp \exp\big(C_1 (\log z)^{C_2} + C_3 z\big)$, in the sense that 
        $$\exists m,M > 0, \exists R>1, \forall z\ge R, \, m \le \ell_+'(z) \big/ \exp\big(C_1 (\log z)^{C_2} + C_3 z\big) \le M,$$
        where $C_1\ge 0$, $C_2>0$ and $C_3\ge 0$.
        Then \Cref{assum:growth} holds.
    \end{enumerate}
\end{proposition}

\begin{example}
    \label{ex:ell}
    Examples of loss functions $\ell$ include: 
    \begin{enumerate}
        \item $\ell:z\mapsto z^p$ where $p\in [1,+\infty)$. In this setting, the minimizers of $\phi$ are called \textit{Fr\'echet $p$-means} of $\mu$.
        In the case $p=1$ they are referred to as \textit{Fr\'echet medians}, and when $p=2$ as \textit{barycenters} or just \textit{Fr\'echet means}.
        The corresponding set of minimizers will be denoted specifically by $M_p(\mu)$.
        When $p=1$, we have $\ell_+'=1$, hence \Cref{assum:bounded} and \ref{assum:growth} hold.
        For any $p>1$, $\ell_+'$ is unbounded but the loss satisfies \ref{assum:growth}.
        Indeed $\ell_+'(z)\asymp \exp((p-1)\log z)$. 
        When $F$ is nonempty, it is usually said that $\mu$ has a finite moment of order $p-1$.

        \item $\ell: z\mapsto z^2 \indic{z \le c} + (2cz - c^2) \indic{z > c}$ where $c>0$.
        It is known as the Huber loss.
        Note that $\sup_{z\ge 0}\ell_+'(z) = 2c <\infty$, hence \ref{assum:bounded} and \ref{assum:growth} hold.

        \item $\ell: z\mapsto 2c^2\big((1+z^2/c^2)^{1/2} -1\big)$ where $c>0$.
        It is known as the pseudo-Huber loss, which is a smooth approximation of the standard Huber loss.
        Similarly to the Huber loss, $\sup_{z\ge 0}\ell_+'(z) = 2c$, hence both \ref{assum:bounded} and \ref{assum:growth} hold.
        
        \item $\ell:z\mapsto \exp(z)$. The derivative $\ell_+'$ is unbounded, however the loss satisfies \ref{assum:growth}.
        \item $\ell:z\mapsto (1+z)\log(1+z) - z$. In that case, $\ell_+':z \mapsto \log(1+z)$ is concave; it is not bounded
        but \ref{assum:growth} is satisfied.
        \item $\ell:z\mapsto \exp(\exp(z))$. The derivative $\ell_+'$ is not bounded and \ref{assum:growth} is not met.
    \end{enumerate}
\end{example}

\subsection{Directions in a metric tree}
\label{sec:directions}

We leverage the geometry of the metric tree $T$ to define a notion of direction relative to a point.
A related concept in a generic metric space is the space of directions \cite[Definition II.3.18]{bridson1999metric}.

\begin{definition}
    \label{def:directions}
    Fix $\alpha\in T$. On the set $T\setd{\alpha}$ we define the equivalence relation $\sim_\alpha$ as follows:
    $$\forall (v,w)\in (T\setd{\alpha})^2,\quad  v\sim_\alpha w \iff [\alpha, v]\cap [\alpha, w]\neq \set{\alpha}.$$
    We let $\dir_\alpha$ denote the quotient set of $T\setd{\alpha}$ by $\sim_\alpha$,
    and we refer to each equivalence class as a \textit{direction at $\alpha$}.
    Given $v\in T\setd{\alpha}$, we write $\Delta_{\alpha \to v}$ for the equivalence class that contains $v$.
\end{definition}

\begin{figure}[H]
    \centering
    \begin{subfigure}{.5\textwidth}
      \centering
      \includegraphics[]{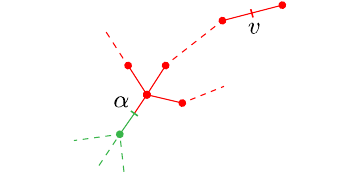}
      \caption{$\alpha$ is in the interior of an edge}
    \end{subfigure}%
    \begin{subfigure}{.5\textwidth}
      \centering
      \includegraphics[]{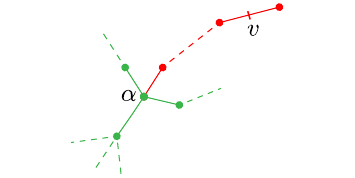}
      \caption{$\alpha$ is a vertex}
    \end{subfigure}
    \caption{Illustration for \Cref{def:directions} in two cases. $\Delta_{\alpha \to v}$ is drawn in red and $T\setminus \Delta_{\alpha \to v}$ is drawn in green.}
    \label{fig:directional-tree}
\end{figure}

Note that an equivalent formulation for $v\sim_\alpha w$ is $\alpha\notin [v,w]$.
\Cref{fig:directional-tree} illustrates the definition in two situations: either $\alpha$ is in the interior of an edge,
or $\alpha$ is a vertex of $T$. We stress that $\alpha$ does not belong to $\Delta_{\alpha \to v}$.

The following proposition shows that there are finitely many directions at $\alpha$.
It also provides an alternative interpretation of the direction $\Delta_{\alpha \to v}$.
Let us introduce the additional terminology of \textit{neighboring vertices}.
If $\alpha$ is a vertex, we refer to the graph-theoretic neighbors of $\alpha$.
If $\alpha$ lies in the interior of an edge, we refer to the endpoints of this edge.

\begin{proposition}
    \label{prop:directions}
    Let $\alpha\in T$.
    \begin{enumerate}
        \item Each direction $\Delta\in \dir_\alpha$ is a nonempty convex subset of $T$.
        \item Letting $v_1,\ldots,v_m$ denote the neighboring vertices of $\alpha$,  
        we have $\dir_\alpha = \set{\Delta_{\alpha\to v_i}:1\le i\le m}$.
        \item The metric space $(T\setd{\alpha},d)$ has finitely many path components and 
        for each $v\in T\setd{\alpha}$, $\Delta_{\alpha \to v}$ is the path component that contains $v$.
    \end{enumerate}
\end{proposition}

In view of \Cref{prop:phi-trees}, $M(\mu)$ is a nonempty compact convex subset of $T$.
By leveraging the concept of directions, we obtain the following more precise statement on the geometry of $M(\mu)$.

\begin{theorem}
    \label{thm:shape}
    The set of Fr\'echet $\ell$-means $M(\mu)$ is a geodesic segment.
\end{theorem}

\begin{example}
    When $p>1$, the loss defining the Fr\'echet $p$-mean is strictly convex, hence by \Cref{prop:phi-trees} $M_p(\mu)$ is a singleton 
    and thus a trivial instance of a geodesic segment.
    In the case of Fr\'echet medians $(p=1)$ or the Huber loss, $\ell$ is not strictly convex. 
    Besides, any geodesic segment can be realized as $M_1(\mu)$ for some $\mu$, as will be seen in \Cref{sec:medians-descriptive}. 
\end{example}

\subsection{Directional derivatives, localization of minimizers and optimality conditions}
\label{sec:calculus-tree}
The following subsection is devoted to locating and characterizing the minimizers of $\phi$.

Given a convex function $f:\R \to (-\infty,+\infty]$, the variations of $f$ with respect to a reference point $\alpha \in \dom f$
are naturally assessed by defining the difference quotient 
$q:\beta\mapsto \frac{f(\beta) - f(\alpha)}{\beta-\alpha}$ and considering its left-hand limit 
$q(\alpha-)$
and its right-hand limit 
$q(\alpha+)$.
In the metric tree $T$, if $\alpha\in \dom \phi$ has $m$ neighboring vertices, then $m$ directions emanate from $\alpha$, as seen in \Cref{prop:directions}.
We adapt the notion of difference quotient to the tree by defining the \textit{metric difference quotient}  
\begin{align}
    \label{eq:diff-quotient-metric}
    Q \colon T\setd{\alpha} &\to (-\infty,+\infty] \\
    \beta &\mapsto \frac{\phi(\beta) - \phi(\alpha)}{d(\beta,\alpha)}. \nonumber
\end{align}
In contrast to $q$, since the tree lacks an ordering, the denominator of $Q$ is always positive.
In order to extend the notions of left and right derivative, we consider a direction $\Delta \in \dir_\alpha$ and 
we assess the limit $\lim_{\beta \to \alpha, \beta \in \Delta} Q(\beta)$.

\begin{proposition}
    \label{prop:derivative}
    Let $\alpha\in \dom \phi$ and let $\Delta\in \dir_\alpha$.
    \begin{enumerate}
        \item The limit $\lim\limits_{\beta \to \alpha, \beta \in \Delta} Q(\beta)$ exists as a quantity in $[-\infty,+\infty]$. We denote it by $\phi'_{\Delta}(\alpha)$.
        \item Assume that $\alpha \in F$ and $\Delta\cap F \neq \emptyset$.
        Then $\phi'_{\Delta}(\alpha)$ is finite with the closed form: 
        \begin{equation}
            \label{eq:derivative}
            \phi'_{\Delta}(\alpha) = 
            \int_{T\setminus \Delta} \ell_+'(d(\alpha,x))\diff\mu(x)
                - \int_{\Delta} \ell_-'(d(\alpha,x)) \diff\mu(x)
                .
        \end{equation}
        In particular, \eqref{eq:derivative} holds whenever $F=T$, \eg as soon as \Cref{assum:bounded} or \ref{assum:growth} is satisfied.
        \item If $\alpha \in F$, then $\phi'_{\Delta}(\alpha)\in (-\infty, +\infty]$. 
        If $\alpha \in F$ and $\phi'_{\Delta}(\alpha) < +\infty$, then \eqref{eq:derivative} holds.
    \end{enumerate}
\end{proposition}

\begin{definition}
    \label{def:directional-derivative}
    \begin{enumerate}
        \item We refer to the limiting value $\phi'_{\Delta}(\alpha)$
        as \textit{the directional derivative of $\phi$ at $\alpha$ along $\Delta$}.
        \item Given $v\in T\setd{\alpha}$, we write $\phi'_{v}(\alpha)$ as a shorthand for $\phi'_{\Delta_{\alpha \to v}}(\alpha)$.
    \end{enumerate}
\end{definition}

\begin{example}
    \begin{enumerate}
        \item For Fr\'echet $p$-means where $p>1$, it holds that 
        $$\phi'_{\Delta}(\alpha) = p\Big(\int_{T\setminus \Delta} d(\alpha,x)^{p-1}\diff\mu(x)
        - \int_{\Delta} d(\alpha,x)^{p-1} \diff\mu(x)\Big).$$
        Regarding the convergence of those integrals, recall that we assume from the outset that $o\in F$, \ie $\int_T d(o,x)^{p-1}\diff\mu(x)$
        is finite.
        The loss at hand satisfies \Cref{assum:growth}, hence $F=T$.

        \item For Fr\'echet medians ($p=1$ above),
        $\phi'_{\Delta}(\alpha) = \mu(T\setminus \Delta) - \mu(\Delta)
        = 1 - 2\mu(\Delta)
        =  2 \mu(T\setminus \Delta) - 1.$
        Remarkably, this directional derivative does not involve the metric $d$; it is expressed solely in terms of $\mu$.
    \end{enumerate}
\end{example}

\begin{remark}
    Assessing the directional derivative of $\phi$ along an edge is not a new idea.
    The computation is performed in \cite{shier1983optimal,brandeau1988para} for the empirical Fr\'echet $p$-mean, 
    in \cite{sturm2003probability} for the Fr\'echet mean on an $m$-spider,
    and in \cite{evans2012phylogenetic} for the Fr\'echet mean.
\end{remark}

\begin{remark}
    The closed form \eqref{eq:derivative} for $\phi'_{\Delta}(\alpha)$
    will be a pivotal ingredient in some of the proofs later on.
    We next demonstrate that \eqref{eq:derivative} may cease to hold when $\Delta \cap F=\emptyset$.
    Indeed, consider the case where $T$ is the double ray (which we identify with $\R$), $\ell:z\mapsto \exp(\exp(z))$
    and define $\mu$ as the distribution of the random variable $X = -\log Y$, where $Y$ has density $f_Y:y\mapsto C e^{-y}y^{-3} \indic{y\ge 1}$
    and $C$ is the normalizing constant.
    The origin we choose is $o = 0$ and we observe that $\int_T \ell_+'(d(o,x)) \diff \mu(x) = \E{Ye^Y} = C \int_{1}^{\infty} y^{-2}\diff y < +\infty$, hence $o\in F$.

    Pick $\beta\in [0,+\infty)$ and note that $X$ is supported on $(-\infty,0]$ hence $d(\beta, X) = \beta - X = \beta + \log Y$ \as
    and $\int_T \ell(d(\beta, x))\diff\mu(x) = C\int_{1}^{\infty}\exp((e^\beta-1)y) y^{-3} \diff y$, which is finite if and only if $\beta=0$.
    Thus $\phi(\beta) = +\infty$ for each $\beta>0$.

    We consider $\Delta=(0,\infty)$ as the direction at $\alpha = 0$ and 
    it follows from $\phi(\beta) = +\infty$ that $Q(\beta) = +\infty$ for each $\beta \in \Delta$, hence $\phi'_{\Delta}(0) = +\infty$.
    Observe that $\int_T \ell_+'(d(\beta,x))\diff\mu(x) = e^\beta\E{Y \exp(e^\beta Y)}=+\infty$ for each $\beta \in \Delta$, 
    hence $\Delta \cap F=\emptyset$.
    Besides, $\int_{T\setminus \Delta} \ell_+'(d(0,x))\diff\mu(x) = \E{Ye^Y} < \infty$ and 
    $\int_{\Delta} \ell_-'(d(0,x)) \diff\mu(x) = 0$, thus the right-hand side of \eqref{eq:derivative} is well-defined
    and finite, but it is not equal to $\phi'_{\Delta}(0)$.
\end{remark}

Next, we leverage the geometry of $T$ and the geodesic convexity of $\phi$ to show a connection between the sign of the directional derivative and the location of the minimizers of $\phi$.

\begin{proposition}
    \label{prop:sign-derivative}
    Let $\alpha\in \dom \phi$ and $\Delta\in \dir_\alpha$.
    \begin{enumerate}
        \item If $\phi'_\Delta(\alpha)<0$, then $M(\mu)\subset \Delta$.
        \item If $\phi'_\Delta(\alpha)>0$, then $M(\mu)\subset T\setminus \Delta$.
        \item If $\phi'_\Delta(\alpha)=0$, then $\alpha\in M(\mu)$.
    \end{enumerate}
\end{proposition}

\begin{remark}
    Given a convex function $f:\R\to\R$, the subdifferential of $f$ at $\alpha\in \R$ is the compact interval $[f_-'(\alpha), f_+'(\alpha)]$.
    Therefore, if $f_-'(\alpha)=0$ or $f_+'(\alpha)=0$, then $\alpha$ is a global minimizer of $f$.
    \Cref{prop:sign-derivative} shows remarkably that a similar fact holds in a metric tree:
    if the directional derivative is zero along a single direction, then $\alpha$ is a minimizer of $\phi$. 
\end{remark}

As a consequence, we obtain the following first-order optimality conditions. 

\begin{proposition}
    \label{prop:optimality}
    Let $\alpha \in \dom \phi$.
    \begin{enumerate}
        \item $\alpha \in M(\mu)$ if and only if for every $\Delta\in \dir_\alpha$, $\phi'_\Delta(\alpha)\ge 0$.
        \item If for every $\Delta\in \dir_\alpha$, $\phi'_\Delta(\alpha)>0$, then $\alpha$ is the unique minimizer of $\phi$, i.e., $M(\mu)=\{\alpha\}$.

    \end{enumerate}
\end{proposition}

\begin{remark}
    \label{rem:descent}
    Any $\Delta\in \dir_\alpha$ with $\phi'_{\Delta}(\alpha)<0$ is called a \textit{descent direction at $\alpha$}.
    By \Cref{prop:sign-derivative}, there is at most one descent direction at $\alpha$.
    By \Cref{prop:optimality} we obtain the following alternative, which is well-known in convex optimization on vector spaces:
    either there exists a descent direction at $\alpha$, or $\alpha\in M(\mu)$. 
\end{remark}

\begin{example}
    Assume that $\alpha$ has exactly two neighboring vertices, \ie $\dir_\alpha = \set{\Delta_1,\Delta_2}$.
    For Fr\'echet medians, by the first item of \Cref{prop:optimality},
    $\alpha$ is in $M_1(\mu)$ if and only if both 
    $1 - 2\mu(\Delta_{1}) \ge 0$ and $1 - 2\mu(\Delta_{2}) \ge 0$,
    which in turn is equivalent to
        $\mu(\Delta_{1}\cup\set{\alpha}) \ge \frac 12 \ \text{ and } \ \mu(\Delta_{2}\cup\set{\alpha}) \ge\frac 12$.
    This last optimality condition is reminiscent of the classical characterization of a median on $\R$ as any $m\in \R$ that satisfies
    both $\mu\left((-\infty,m]\right)\ge 1/2$ and $\mu\left([m,+\infty)\right)\ge 1/2$.
\end{example}

Lastly, the following localization property involves the convex hull of the support of $\mu$.
The convex hull of a subset $A$, denoted by $\conv(A)$, is the smallest geodesically convex subset of $T$ that contains $A$.
On the real line, it is well-known that if $A$ is closed then $\conv(A)$ is also closed (remarkably, this implication does not hold in higher dimensions).
We show that this property remains true on a metric tree.

\begin{proposition}
    \label{prop:localization-support}
    \begin{enumerate}
        \item Let $A$ be a closed subset of $T$. Then $\conv(A)$ is closed.
        \item The subset $\conv(\supp \mu)$ is closed and $M(\mu)\subset \conv(\supp \mu)$.
        \item Assume that $\supp \mu\subset G$ where $G$ is a convex subset of $T$. Then $M(\mu)\subset G$.
    \end{enumerate}
\end{proposition}

\section{Stickiness}
\label{sec:stickiness}

\subsection{Stickiness via directional derivatives}

A point $\alpha\in T$ can be classified according to the signs of directional derivatives at $\alpha$.
Such a classification gives rise to the concept of stickiness, which was introduced in \cite{hotz2013sticky} and further explored in 
\cite{huckemann2015sticky,bhat2017omnibus,lammers2023sticky}.

\begin{definition}
    \label{def:sticky}
    Let $\alpha\in \dom \phi$ with $\dir_\alpha = \set{\Delta_1,\ldots,\Delta_m}$.
    Depending on which of the following disjoint and exhaustive conditions is satisfied, we say that $\alpha$ is:
    \begin{itemize}
        \item \textit{sticky} if $\phi_{\Delta_i}'(\alpha) > 0$ for every $i\in [m]$, 
        \item \textit{partly sticky} if $\phi_{\Delta_i}'(\alpha) \ge 0$ for every $i\in [m]$, and $\phi_{\Delta_i}'(\alpha)=0$ for at least one $i\in [m]$, 
        \item \textit{nonsticky} if there exists some $i\in [m]$ such that $\phi_{\Delta_i}'(\alpha) < 0$.
    \end{itemize}
\end{definition}

\begin{remark}
    \begin{enumerate}
        \item By the first item of \Cref{prop:optimality}, $\alpha$ is in $M(\mu)$ if and only if it is either sticky or partly sticky.
        Put differently, $\alpha\notin M(\mu)$ if and only if $\alpha$ is nonsticky.
        \item In view of \Cref{rem:descent}, if $\alpha$ is nonsticky, there is a unique direction along which the derivative is negative.
        \item By the second item of \Cref{prop:optimality}, if $\phi$ has more than one minimizer, then all its minimizers are partly sticky.
    \end{enumerate}
\end{remark}

\begin{remark}
    Originally Hotz et al.\ \cite{hotz2013sticky} defined stickiness in the setting of the Fr\'echet mean ($\ell:z\mapsto z^2$) on an open book. 
    In \cite[Definition 2.10]{hotz2013sticky}, stickiness is defined according to the sign of the quantity 
    $\int_{\Delta} d(\alpha,x) \diff\mu(x) - \int_{T\setminus \Delta} d(\alpha,x) \diff\mu(x)$,
    where $\Delta$ is a direction in $\dir_\alpha$. With our notation, it can be rewritten as $-\frac 12 \phi_{\Delta}'(\alpha)$.
\end{remark}

According to the stickiness regime of the point $\alpha$, we obtain different localization constraints for $M(\mu)$, 
as seen in the following theorem.

\begin{theorem}
    \label{thm:localization-sticky}
    Let $\alpha\in \dom \phi$ with $\dir_\alpha = \set{\Delta_1,\ldots,\Delta_m}$.
    \begin{enumerate}
        \item If $\alpha$ is sticky, then $M(\mu)=\set{\alpha}$.
        \item If $\alpha$ is partly sticky, we let $I_0 = \set{i\in[m]: \phi_{\Delta_i}'(\alpha)=0}$.
        \begin{enumerate}
            \item The following inclusions hold: $\set{\alpha} \subset M(\mu) \subset \set{\alpha} \cup \bigcup_{i\in I_0} \Delta_{i}$.
            \item Assume that $\mu\neq \delta_\alpha$. Then $I_0$ has either one or two elements. 
            \item If $|I_0|=2$, then $\mu\left(\set{\alpha}\cup \bigcup_{i\in I_0} \Delta_{i}\right) = 1$.
        \end{enumerate}
        \item If $\alpha$ is nonsticky with $\phi_{\Delta_i}'(\alpha) < 0$ for some $i\in \{1,\ldots,m\}$, then $M(\mu)\subset \Delta_{i}$.
    \end{enumerate}
\end{theorem}

Motivated by the result on the number of vanishing derivatives in the partly sticky case, we introduce the following terminology.

\begin{definition}
    \label{def:partly-sticky-sided}
    Let $\alpha\in \dom \phi$ be partly sticky and let us borrow the notations of \Cref{thm:localization-sticky}.
    When $I_0$ is a singleton, we say that $\alpha$ is \textit{one-sided partly sticky}.
    When $I_0$ has two elements, we say that $\alpha$ is \textit{two-sided partly sticky}.
\end{definition}

The following proposition relates the stickiness regime of $\alpha$ to the geometry of the tree.
The \textit{degree of a vertex} refers to the number of edges incident to it.
We say that a vertex is \textit{branching} if its degree is at least $3$. 

\begin{proposition}
    \label{prop:sticky-geometry}
    Let $\alpha\in \operatorname{int}(F)$, write
    $\dir_\alpha=\set{\Delta_1,\ldots,\Delta_m}$, and assume that
    $\ell$ is differentiable at $d(\alpha,x)$ for $\mu$-almost every
    $x\in T\setminus\set{\alpha}$, that is 
    \begin{equation}
        \label{eq:differentiability}
        \mu\big(
        \set[\big]{x\in T\setminus\set{\alpha}: \ell_-'(d(\alpha,x)) \neq \ell_+'(d(\alpha,x))}
        \big) = 0.
    \end{equation}
    
    Then:
    \begin{enumerate}
        \item 
        $\sum_{i=1}^m \phi'_{\Delta_i}(\alpha)
        =
        m \ell_+'(0)\mu(\set{\alpha})
        +
        (m-2) \int_{T\setminus\set{\alpha}} \ell_+'(d(\alpha,x)) \diff \mu(x)$.
        
        \item If $\alpha$ is sticky and $\ell_+'(0)\mu(\set{\alpha})=0$,
        then $\alpha$ is a branching vertex.
        \item If $m=2$ and $\alpha\in M(\mu)$, then $\alpha$ is
        two-sided partly sticky if and only if
        $\ell_+'(0)\mu(\set{\alpha})=0$.
    \end{enumerate}
\end{proposition}

\begin{remark}
    We require $\alpha \in \operatorname{int}(F)$ so that $F$ intersects each direction at $\alpha$ and the closed form \eqref{eq:derivative} is available.
    As seen in \Cref{lemma:interior-F}, this topological condition admits a more explicit formulation.
    In the case where $\alpha$ is a leaf, $\alpha\in \operatorname{int}(F) \iff \alpha \in F$.
    Otherwise, it is equivalent to $\E{\ell_+'(d(\alpha, X) + \rho)} < +\infty$ for some $\rho \in (0, +\infty)$.
\end{remark}

\begin{remark}
    By \Cref{prop:sticky-geometry}, if $\mu$ is atomless, 
    any sticky Fr\'echet median of $\mu$ is located at a branching vertex.
    Similarly, if $\mu$ has a finite moment of order $p-1$ (where $p>1$),
    a sticky Fr\'echet $p$-mean of $\mu$ must be at a branching vertex.
\end{remark}

\subsection{Some examples}
\label{sec:examples-stickiness}

Let us illustrate the stickiness phenomenon in some concrete situations.
We first elucidate the case of a degenerate probability measure.

\begin{example}
    \label{ex:dirac}
    Let $c\in T$ and $\mu = \delta_c$.
    Then $M(\mu)=\set{c}$, and in every direction $\Delta\in \dir_c$ we have $\phi'_\Delta(c)=\ell_+'(0)$.
    Thus $c$ is sticky if and only if $\ell_+'(0)>0$.
    When $\ell_+'(0)=0$, $I_0$ has the same cardinality as $\dir_c$, which can be arbitrarily large.
\end{example}

Next, we turn to the case where $T$ is geometrically identical to an interval of $\R$.

\begin{proposition}
    \label{prop:stickiness-R}
    Assume that $T$ is isometric to a closed interval of positive length,
    and that $F=T$ (\eg \ref{assum:bounded} or \ref{assum:growth} is in force). 
    Suppose additionally that either
    \begin{enumerate}[label=(\alph*)]
        \item $\mu$ is atomless (\eg $\mu$ has a density \wrt the Lebesgue measure $\lambda$); or
        \item $\ell$ is differentiable on $(0,+\infty)$ and $\ell_+'(0)=0$.
    \end{enumerate}
    Let $\alpha \in M(\mu)$. Then:
    \begin{enumerate}
        \item For every $\Delta \in \dir_\alpha$, $\phi'_\Delta(\alpha) = 0$.
        \item If $\alpha$ is a leaf, then $\mu = \delta_\alpha$ and $\alpha$ is one-sided partly sticky.
        \item If $\alpha$ is not a leaf, then $\alpha$ is two-sided partly sticky.
    \end{enumerate}
\end{proposition}

In particular, when $T = \R$, 
any Fr\'echet median ($\ell:z\mapsto z$) of an atomless $\mu$ is two-sided partly sticky. 
The same conclusion holds for the Fr\'echet mean ($\ell:z\mapsto z^2$) of any $\mu$ with a finite first moment.
Consequently, $T$ being isometric to an interval prevents stickiness in a wide range of circumstances.
Still, one cannot entirely rule out stickiness in this setting.
As \Cref{ex:dirac} reveals, the Fr\'echet median of a Dirac mass $\delta_c$ is sticky.
The reader may verify that for $\mu = \frac 34 \delta_{\alpha_1} + \frac 14 \delta_{\alpha_2}$, the point $\alpha_1$ 
is a sticky Fr\'echet median of $\mu$, whatever the geometry of $T$.

Given any metric tree $T$ with a branching vertex, we may find a distribution $\mu$ on $T$ whose Fr\'echet $\ell$-mean is sticky,
irrespective of the loss $\ell$. 

\begin{example}
    Let $c$ be a branching vertex in $T$, 
    with $\dir_c=\set{\Delta_1,\ldots,\Delta_m}$ (where $m\ge 3$), 
    and let $\rho>0$ denote the distance from $c$ to its closest neighboring vertex.
    For each $i \in [m]$, we let $\alpha_i \in \Delta_i$ be such that $d(c, \alpha_i) = \rho$, 
    and we define $\mu=\frac 1m \sum_{i=1}^m \delta_{\alpha_i}$.
    Then, for every $i\in[m]$,
    $$
    \phi'_{\Delta_i}(c)
    =
    \frac 1m \bigl((m-1)\ell_+'(\rho)-\ell_-'(\rho)\bigr)
    \ge
    \frac{m-2}{m}\ell_+'(\rho)>0.
    $$
    Thus $c$ is sticky.
\end{example}

Next, we provide an example of a continuous $\mu$ whose Fr\'echet $\ell$-mean is sticky, whatever the loss $\ell$.

\begin{example}
    Let $m\ge 3$ and suppose that $T$ is an $m$-spider with center vertex $c$ and leaves $v_1,\ldots,v_m$.
    Suppose that $d(c, v_i) = 2$ and define $\alpha_i$ as the midpoint of $c$ and $v_i$ for each $i\in [m]$.
    Let $\mu$ denote the uniform measure on the subset $\bigcup_{i=1}^m [c,\alpha_i]$:
    $\mu$ has density $x\mapsto \frac 1m \sum_{i=1}^{m} \indic{[c,\alpha_i]}(x)$ with respect to the Lebesgue measure $\lambda$ on $T$.
    Then 
    \begin{align*}
        \phi_{v_1}'(c) = \frac{1}{m} \left(-\int_0^1 \ell_-'(z) \diff z  + \sum_{i=2}^{m}\int_0^1 \ell_+'(z) \diff z \right)
        = \frac{m-2}{m} \int_0^1 \ell_+'(z) \diff z > 0,
    \end{align*}
    where the second equality follows from the last item of \Cref{lemma:integral-rep}.
    By symmetry, all the directional derivatives $\phi_{v_i}'(c)$ are equal, hence positive, and $c$ is sticky.
\end{example}

\subsection{Stickiness via robustness}

The earlier definition of stickiness involves the signs of the directional derivatives at $\alpha$.
It is therefore stated in terms of the landscape of the objective function $\phi$ around $\alpha$.
The following theorem provides an equivalent formulation of stickiness:
$\alpha$ is sticky if and only if the equality $M(\nu) = \set{\alpha}$ holds for every measure $\nu$ that is sufficiently close to $\mu$.
The notion of stickiness thus has an interpretation in terms of robustness.

We quantify the closeness between two probability measures $\nu_1,\nu_2$ on $T$ using the total variation distance defined as 
$\tv(\nu_1,\nu_2) = \sup_{B\in \mathcal B(T)} |\nu_1(B)-\nu_2(B)|$.
We denote by $\mathcal P_1$ the set of probability measures on $T$ that have a finite first moment.
For measures in $\mathcal P_1$ we may also employ the $1$-Wasserstein distance
$W_1(\nu_1,\nu_2) = \sup\set{\int_T f(x) \diff\nu_1(x)-\int_T f(x) \diff\nu_2(x) : f:T\to \R \text{ is } 1\text{-Lipschitz}}$.
Total variation is stronger than $1$-Wasserstein when $T$ is bounded, in the sense that $W_1(\nu_1,\nu_2)\le \diam(T) \tv(\nu_1,\nu_2)$
\cite[Theorem~6.15]{villani2009optimal}.
Note that closeness in $W_1$ need not imply closeness in total variation.

\begin{theorem}[Stickiness and robustness to perturbations]
    \label{thm:robust-sticky}
    Let $\alpha\in \dom \phi$.
    \begin{enumerate}
        \item  Assume that $T$ is bounded or $\ell_+'$ is bounded.
        Then $\alpha$ is sticky if and only if 
        there exists $\varepsilon > 0$ such that for every probability measure $\nu$ satisfying $\tv(\nu, \mu)\le \varepsilon$, we have $M(\nu) = \set \alpha$.
        \item Assume that $\mu\in \mathcal P_1$, 
        that $\ell$ is differentiable on $I$ (where $I = [0,\diam(T)]$ if $T$ is bounded and $I=[0,\infty)$ otherwise), 
        that $\ell'$ is Lipschitz continuous on $I$
        and $\ell'(0)=0$.
        Then $\alpha$ is sticky if and only if 
        there exists $\varepsilon > 0$ such that for every $\nu\in \mathcal P_1$ satisfying 
        $W_1(\nu, \mu)\le \varepsilon$, we have $M(\nu) = \set \alpha$.
        
    \end{enumerate}
\end{theorem}

\begin{remark}
    For the total variation statement in \Cref{thm:robust-sticky}, it is not possible to replace the assumption that $T$ is bounded 
    with the weaker assumption that $\supp \mu$ is bounded.
    Indeed, consider the setting where $T$ is the open book composed of $3$ rays, with center $c$ and $\dir_c = \set{\Delta_1, \Delta_2, \Delta_3}$,
    define $v_i\in \Delta_i$ such that $d(v_i, c) = 1$ for each $i\in \set{1,2,3}$ and let $\mu = \frac 13 \sum_{i=1}^{3}\delta_{v_i}$.
    The loss is $\ell:z\mapsto z^2$.
    By the symmetry argument developed in \Cref{sec:examples-stickiness}, $c$ is sticky.
    Next, we demonstrate that the TV robustness statement does not hold in this case.
    Fix an arbitrary $\varepsilon \in (0,1]$. We consider $R$ such that $R>\frac{1-\varepsilon}{3\varepsilon}$,
    $w\in \Delta_1$ such that $d(w,c) = R$ and we define the perturbed measure $\nu= (1-\varepsilon)\mu + \varepsilon \delta_w$, so that $\tv(\nu, \mu)\le \varepsilon$.
    The objective function $\tphi$ corresponding to $\nu$ satisfies  $\tphi_{\Delta_1}'(c) = (1-\varepsilon)\phi_{\Delta_1}'(c) + \varepsilon (-2d(c,w)) 
    = \frac 23(1-\varepsilon) -2 \varepsilon R < 0$.
    This directional derivative being negative implies that $c\notin M(\nu)$, 
    hence the $\tv$ robustness statement in \Cref{thm:robust-sticky} does not hold when only $\supp \mu$ is bounded, while $T$ is unbounded.
    However, this specific example fulfills the conditions in the second item of \Cref{thm:robust-sticky}, hence $W_1$ robustness holds.
\end{remark}

\begin{remark}
    For the specific loss $\ell:z\mapsto z^2$, 
    connections between stickiness and robustness under perturbations have already appeared in the literature.
    The results stated in \cite[Section~7]{huckemann2015sticky}, \cite{bhat2017omnibus} and \cite[Section~5]{lammers2023types}
    cover stratified spaces.
    The general setting of a CAT($\kappa$) space was studied by Lammers et al.\ \cite{lammers2023sticky}.
    When the space is a metric tree without leaves, Theorem 4.22 in \cite{lammers2023sticky} implies our \Cref{thm:robust-sticky},
    albeit only in the special case of the squared loss.
\end{remark}

\section{Estimation of Fr\'echet \texorpdfstring{$\ell$}{ℓ}-means and statistical results}
\label{sec:stats}

\subsection{Estimation setting}
\label{sec:estimation-trees}

\begin{definition}
    Let $X,X_1,X_2,\ldots$ be a sequence of \iid $T$-valued random elements defined on some probability space $(\Omega,\mathcal A, \mathbb P)$, 
    each with distribution $\mu$.
    For each $n\ge 1$, we define the \textit{empirical measure} $\hmu_n = \tfrac 1n \sum_{k=1}^n \delta_{X_k}$ and
    the \textit{empirical objective function} $\hphi_n:\alpha \mapsto \frac 1n \sum_{k=1}^n [\ell(d(\alpha, X_k)) - \ell(d(o, X_k))]$.
    Minimizers of $\hphi_n$ (\ie elements of $M(\hmu_n)$) are called \textit{empirical (or sample) Fr\'echet $\ell$-means}.
    Given $\alpha\in T$ and $\Delta \in \dir_\alpha$, the directional derivative of $\hphi_n$ at $\alpha$ along $\Delta$
    will be denoted by $\hphi_\Delta'(\alpha)$ and referred to as an \textit{empirical directional derivative}
    (the integer $n$ is clear from the context and is therefore omitted from the notation).
\end{definition}

The quantities $X_n$, $\hphi_n$ and $\hmu_n$ are random. When needed, the dependence on $\omega\in \Omega$
will be made explicit by adding a superscript, \eg $X_n^{\omega},\hphi_n^{\omega},\hmu_n^{\omega}$.

We assume throughout that the probability space $(\Omega,\mathcal A, \mathbb P)$ is complete. 
By \cite[Proposition~5.3.9, Proposition 5.3.13 (iii)]{molchanov2017theory}, 
$M(\hmu_n)$ is a \textit{measurable closed subset} of $T$, \ie for each compact subset $K\subset T$, 
the subset $\{\omega \in \Omega: M(\hmu_n^{\omega})\cap K\neq\emptyset\}$ is in $\mathcal A$. 
Moreover, let us note that for each Borel subset $B\subset T$, 
$\set{M(\hmu_n)\subset B}$ coincides with $\set{M(\hmu_n)\cap (T\setminus B) = \emptyset}$,
the latter being in $\mathcal A$ according to \cite[Theorem~1.3.3(iii)]{molchanov2017theory}.

By \Cref{lemma:Borel-selection}, there exists a measurable selection of $M(\hmu_n)$, 
\ie a minimizer $\halpha_n$ of $\hphi_n$ that is measurable.
More precisely, $\halpha_n$ is a Borel function of $(X_1,\ldots,X_n)$.
This result is most useful when $M(\hmu_n)$ is not a singleton, which may happen when $\ell$ is not strictly convex.

\subsection{A law of large numbers for sets}

In terms of sets, estimation is successful if the random set $M(\hmu_n)$ gets closer in some sense to the true set $M(\mu)$ as $n\to \infty$. 
Two modes of convergence are usually considered for the sequence $(M(\hmu_n))_{n\ge 1}$.

\begin{definition}
    \label{def:consistency-sets}
    \begin{enumerate}
        \item $(M(\hmu_n))_{n\ge 1}$ is \textit{strongly consistent in outer limit} \cite{schotz2022strong} 
        (alternatively, \textit{in the sense of Ziezold} \cite{huckemann2011intrinsic}) if
        $$\P[\big]{\bigcap_{n\ge 1} \overline{\bigcup_{p\ge n} M(\hmu_p)} \subset M(\mu)} = 1.$$

        \item $(M(\hmu_n))_{n\ge 1}$ is \textit{strongly consistent in one-sided Hausdorff distance} \cite{schotz2022strong} 
        (alternatively, \textit{in the sense of Bhattacharya--Patrangenaru} \cite{huckemann2011intrinsic}) if
        $$\P[\big]{\sup_{\alpha\in M(\hmu_n)} \inf_{\beta \in M(\mu)} d(\alpha,\beta) \xrightarrow[n\to\infty]{} 0} = 1.$$
    \end{enumerate}
\end{definition}

Each of these statements is regarded as a set-valued strong law of large numbers.
A caveat about \Cref{def:consistency-sets} is that these notions of closeness between $M(\hmu_n)$ and $M(\mu)$ are only one-sided:
there might exist some $\alphas \in M(\mu)$ such that the distance of $\alphas$ to $M(\hmu_n)$ remains bounded away from $0$ with positive probability.
However, in the case of a unique Fr\'echet $\ell$-mean (\ie $M(\mu)=\set{\alphas}$), 
strong consistency in one-sided Hausdorff distance implies strong consistency in the usual sense:
$d(\halpha_n, \alphas)\to 0$ almost surely
for any sequence of measurable selections $(\halpha_n)_{n\ge 1}$.

\begin{theorem}[Law of large numbers]
    \label{thm:lln-sets}
    $(M(\hmu_n))_{n\ge 1}$ is strongly consistent in both senses of \Cref{def:consistency-sets}.
\end{theorem}

\begin{remark}
    Remarkably, \Cref{thm:lln-sets} requires no assumption beyond our standing integrability condition $o\in F$,
    \ie $\int_T \ell_+'(d(o,x)) \diff \mu(x) < \infty$.
    Thus, our result is more general than what follows from the results 
    \cite[Corollary~4.3 and Corollary~4.4]{schotz2022strong} or \cite[Theorem 3.1]{aveni2026uniform}, 
    since no additional growth assumption on $\ell_+'$ or $\ell$ is required.
\end{remark}

\subsection{Empirical stickiness}

We leverage the concept of stickiness introduced in \Cref{sec:stickiness} to describe in more detail how $(M(\hmu_n))_{n\ge 1}$ converges to $M(\mu)$.
When the directional derivative $\phi_{\Delta_i}'(\alpha)$ is nonzero, \Cref{prop:sign-derivative} helps to locate the minimizers of $\phi$.
As $n$ grows, it is expected that the empirical counterpart $\hphi_{\Delta_i}'(\alpha)$ approaches $\phi_{\Delta_i}'(\alpha)$ with high probability,
which entails that it becomes nonzero and has the same sign as $\phi_{\Delta_i}'(\alpha)$.
Since there are only finitely many directions at $\alpha$, 
it is then possible to obtain identical localization constraints on the empirical Fr\'echet $\ell$-means.
The following results make this intuition precise. 

We start with a sticky law of large numbers, which is an asymptotic result.
When $\alpha$ is sticky, with probability $1$ the sets of empirical minimizers $M(\hmu_n)$ eventually collapse to $\set{\alpha}$. %
This justifies the use of the adjective \say{sticky}.

\begin{theorem}[Sticky law of large numbers]
    \label{thm:sticky-lln}
    Let $\alpha\in \operatorname{int}(F)$ and write $\dir_\alpha = \set{\Delta_1,\ldots, \Delta_m}$.
    \begin{enumerate}
        \item If $\alpha$ is sticky, then with probability $1$, $M(\hmu_n)=\set{\alpha}$ for all large enough $n$.

        \item If $\alpha$ is partly sticky, let $I_0$ be as in \Cref{thm:localization-sticky}.
        \begin{enumerate}
            \item If $I_0 = \set{i}$, then with probability $1$, $M(\hmu_n)\subset \set{\alpha} \cup \Delta_{i}$ for all large enough $n$.
            
            \item If $I_0 = \set{i,i'}$, then with probability $1$, 
            $M(\hmu_n)\subset \set{\alpha} \cup \Delta_{i} \cup \Delta_{i'}$ for all $n\ge 1$.
        \end{enumerate}

        \item If $\alpha$ is nonsticky with $\phi_{\Delta_i}'(\alpha) < 0$ for some $i\in\{1,\ldots,m\}$, then with probability $1$, $M(\hmu_n)\subset \Delta_{i}$ for all large enough $n$.
    \end{enumerate}
\end{theorem}

\begin{remark}
    We require $\alpha \in \operatorname{int}(F)$ so that $F$ intersects each direction at $\alpha$ and the closed form \eqref{eq:derivative} is available.
    As seen in \Cref{lemma:interior-F}, this topological condition admits a more explicit formulation.
    In the case where $\alpha$ is a leaf, $\alpha\in \operatorname{int}(F) \iff \alpha \in F$.
    Otherwise, it is equivalent to $\E{\ell_+'(d(\alpha, X) + \rho)} < +\infty$ for some $\rho \in (0, +\infty)$.
    For the specific item 2.b), the proof reveals that $\alpha \in \operatorname{int}(F)$ may be relaxed to $\alpha \in \dom \phi$.
\end{remark}

For convenience, let us say that $\alpha$ is \textit{directionally sticky} when it is sticky in the sense of \Cref{def:sticky}
and that $\alpha$ is \textit{sample sticky} when $\P{\exists N\ge 1, \forall n\ge N, \alpha \in M(\hmu_n)}=1$.
We provide a converse to the first item of \Cref{thm:sticky-lln}.

\begin{proposition}
    \label{prop:converse-empirical-stickiness}
    Let $\alpha\in F$ be sample sticky.
    Then $\alpha$ is directionally sticky or $\mu = \delta_\alpha$.
\end{proposition}

We next turn to non-asymptotic results.
A random variable $Z$ is said to be sub-Gaussian with variance proxy $\sigma^2$ if $\E{e^{(Z-\E{Z})t}} \le e^{\frac{\sigma^2t^2}2}$ for every $t\in \R$. 
In the case where the random variable $\ell_+'(d(\alpha,X))$ is sub-Gaussian, stickiness has a non-asymptotic formulation.

\begin{theorem}[Non-asymptotic empirical stickiness]
    \label{thm:sticky-nonasymp-subgaussian}
    Let $\alpha\in \operatorname{int}(F)$ with $\dir_\alpha = \set{\Delta_1,\ldots, \Delta_m}$ and let $n\ge 1$ be fixed.
    Assume that $\ell_+'(d(\alpha,X))$ is sub-Gaussian.
    Then for each $i\in [m]$,
    $\indic{T\setminus \Delta_{i}}(X) \ell_+'(d(\alpha,X)) 
        - \indic{\Delta_{i}}(X) \ell_-'(d(\alpha,X))$ is sub-Gaussian and we 
    denote by $\sigma_i^2$ a variance proxy.
    \begin{enumerate}
        \item If $\alpha$ is sticky, then 
        $$\P{M(\hmu_n)=\set{\alpha}} \ge 1 - \sum_{i=1}^m \exp\bigg(-\frac{n\phi_{\Delta_i}'(\alpha)^2}{2\sigma_i^2}\bigg).$$

        \item If $\alpha$ is one-sided partly sticky with $I_0=\{i_0\}$, then
            $$ \P[\big]{M(\hmu_n)\subset \set{\alpha} \cup \Delta_{i_0}} 
            \ge 1 - \sum_{i\neq i_0} \exp\bigg(-\frac{n\phi_{\Delta_i}'(\alpha)^2}{2\sigma_i^2}\bigg).$$
        
        \item If $\alpha$ is nonsticky with $\phi_{\Delta_i}'(\alpha) < 0$ for some $i\in\{1,\ldots,m\}$, then 
        $$\P{M(\hmu_n)\subset \Delta_{i}} \ge 1 -  \exp\bigg(-\frac{n\phi_{\Delta_i}'(\alpha)^2}{2\sigma_i^2}\bigg).$$ 
    \end{enumerate}
\end{theorem}

\begin{remark}
    Let us be more precise in the bounded case where \ref{assum:bounded} holds. Let $L = \ell_+'\big(\sup_{x\in \supp \mu} d(\alpha,x)\big)$,
    with the convention that $L = \lim_{z\to\infty}\ell_+'(z)$ when $\supp \mu$ is unbounded. 
    In that case, the random variable $\indic{T\setminus \Delta_{i}}(X) \ell_+'(d(\alpha,X)) 
        - \indic{\Delta_{i}}(X) \ell_-'(d(\alpha,X))$ takes values in $[-L,L]$, hence it is sub-Gaussian with variance proxy $\sigma_i^2 = L^2$.
\end{remark}

\subsection{Central limit theorems in the partly sticky case}
\label{sec:CLT}

In this section we assume that $M(\mu)=\set{\alphas}$, \ie there is a unique Fr\'echet $\ell$-mean $\alphas$. 
We are thus in the classical setting of parameter estimation. 
However, the empirical set $M(\hmu_n)$ may not be a singleton; 
we consider therefore an arbitrary sequence of measurable selections $(\halpha_n)_{n\ge 1}$, as explained in \Cref{sec:estimation-trees}.
Recall that $\alphas$ must be either sticky or partly sticky, since it is in $M(\mu)$. 

In the rest of this section, we assume that $\alphas$ is partly sticky
and $\mu \neq \delta_{\alphas}$, so that by \Cref{thm:localization-sticky} the number of vanishing directional derivatives at $\alphas$ is either $1$ or $2$.
We write $\rhos$ for the distance from $\alphas$ to its closest neighboring vertex.
Let us introduce two integrability assumptions, 
followed by some notation and technical results that will prove useful in stating and establishing central limit theorems.

\begin{assumption}
    \label{assum:clt-1}
    $\E{\ell_+'(d(\alphas, X) + \rho_0)} < \infty$ for some $\rho_0\in (0, \rhos)$.
\end{assumption}

\begin{assumption}
    \label{assum:clt-2}
    $\E{\ell_+'(d(\alphas, X) + \rho_0)^2} < \infty$ for some $\rho_0\in (0, \rhos)$.
\end{assumption}

\begin{remark}
    Let us note that \ref{assum:clt-2} implies \ref{assum:clt-1}.
    Moreover, if \ref{assum:clt-1} holds, then the inclusion $B(\alphas, \rho_0) \subset F$ 
    follows from the triangle inequality and the monotonicity of $\ell_+'$, hence $\alphas \in \operatorname{int}(F)$.
\end{remark}

\begin{definition}
    \label{def:Gamma-S}
    In the one-sided case, we denote by $\Delta_+\in \dir_{\alphas}$ the direction such that $\phi_{\Delta_+}'(\alphas)=0$.
    In the two-sided case, we denote by $\Delta_-, \Delta_+\in \dir_{\alphas}$ the directions with vanishing derivative.
    Given $\rho_0\in (0, \rhos)$,
    we denote by $\gamma_+$ a unit-speed geodesic such that $\gamma_+(0)=\alphas$ and $\gamma_+((0,\rho_0]) \subset \Delta_+$.
    In the two-sided case, we additionally define a unit-speed geodesic
    $\gamma_-$ such that $\gamma_-(0)=\alphas$ and $\gamma_-((0,\rho_0]) \subset \Delta_-$.
    We let $C = [0,+\infty)$ in the one-sided case and $C=\R$ in the two-sided case.
    We introduce the \textit{signed radial coordinate mapping}
    $S:T\to \R$, $x \mapsto (\indic{\Delta_+}(x)-\indic{T\setminus \Delta_+}(x))d(\alphas, x)$,
    as well as $\hm_n = S(\halpha_n)$
    and the geodesic $\begin{aligned}[t]
        \Gamma \colon C\cap [-\rho_0,\rho_0] &\to T \\
        z &\mapsto \begin{cases}
            \gamma_+(z) &\text{ if } z\ge 0,
            \\\gamma_-(|z|) &\text{ if } z<0.
        \end{cases} 
    \end{aligned}$

\end{definition}

The following proposition relates $\Gamma$ and $S$.

\begin{proposition}
    \label{prop:Gamma-S}
    \begin{enumerate}
        \item The equality $d(\Gamma_z, x) = |z - S(x)|$ holds for every $z\in C\cap [-\rho_0,\rho_0]$ and every 
        $x\in T$ (resp., $x\in \Delta_- \cup \set{\alphas} \cup \Delta_+$) in the one-sided (resp., two-sided) case.
        \item $S\circ \Gamma = \Id$ and $\Gamma$ is unit-speed, \ie $d(\Gamma_z, \Gamma_{z'}) = |z-z'|$.
    \end{enumerate}
\end{proposition}

The coordinate mapping $S$ flattens the tree by mapping it to an interval of $\R$.
It follows from \Cref{prop:Gamma-S} that the restriction of $S$ to the subset $\Gamma(C\cap [-\rho_0,\rho_0])$ is an isometry.
This fact simplifies the analysis of $\halpha_n$, as it allows us to focus instead on $\hm_n = S(\halpha_n)$, 
which is a simpler object as it lives on the real line.
In order to investigate the properties of $\hm_n$, we introduce some additional functions.

\begin{definition}
    \label{def:Psi}
    Let $\psi:\R \to \R$, $z\mapsto \ell(|z|)$, and let us define the
    \textit{empirical local objective function}
    $\hPsi_n: \R \to \R$, $z\mapsto \frac 1n \sum_{k=1}^{n} [\psi(z-S(X_k)) - \psi(S(X_k))]$.
    Under \Cref{assum:clt-1}, which supplies the value of $\rho_0$ used in \Cref{def:Gamma-S},
    we introduce the 
    \textit{local objective function}
    $\Psi: C\cap [-\rho_0,\rho_0] \to \R$, $z \mapsto \E{\psi(z-S(X)) - \psi(S(X))}.$
\end{definition}

Next, we establish a connection between the local objective $\Psi$ and the original objective $\phi$,
as well as between the empirical counterparts $\hphi_n$ and $\hPsi_n$.
We also provide results on the right derivative of $\Psi$, which plays a major role in the central limit theorems to come.

\begin{proposition}
    \label{prop:Psi}
    \begin{enumerate}
        \item Assume \ref{assum:clt-1}. 
        $\Psi$ is well-defined, convex, with unique minimizer $0$ and 
        $\phi(\Gamma_z) - \phi(\alphas) = \Psi(z)$ for every $z\in C\cap [-\rho_0,\rho_0]$.
        \item Assume \ref{assum:clt-1}. 
        For each $z\in \operatorname{int}(C\cap [-\rho_0,\rho_0])$, $\Psi_+'(z) = \E{\psi_+'(z-S(X))}$ and $\Psi_-'(z) = \E{\psi_-'(z-S(X))}$.
        At $0$, $\Psi_+'(0)=\E{\psi_+'(-S(X))}=0$ and in the two-sided case $\Psi_-'(0)=\E{\psi_-'(-S(X))}=0$.
        \item Assume that $\E{\ell_+'(d(\alphas, X))^2} < +\infty$.
        The quantity $\sigma^2\coloneqq \E{\psi_+'(-S(X))^2}$ belongs to $(0,+\infty)$.
        \item $\P[\big]{\forall n\ge 1, \forall z\in C\cap [-\rho_0,\rho_0],\,\hphi_n(\Gamma_z) - \hphi_n(\alphas) = \hPsi_n(z)} = 1$.
        \item Assume that $\alphas \in \operatorname{int}(F)$.
        With probability $1$, there exists $N\ge 1$ such that for every $n\ge N$ the following statements hold:
        $M(\hmu_n) \subset \Gamma(C\cap [-\rho_0,\rho_0])$, $\hm_n\in C\cap [-\rho_0,\rho_0]$, $\Gamma_{\hm_n} = \halpha_n$ and 
        $\hm_n\in \argmin_{z\in C\cap [-\rho_0,\rho_0]} \hPsi_n(z)$.
    \end{enumerate}
\end{proposition}

\subsubsection{The two-sided partly sticky case}
\label{sec:two-sided-clt}

We start with the two-sided partly sticky case where $\phi_{\Delta_-}'(\alphas) = \phi_{\Delta_+}'(\alphas)=0$. 
As seen in Theorems~\ref{thm:localization-sticky} and \ref{thm:sticky-lln}, two derivatives being zero results in the following localization constraints:
$\supp \mu \subset \Delta_-\cup\set{\alphas}\cup\Delta_+$ 
and $\P{\forall n\ge 1, \,M(\hmu_n) \subset \Delta_-\cup\set{\alphas}\cup\Delta_+} = 1$.
The following proposition adds to the picture: asymptotically,
the whole set of empirical minimizers $M(\hmu_n)$
has an equal probability of lying in $\Delta_-$ or $\Delta_+$.

\begin{proposition}
    \label{prop:proba-localization-two-sided}
    Assume that $\E{\ell_+'(d(\alphas, X))^2} < +\infty$.
    The probabilities $\P{M(\hmu_n) \subset \Delta_+}$, $\P{\halpha_n \in \Delta_+}$,
    $\P{M(\hmu_n) \subset \Delta_-}$ and $\P{\halpha_n \in \Delta_-}$ 
    go to $1/2$ as $n\to \infty$.
\end{proposition}

The statement of the following central limit theorem makes use of the notations introduced in 
Definitions~\ref{def:Gamma-S} and \ref{def:Psi}.
In the two-sided case, $C\cap [-\rho_0,\rho_0] = [-\rho_0,\rho_0]$, moreover the local objective $\Psi$ is convex on this interval and $0$ is its unique minimizer,
hence $\Psi_+'$ is increasing and, for $z\in (-\rho_0,\rho_0)$, $\Psi_+'(z)=0\iff z=0$.

\begin{theorem}[Two-sided central limit theorem]
    \label{thm:two-sided-clt}
    Assume the following: 
    \begin{enumerate}[label=(\alph*)]
        \item \Cref{assum:clt-2} holds.
        \item $\Psi_+'$ satisfies
        $$\frac{\Psi_+'(z)}{z^{a_+}} \xrightarrow[z\downarrow 0]{} K_+  \quad \text{and} \quad  \frac{\Psi_+'(z)}{|z|^{a_-}} \xrightarrow[z\uparrow 0]{} -K_-,$$
        where $a_-, a_+,K_-,K_+$ are positive constants.
    \end{enumerate}
    Let $Z\sim N(0,1)$ and let $\sigma^2\in (0,\infty)$ be as in \Cref{prop:Psi}. 
    Then:
    \begin{enumerate}
        \item $\forall t> 0, \quad \P*{n^{1/(2a_+)} \hm_n \le t} \xrightarrow[n\to\infty]{} \P*{Z\le \frac{K_+ }{\sigma}t^{a_+}}$.
        \item $\forall t< 0, \quad \P*{n^{1/(2a_-)} \hm_n \le t} \xrightarrow[n\to\infty]{} \P*{Z\le -\frac{K_- }{\sigma}|t|^{a_-}}$.
        \item Assume additionally that $a_-=a_+=a$. Then $\big(n^{1/(2a)} \hm_n\big)_{n\ge 1}$ converges in distribution to
        $$\frac{\sgn(Z)+1}{2} \Big(\frac{\sigma}{K_+}|Z|\Big)^{1/a}  
        + \frac{\sgn(Z)-1}{2} \Big(\frac{\sigma}{K_-}|Z|\Big)^{1/a},$$
        and $\big(n^{1/(2a)} d(\halpha_n, \alphas)\big)_{n\ge 1}$ converges in distribution to
        $$\frac{1+\sgn(Z)}{2} \Big(\frac{\sigma}{K_+}|Z|\Big)^{1/a}  
        + \frac{1-\sgn(Z)}{2} \Big(\frac{\sigma}{K_-}|Z|\Big)^{1/a}.$$
    \end{enumerate}
\end{theorem}

\begin{remark}
    \label{rem:implies-A4}
    The boundedness assumption \ref{assum:bounded} clearly implies \ref{assum:clt-2}.
    If simultaneously \ref{assum:growth} holds and $\E{\ell_+'(d(o,X))^2} < \infty$,
    then \ref{assum:clt-2} is fulfilled, as seen in \Cref{lemma:assum-growth-second-moment}.
\end{remark}

\begin{remark}
    Condition (b) in \Cref{thm:two-sided-clt} is stated in terms of the local objective, which is defined through $\ell$ and the signed coordinate $S$.
    However, this condition admits an equivalent, more intrinsic formulation, in terms of 
    the local behavior of the original objective $\phi$ along the two branches through $\alphas$ corresponding to $\Delta_+$ and $\Delta_-$.
    More precisely, for $\varepsilon \in \set{+,-}$, let us define
    $\Phi^\varepsilon: [0,\rho_0] \to \R$, $t\mapsto \phi(\gamma_\varepsilon(t))$.
    As seen in \Cref{lemma:equiv-assum}, condition (b) is equivalent to the convergence $(\Phi^\varepsilon)_+'(t)/t^{a_\varepsilon} \to K_\varepsilon$ as $t\downarrow 0$.
\end{remark}

\begin{remark}
    When $a_- = a_+ = a \ge 1$, the limiting distribution for $n^{1/(2a)} \hm_n$ is known as a two-piece power normal distribution \cite{sulewski2021two}
    with location parameter $0$, left scale parameter $(\sigma/K_-)^{1/a}$, right scale parameter $(\sigma/K_+)^{1/a}$
    and shape parameter $a$.
\end{remark}

We specialize the preceding theorem to the regular setting,
where $\Psi_+'$ is differentiable at $0$ with positive derivative.
In that case, $\halpha_n$ is $\sqrt n$-consistent and the limiting distribution for $\sqrt n \hm_n$ is Gaussian.
The following corollary gives readily verifiable conditions on the loss and on the local behavior of $\mu$
that guarantee differentiability of $\Psi_+'$, as well as an explicit formula for the asymptotic variance.

\begin{corollary}[The regular case with a smooth loss]
    \label{corol:clt-regular}
    Assume the following:
    \begin{enumerate}[label=(\alph*)]
        \item $\E{\ell_+'(d(\alphas, X))^2} < +\infty$.
        \item $\ell$ is twice differentiable on $[0,+\infty)$
        and $\ell''$ is bounded, \ie $L\coloneqq \sup_{z\ge 0} \ell''(z)\in [0,\infty)$.
        \item In the case where $\ell'(0)>0$, 
        there exists $\varepsilon>0$ such that the restriction $\mu_{\mid B(\alphas, \varepsilon) \cap (\Delta_-\cup\set{\alphas}\cup\Delta_+)}$
        is absolutely continuous \wrt the Lebesgue measure $\lambda_{\mid B(\alphas, \varepsilon) \cap (\Delta_-\cup\set{\alphas}\cup\Delta_+)}$, with a density $f$ such that
        $$\lim_{z\downarrow 0} f(\gamma_+(z)) = \lim_{z\downarrow 0} f(\gamma_-(z)).$$
        We denote by $\eta$ this common limit and we assume that $\eta< +\infty$.
    \end{enumerate}
    Define $$K \coloneqq \begin{cases}
        \E{\ell''(d(\alphas, X))} & \text{ if } \ell'(0) = 0,\\
        \E{\ell''(d(\alphas, X))} + 2 \ell'(0) \eta & \text{ if } \ell'(0) > 0,
    \end{cases}$$
    and assume that $K>0$.

    Then $\Psi_+'(z) = Kz + o(z)$ as $z\to 0$.
    Consequently, \Cref{thm:two-sided-clt} applies with 
    $a_+=a_-=1$, $K_+=K_-=K$ and $\sigma^2 = \E{\ell'(d(\alphas, X))^2}$.
    In particular, with $Z\sim N(0,1)$, 
    $$\sqrt n \hm_n \xrightarrow[]{(d)} \frac \sigma K Z
    \quad \text{and}\quad  
    \sqrt n d(\halpha_n, \alphas) \xrightarrow[]{(d)} \frac \sigma K |Z|.
    $$
\end{corollary}

\begin{example}
    \begin{enumerate}
        \item For the Fr\'echet mean ($\ell:z\mapsto z^2$), 
        condition (a) of \Cref{corol:clt-regular} is equivalent to $\E{d^2(\alphas, X)}<+\infty$, \ie $X$ has a finite second moment.
        Condition (b) is satisfied, hence the corollary applies with 
        $K = 2$, $\sigma^2=4\E{d^2(\alphas, X)}$ and $(\sqrt n \hm_n)_{n\ge 1}$ converges in distribution to $N(0, \E{d^2(\alphas, X)})$.

        Consider the case where $T$ is the double ray and identify $T$ with $\R$.
        Then $\alphas$ is equal to $\E{X}$, $\halpha_n$ is the sample mean, $\E{d^2(\alphas, X)} = \V{X}$ and
        we recover the classical Lindeberg--L\'evy central limit theorem.

        \item For Fr\'echet medians ($\ell:z\mapsto z$), conditions (a) and (b) of \Cref{corol:clt-regular} are fulfilled.
        Under condition (c), we have $K = 2\eta$ and $\sigma^2=1$. Provided that $\eta > 0$, $(\sqrt n \hm_n)_{n\ge 1}$ converges in distribution to $N(0, 1/(4\eta^2))$.
        For convenience, we denote by $\F$ the CDF of $S(X)$.
        As seen in the proof of the corollary, condition (c) together with the requirement that $K>0$ may be replaced 
        by the assumption that $\F$ is differentiable at $0$, with $\F'(0)>0$.
        In the case of the double ray, we recover the classical condition for asymptotic normality of empirical medians 
        (see, \eg \cite[Corollary A p.~77]{serfling1980approx}).

        \item For the pseudo-Huber loss $\ell: z\mapsto 2c^2((1+z^2/c^2)^{1/2} -1)$,
        conditions (a) and (b) of \Cref{corol:clt-regular} are met, and $\ell'(0) = 0$.
        Therefore, the corollary applies with
        $$K = 2\E[\Big]{\Big(1+\frac{d^2(\alphas, X)}{c^2}\Big)^{-3/2}} \quad \text{and} \quad 
        \sigma^2 = 4\E*{\frac{d^2(\alphas, X)}{1+\frac{d^2(\alphas, X)}{c^2}}}.$$
    \end{enumerate}
\end{example}

\subsubsection{The one-sided partly sticky case}
Let us turn to the one-sided partly sticky case where only $\phi_{\Delta_+}'(\alphas)=0$. 
For each direction $\Delta \in \dir_\alphas \setd{\Delta_+}$,
the probability $\P{\hphi_{\Delta}'(\alphas)>0}$ goes to $1$ as $n$ grows.
As revealed in \Cref{lemma:clt}, $\P{\hphi_{\Delta_+}'(\alphas)>0}$
and $\P{\hphi_{\Delta_+}'(\alphas)\ge 0}$ both go to $1/2$.
As a straightforward consequence, the event $\set{M(\hmu_n) = \set{\alphas}}$
has probability converging to $1/2$.
In words, empirical Fr\'echet $\ell$-means exhibit a partial collapse, 
similar in spirit to the sticky case, but occurring with limiting probability $1/2$.

\begin{proposition}
    \label{prop:one-sided-collapse}
    Assume that $\alphas \in \operatorname{int}(F)$ and $\E{\ell_+'(d(\alphas, X))^2} < +\infty$.
    The probabilities $\P{M(\hmu_n) = \set{\alphas}}$, $\P{\halpha_n = \alphas}$,
    $\P{M(\hmu_n) \subset \Delta_+}$ and 
    $\P{\halpha_n \in \Delta_+}$ go to $1/2$ as $n\to \infty$.
\end{proposition}

Unlike the two-sided case, one-sided partial stickiness does not guarantee that the support of $\mu$ 
is a subset of $\set{\alphas}\cup \Delta_+$.
Although the probability $\P{\halpha_n \notin \set{\alphas}\cup \Delta_+}$ may be positive for each $n$, this is harmless for the central limit theorem below, 
because the event $\set{\halpha_n \in \Gamma([0,\rho_0])}$ has probability going to $1$.
We reuse the notations introduced in Definitions~\ref{def:Gamma-S} and \ref{def:Psi}.

\begin{theorem}[One-sided central limit theorem]
    \label{thm:one-sided-clt}
    Assume the following: 
    \begin{enumerate}[label=(\alph*)]
        \item \Cref{assum:clt-2} holds.
        \item $\Psi_+'$ satisfies
        $\frac{\Psi_+'(z)}{z^{a_+}} \xrightarrow[z\downarrow 0]{} K_+$,
        where $a_+, K_+$ are positive constants.
    \end{enumerate}
    Let $Z\sim N(0,1)$ and let $\sigma^2\in (0,\infty)$ be as in \Cref{prop:Psi}. 
    Then:
    \begin{enumerate}
        \item $\forall t> 0, \quad \P*{n^{1/(2a_+)} \hm_n \le t} \xrightarrow[n\to\infty]{} \P[\big]{Z\le \frac{K_+}{\sigma} t^{a_+}}$.
        \item $\forall t< 0, \quad \P*{n^{1/(2a_+)} \hm_n \le t} \xrightarrow[n\to\infty]{} 0$.
        \item $\big(n^{1/(2a_+)} \hm_n\big)_{n\ge 1}$ and $\big(n^{1/(2a_+)} d(\halpha_n, \alphas)\big)_{n\ge 1}$ converge weakly to $\big(\frac{\sigma}{K_+}\max(Z,0)\big)^{1/a_+}$.
    \end{enumerate}
\end{theorem}

The limiting distribution of $n^{1/(2a_+)} d(\halpha_n, \alphas)$ thus puts half of its mass at $0$.
This agrees with the convergence of $\P{\halpha_n = \alphas}$ to $1/2$ established in \Cref{prop:one-sided-collapse}.
Conditionally on the event $\set{\halpha_n \neq \alphas}$, this rescaled distance converges in distribution to $(\sigma |Z| / K_+)^{1/a_+}$.
Therefore, conditionally on non-collapse, $\halpha_n$ has fluctuations of order $n^{-1/(2a_+)}$ into the direction $\Delta_+$.

\section{Estimating the stickiness regime \texorpdfstring{of $\alphas$}{}}
\label{sec:estimate-stickiness}

We assume in this section that there is a unique Fr\'echet $\ell$-mean $\alphas$, \ie $M(\mu)=\set{\alphas}$.
We write $\dir_{\alphas} = \set{\Delta_1,\ldots,\Delta_m}$ and $I_0 = \set{i\in [m]: \phi_{\Delta_i}'(\alphas) = 0}$.
We assume additionally that $\mu \neq \delta_{\alphas}$, which entails by \Cref{thm:localization-sticky} that $|I_0| \in \set{0,1,2}$.
We denote by $\rhos$ the distance from $\alphas$ to its closest neighboring vertex.

Our goal is to determine from the data whether the unknown parameter $\alphas$ is sticky, one-sided partly sticky or two-sided partly sticky.
In other words, we wish to estimate the stickiness regime of $\alphas$, or equivalently, estimate the cardinality of $I_0$.

In what follows, we introduce two estimators of the stickiness regime, 
based respectively on pairwise collisions of empirical Fr\'echet $\ell$-means and on empirical difference quotients.
The consistency guarantees for both estimators rely crucially on the asymptotic theory developed in \Cref{sec:stats}, 
in particular on the laws of large numbers and central limit theorems established there.

\subsection{Via pairwise collisions}
\label{sec:collisions}

Recall from \Cref{sec:stats} that in the sticky case, empirical Fr\'echet $\ell$-means collapse to $\alphas$,
while in the one-sided partly sticky case they exhibit partial collapse, and 
there is no collapse in the two-sided partly sticky case.
This motivates the following approach.
We split the $n$-sample into $j_n$ disjoint blocks of equal size, 
and for each block (indexed by $j$) we compute a corresponding empirical Fr\'echet $\ell$-mean $\halpha_{n,j}$.
We then consider the proportion $\hp_n$ of pairwise collisions among the $\halpha_{n,j}$.
For large $n$, we expect in the sticky case that this proportion is close to $1$;
in the one-sided (resp., two-sided) case it should be close to $1/4$ (resp., $0$).
This rationale leads to the following definition.

\begin{definition}
    \label{def:estimator-collisions}
    For each $n \geq 2$, let $j_n \in [n]\setd{1}$ and let $B_{n,1},\ldots,B_{n,j_n}$ be disjoint deterministic subsets of $[n]$, each with cardinality $b_n = \floor{n/j_n}$.
    For each $j\in [j_n]$, let $\hmu_{n,j} = \frac{1}{b_n} \sum_{k\in B_{n,j}} \delta_{X_k}$ and 
    let $\halpha_{n,j} = f_{b_n}((X_k)_{k\in B_{n,j}})$ be the measurable selection
    from $M(\hmu_{n,j})$ introduced in \Cref{lemma:Borel-selection}.
    We define $\displaystyle \hp_n = \frac{2}{j_n(j_n-1)}\sum_{1\le j < j' \le j_n} \indic{\halpha_{n,j} = \halpha_{n,j'}}$
    and the estimator of the stickiness regime
    $$\hR_n = \begin{cases}
        0 &\text{ if } \hp_n \in (5/8, 1], \\
        1 &\text{ if } \hp_n \in (1/8, 5/8], \\
        2 &\text{ if } \hp_n \in [0, 1/8].
    \end{cases}$$
\end{definition}

We state some guarantees on $\hR_n$, depending on the stickiness regime of $\alphas$
and asymptotic conditions involving the number of blocks $j_n$.

\begin{theorem}
    \label{thm:estimator-stickiness}
    \begin{enumerate}
        \item In the case where $\alphas$ is sticky:
        \begin{enumerate}
            \item Assume that $\alphas \in \operatorname{int}(F)$ and $n/j_n \to +\infty$. Then $\P{\hR_n = 0} \to 1$.
            \item In addition to the assumptions of part a), assume that $\E{\ell_+'(d(\alphas, X))^2} < +\infty$. 
            Then $\P{\exists N\geq 2, \forall n\geq N, \hR_n = 0} = 1$.
        \end{enumerate}
        \item In the case where $\alphas$ is one-sided partly sticky:
        \begin{enumerate}
            \item Assume that $\alphas \in \operatorname{int}(F)$, $\E{\ell_+'(d(\alphas, X))^2} < +\infty$, $j_n\to +\infty$ and $n/j_n \to +\infty$. 
            Then $\P{\hR_n = 1} \to 1$.
            \item In addition to the assumptions of part a), assume that $\sum_{n=2}^\infty e^{-j_n/512} < +\infty$. 
            Then $\P{\exists N\geq 2, \forall n\geq N, \hR_n = 1} = 1$.
        \end{enumerate}
        \item Consider the case where $\alphas$ is two-sided partly sticky.
        Assume the following:
    \begin{enumerate}[label=(\alph*)]
        \item Assumption~\ref{assum:clt-2} holds: $\E{\ell_+'(d(\alphas, X) + \rho_0)^2} < +\infty$ for some $\rho_0\in (0, \rhos)$.
        \item For each $i\in [m]$, let $\gamma_i:[0,\rho_0] \to T$ be the unit-speed geodesic such that $\gamma_i(0) = \alphas$ and $\gamma_i(\rho_0) \in \Delta_i$.
        For every $i\in I_0$,
        the function $\Phi_i:[0,\rho_0] \to \R$, $t\mapsto \phi(\gamma_i(t))$ satisfies 
        $$\frac{(\Phi_i)_+'(t)}{t^{a_i}} \xrightarrow[t\downarrow 0]{} K_i,$$
        where $a_i$ and $K_i$ are positive constants.
        \item $n/j_n \to +\infty$.
    \end{enumerate}
     Then $\P{\hR_n = 2} \to 1$.
    \end{enumerate}
\end{theorem}

As an immediate consequence, we obtain the following consistency guarantee, which is valid irrespective of the stickiness regime of $\alphas$.

\begin{corollary}[Consistency of $\hR_n$]
    \label{corol:estimator-stickiness}
    Assume conditions 3(a),(b),(c) of \Cref{thm:estimator-stickiness} as well as $j_n \to +\infty$.
    Then $\P{\hR_n = |I_0|} \xrightarrow[n\to\infty]{} 1$. Equivalently, $\hR_n \xrightarrow[n\to \infty]{L^1} |I_0|$.
\end{corollary}

Although the preceding estimator is conceptually simple, it requires computing $j_n$ empirical Fr\'echet $\ell$-means, 
with $j_n\to +\infty$ for consistency across all three stickiness regimes. 
This may be computationally prohibitive when each empirical minimization problem is itself costly, 
particularly since these minimizers are used only to detect pairwise collisions. 
We next propose an alternative estimator that requires computing only one empirical Fr\'echet $\ell$-mean.

\subsection{Via empirical difference quotients}
\label{sec:empirical-difference-quotients}

Let us briefly explain our method. We start by splitting the sample into two parts, $A_n$ and $B_n$. 
We consider the empirical Fr\'echet $\ell$-mean $\halpha_{A_n}$ obtained using only the observations from $A_n$.
Next, we introduce the sphere $\S_n$ centered at $\halpha_{A_n}$ and with a radius chosen so that, with high probability, $\S_n$ intersects each direction $\Delta_i$ exactly once. 
Based on the sphere points and the observations in $B_n$, we compute empirical difference quotients which act as surrogates for the directional derivatives $\phi_{\Delta_i}'(\alphas)$,
and our estimator of $|I_0|$ is defined by comparing these surrogates to a suitable threshold.
The following definition provides the details of our construction.

\begin{definition}
    \label{def:estimator-difference-quotients}
    For each $n\ge 2$, let $A_n$ and $B_n$ be two nonempty deterministic subsets of $[n]$ that form a partition thereof.
    Let $\hmu_{A_n} = \frac{1}{|A_n|} \sum_{k\in A_n} \delta_{X_k}$ and 
    let $\halpha_{A_n} = f_{|A_n|}((X_k)_{k\in A_n})$ be the measurable selection
    from $M(\hmu_{A_n})$ introduced in \Cref{lemma:Borel-selection}.
    Given a sequence $(r_n)_{n\ge 2}$ of positive real numbers, we define the sphere $\S_n=\set{\beta \in T: d(\beta, \halpha_{A_n}) = r_n}$.
    We introduce the empirical objective $\hphi_{B_n}: \alpha \mapsto \frac 1{|B_n|} \sum_{k\in B_n} [\ell(d(\alpha, X_k)) - \ell(d(o, X_k))]$ and 
    we denote by $(\tau_n)_{n\ge 2}$ a sequence of positive thresholds.
    The estimator of the stickiness regime is  
    $$\hR_n' \coloneqq \min\Big(2, \Big|\set[\Big]{\beta \in \S_n : \frac{\hphi_{B_n}(\beta) - \hphi_{B_n}(\halpha_{A_n})}{r_n}\le \tau_n }\Big|\Big).$$
\end{definition}

The measurability of $\hR_n'$ is addressed in \Cref{lemma:estimator-measurable}.
Let us briefly explain the requirements that $|A_n|$, $|B_n|$, $r_n$, $\tau_n$ should satisfy 
for $\hR_n'$ to be a sensible estimator of $|I_0|$.
This will also shed some light on the role of the assumptions in the theorem below.
First, we need the sphere $\S_n$ (which is finite by \Cref{lemma:ball}) to correctly discover the local geometry at $\alphas$, 
in the sense that $\S_n = \set{\beta_i: i\in [m]}$, 
where each $\beta_i$ lies on the segment that connects $\alphas$ to its neighboring vertex in $\Delta_i$.
In order for the sphere to intersect only those edges, we require that $r_n\to 0$.
However, $r_n$ cannot be too small relative to $d(\halpha_{A_n}, \alphas)$,
otherwise $\S_n$ may miss some directions in $\dir_{\alphas}$.
Quantitatively, we require that $d(\halpha_{A_n}, \alphas) = o_\PP(r_n)$,
which in turn is guaranteed by the central limit theorems \ref{thm:two-sided-clt} and \ref{thm:one-sided-clt},
provided that $|A_n|^{-1/(2a_i)}\ll r_n$ for each $i\in I_0$.
Second, we leverage the independence of the subsamples by conditioning on the $(X_k)_{k\in A_n}$,
thus freezing the randomness in $\halpha_{A_n}$. 
Each empirical difference quotient $(\hphi_{B_n}(\beta_i) - \hphi_{B_n}(\halpha_{A_n}))/r_n$
then has fluctuations of order $|B_n|^{-1/2}$ around its mean
$(\phi(\beta_i) - \phi(\halpha_{A_n}))/r_n$; we consequently require $|B_n|^{-1/2} \ll \tau_n$.
For each $i\in I_0$, $(\phi(\beta_i) - \phi(\halpha_{A_n}))/r_n$ has an upper bound of order $r_n^{a_i}$,
hence we also impose $r_n^{a_i} \ll \tau_n$.
Moreover, we add $\tau_n \to 0$, so that directions with positive population derivative are not selected.

\begin{theorem}[Consistency of $\hR_n'$]
    \label{thm:estimator-stickiness-2}
    Assume the following:
    \begin{enumerate}[label=(\alph*)]
        \item Assumption~\ref{assum:clt-2} holds: $\E{\ell_+'(d(\alphas, X) + \rho_0)^2} < \infty$ for some $\rho_0\in (0, \rhos)$.
        \item For each $i\in [m]$, let $\gamma_i:[0,\rho_0] \to T$ be the unit-speed geodesic such that $\gamma_i(0) = \alphas$ and $\gamma_i(\rho_0) \in \Delta_i$.
        For every $i\in I_0$,
        the function $\Phi_i:[0,\rho_0] \to \R$, $t\mapsto \phi(\gamma_i(t))$ satisfies 
        $$\frac{(\Phi_i)_+'(t)}{t^{a_i}} \xrightarrow[t\downarrow 0]{} K_i,$$
        where $a_i$ and $K_i$ are positive constants.
        \item The following asymptotics hold:
        $|A_n|\to \infty$, $|B_n|^{-1/2} = o(\tau_n)$, $r_n \to 0$, $\tau_n \to 0$ and for each $i\in I_0$,
        $|A_n|^{-1/(2a_i)} = o(r_n)$, $r_n^{a_i}=o(\tau_n)$.
    \end{enumerate}
    Then $\P{\hR_n' = |I_0|} \xrightarrow[n\to\infty]{} 1$. Equivalently, $\hR_n' \xrightarrow[n\to \infty]{L^1} |I_0|$.
\end{theorem}

\begin{remark}
    An explicit tuning that satisfies condition (c) is $|A_n|\sim |B_n| \sim n/2$, $r_n = 1/\log n$ and $\tau_n = 1/\log \log (n+1)$.
    Remarkably, it is valid irrespective of the values of the exponents $a_i$.
    In the regular case where each $a_i$ equals $1$, we may pick $r_n = n^{-1/3}$ and $\tau_n = n^{-1/4}$ for instance.
\end{remark}

\section{The special case of Fr\'echet medians}
\label{sec:medians}

We now restrict our focus to Fr\'echet medians, \ie the case where the loss is $\ell:z\mapsto z$.
Among the operations-research community, the median case has generated the most interest,
as it is the most intuitive in applications:
the practitioner looks for a new facility on the network that minimizes the average travel distance to the demand.

\subsection{Further descriptive results}
\label{sec:medians-descriptive}

By \Cref{thm:shape}, the set of medians $M_1(\mu)$ is a geodesic segment.
Conversely, any geodesic segment $[\alpha, \alpha']$ is the set of medians of some distribution:
$[\alpha, \alpha'] = M_1(\mu)$ where $\mu = \frac 12 (\delta_{\alpha} + \delta_{\alpha'})$ is the uniform distribution on $\set{\alpha, \alpha'}$.
Indeed, for this distribution, 
the directional derivatives at $\alpha$ and $\alpha'$ are all nonnegative, hence $\set{\alpha, \alpha'} \subset M_1(\mu)$.
Since $M_1(\mu)$ is convex, it follows that $[\alpha, \alpha'] \subset M_1(\mu)$. 
The reverse inclusion holds by \Cref{prop:localization-support}.

Besides, for a discrete measure with uniform weights $\nu = \frac 1m \sum_{i=1}^m \delta_{x_i}$ on the real line $\R$,
it is well-known that $M_1(\nu)$ is a compact interval of the form $[x_i, x_{i'}]$, possibly with $i=i'$.
We provide a generalization of this fact on a tree.
It is true for any $\mu$, not only for discrete distributions.
Regarding terminology, the \textit{degree of a vertex} refers to the number of edges incident to it.
We say that a vertex is \textit{branching} if its degree is $\ge 3$ and we denote by $D$ the set of branching vertices. 

\begin{proposition}
    \label{prop:endpoints-median}
    Let $\alpha_-$ and $\alpha_+$ denote the endpoints of $M_1(\mu)$, \ie $M_1(\mu) = [\alpha_-,\alpha_+]$.  
    The following inclusion holds: $\set{\alpha_-,\alpha_+}\subset D \cup \supp \mu$.
\end{proposition}

\begin{example}
    When $T$ is the double ray (\ie $T$ is isometric to the real line), $D$ is empty and $\set{\alpha_-,\alpha_+}\subset \supp \mu$.
    We next consider the metric tree $T$ from \Cref{fig:endpoints-median-example}.
    The points $\alpha_1,\alpha_2,\alpha_3,\alpha_4$ are chosen arbitrarily in
    $(v_1,v_3)$, $(v_1,v_4)$, $(v_2,v_5)$, $(v_2,v_6)$ respectively.
    The probability distribution is $\mu = \frac{1}{4}\sum_{i=1}^{4}\delta_{\alpha_i}$.
    Note that $\phi_{v_3}'(v_1) = \phi_{v_4}'(v_1) = \frac 12$ and $\phi_{v_2}'(v_1) = 0$. 
    Hence $v_1\in M_1(\mu)$ and there are no medians on the outer left branches. Similarly, $v_2\in M_1(\mu)$ and there are no medians on the outer right branches.
	Consequently, $M_1(\mu) = [v_1,v_2]$, with $\set{v_1,v_2}\subset D$, which is disjoint from $\supp \mu$.

    \begin{figure}[H]
        \centering
        \begin{tikzpicture}[
            line cap=round,
            line join=round,
            tree edge/.style={line width=0.8pt},
            vertex/.style={circle, fill=black, inner sep=1.5pt},
            tree label/.style={font=\small, inner sep=1pt},
            alpha tick/.style={line width=1.0pt}
        ]

        \coordinate (v1) at (0,0);
        \coordinate (v2) at (4,0);
        \coordinate (v3) at (-1.414,1.414);
        \coordinate (v4) at (-1.414,-1.414);
        \coordinate (v5) at (5.414,1.414);
        \coordinate (v6) at (5.414,-1.414);

        \coordinate (a1) at ($(v1)!0.8!(v3)$);
        \coordinate (a2) at ($(v1)!0.4!(v4)$);
        \coordinate (a3) at ($(v2)!0.42!(v5)$);
        \coordinate (a4) at ($(v2)!0.62!(v6)$);

        \draw[tree edge] (v1) -- (v2);
        \draw[tree edge] (v1) -- (v3);
        \draw[tree edge] (v1) -- (v4);
        \draw[tree edge] (v2) -- (v5);
        \draw[tree edge] (v2) -- (v6);

        \node[vertex] at (v1) {};
        \node[vertex] at (v2) {};
        \node[vertex] at (v3) {};
        \node[vertex] at (v4) {};
        \node[vertex] at (v5) {};
        \node[vertex] at (v6) {};

        \draw[alpha tick] ($(a1)+(-0.08,-0.08)$) -- ($(a1)+(0.08,0.08)$);
        \draw[alpha tick] ($(a2)+(-0.08,0.08)$)  -- ($(a2)+(0.08,-0.08)$);
        \draw[alpha tick] ($(a3)+(-0.08,0.08)$)  -- ($(a3)+(0.08,-0.08)$);
        \draw[alpha tick] ($(a4)+(-0.08,-0.08)$) -- ($(a4)+(0.08,0.08)$);

        \node[tree label, anchor=east]       at ($(v1)+(-0.18,0.00)$)   {$v_1$};
        \node[tree label, anchor=west]       at ($(v2)+(0.18,0.00)$)    {$v_2$};

        \node[tree label, anchor=south east] at ($(v3)+(-0.06,0.10)$)   {$v_3$};
        \node[tree label, anchor=north east] at ($(v4)+(-0.06,-0.10)$)  {$v_4$};
        \node[tree label, anchor=south west] at ($(v5)+(0.06,0.10)$)    {$v_5$};
        \node[tree label, anchor=north west] at ($(v6)+(0.06,-0.10)$)   {$v_6$};

        \node[tree label, anchor=west]       at ($(a1)+(0.15,0.1)$)    {$\alpha_1$};
        \node[tree label, anchor=east]       at ($(a2)+(0.4,-0.3)$)   {$\alpha_2$};
        \node[tree label, anchor=east]       at ($(a3)+(0.1,0.3)$)   {$\alpha_3$};
        \node[tree label, anchor=east]       at ($(a4)+(-0.13,-0.20)$)  {$\alpha_4$};

        \end{tikzpicture}
        \caption{A metric tree $T$ such that $D = \set{v_1, v_2}$.}
        \label{fig:endpoints-median-example}
    \end{figure}
\end{example}

\begin{remark}
    Hakimi \cite{hakimi1964opti} states a weaker statement: when $\mu$ is discrete and supported on $\mathcal V$, 
    he proves that $M_1(\mu) \cap \mathcal V \neq \emptyset$, where $\mathcal V$ denotes the set of all vertices in $T$.
\end{remark}

A measure $\nu$ on $\R$ has at least two medians if and only if there exist $m_1<m_2$ such that 
$\nu\left((-\infty,m_1]\right) = \nu\left([m_2,+\infty)\right)=\frac 12$ \cite[Corollary 2.6]{romon2022quantiles}.
The next proposition is an extension of this fact to metric trees.
As a byproduct, we derive a tractable guarantee that $M_1(\mu)$ is a singleton.

\begin{proposition}
    \label{prop:unique-median}
    \begin{enumerate}
        \item Assume that $M_1(\mu) = [\alpha_-,\alpha_+]$ with $\alpha_- \neq \alpha_+$.
        Then the subsets $T\setminus \Delta_{\alpha_- \to \alpha_+}$ and $T\setminus \Delta_{\alpha_+ \to \alpha_-}$
        are disjoint, closed, convex and each has probability $1/2$ \wrt $\mu$. 
        \item If $\alpha_-\neq \alpha_+$ satisfy $\mu(T\setminus \Delta_{\alpha_- \to \alpha_+}) = \mu(T\setminus \Delta_{\alpha_+ \to \alpha_-}) = 1/2$,
        then $[\alpha_-,\alpha_+] \subset M_1(\mu)$.
        \item $\mu$ has more than one Fr\'echet median if and only if there exist two disjoint closed convex subsets $C_1,C_2\subset T$ such that $\mu(C_1)=\mu(C_2)=\frac 12$.
    \end{enumerate}
\end{proposition}

\begin{remark}
    If the support of $\mu$ is a convex subset of $T$ (\eg if $\mu$ has a density $f$ with respect to the Lebesgue measure on $T$, 
    such that $f$ is positive on a convex subset $C$ and zero on $T\setminus C$), \Cref{prop:unique-median} entails that $\mu$ has a unique Fr\'echet median.
\end{remark}

\begin{remark}
    If $n$ is odd and $\mu = \frac 1n \sum_{k=1}^{n} \delta_{x_k}$,
    then $\mu$ cannot put mass $1/2$ on any subset of $T$.
    By \Cref{prop:unique-median}, $\mu$ thus has a unique Fréchet median.
\end{remark}

\subsection{Confidence regions for Fr\'echet medians}

We recall that in the median case, $\phi_{\Delta}'(\alpha)$ has the closed form $\phi_{\Delta}'(\alpha) = 1-2\mu(\Delta)$,
which depends only on the mass that $\mu$ puts on $\Delta$; it does not involve the metric $d$.
In order to control the empirical directional derivatives at a given Fr\'echet median $\alphas$, 
we make use of the following lemma regarding the binomial distribution.

\begin{lemma}
    \label{lemma:binomial}
    Let $n\ge 1$, $k > n/2$ be integers, and 
    define the function $$f:[0,1] \to [0,1], \quad p\mapsto \P{\Sigma_{n,p} \ge k},$$
    where $\Sigma_{n,p}\sim \Bin(n,p)$.
    Then $f$ is convex on $[0,1/2]$, vanishes at $0$, and therefore
    $\forall p\in [0,1/2]$, $\P{\Sigma_{n,p} \ge k} \le 2p \P{\Sigma_{n,1/2} \ge k}$.
\end{lemma}

We are now ready to state a deviation bound for the empirical derivatives at any $\alphas\in  M_1(\mu)$.

\begin{proposition}
    \label{prop:control-derivatives-median}
    Let $\alphas\in M_1(\mu)$ with $\dir_{\alphas} = \set{\Delta_1,\ldots,\Delta_m}$.
    For any $n\ge 1$ and $t>0$, 
    \begin{equation}
        \label{eq:control-derivatives-median}
        \P*{\min_{1\le i\le m} \hphi_{\Delta_i}'(\alphas)\le -t}
        \le 2\P*{\Sigma \ge \ceil*{\tfrac{n(1+t)}{2}}}
        \le 2e^{-\frac{nt^2}2},
    \end{equation}
    where $\Sigma \sim \Bin(n, 1/2)$.
\end{proposition}

\begin{remark}
    A naive union bound over the $m$ directions in $\dir_{\alphas}$ would lead to an extra factor $m$ in the upper bounds.
    By leveraging \Cref{lemma:binomial}, we obtain upper bounds  
    that are independent of the number of directions at $\alphas$.
\end{remark}

\begin{remark}
    In the median case, (empirical) directional derivatives belong to the interval $[-1,1]$.
    Consequently, when $t$ is larger than $1$, the two probability terms in \eqref{eq:control-derivatives-median} vanish.
\end{remark}

According to \Cref{prop:optimality}, $\alphas\in M_1(\mu)$ if and only if all the directional derivatives at $\alphas$ are nonnegative.
To obtain a data-driven confidence region for $M_1(\mu)$, we may therefore consider the points in the tree where all the empirical directional derivatives are larger than $-t$ for some threshold $t>0$.
This motivates the following definition.

\begin{definition}
    \label{def:region}
    Given $n\ge 1$ and a threshold $t\in (0,1]$, we define the random subset $$\hC_t = \set[\big]{\alpha\in T: \forall \Delta \in \dir_\alpha, \, \hphi_{\Delta}'(\alpha) > -t},$$
    as well as the integer $s = \ceil*{\frac{n(1+t)}{2}}$.
\end{definition}

Next, we state some properties enjoyed by the region $\hC_t$.

\begin{proposition}
    \label{prop:region}
    \begin{enumerate}
        \item $\hC_t$ satisfies 
        $\hC_t = \set[\big]{\alpha\in T: \forall \Delta \in \dir_\alpha, \, \sum_{k=1}^{n} \indic{X_k\in \Delta} < s}$ and
        \begin{equation}
            \label{eq:region-equality}
            \hC_t = \bigcap_{\substack{K \subset [n]\\|K|=s}} \conv(\set{X_k:k\in K}).
        \end{equation}
        \item $\hC_t$ is nonempty, compact and convex, with $M_1(\hmu_n) \subset \hC_t \subset \conv(\set{X_k:1\le k\le n})$.
        It is a measurable closed subset, in the sense of \Cref{sec:estimation-trees}.
        \item Let $\Sigma \sim \Bin(n, 1/2)$. 
        Then $\P{M_1(\mu) \not\subset \hC_t} \underset{(*)}{\le} 2 \P{\Sigma\ge s}$.
        If $\mu$ has more than one median, then $(*)$ is an equality.
    \end{enumerate}
\end{proposition}

\begin{remark}
    The right-hand side $2 \P{\Sigma\ge s}$ is independent of $\mu$, hence the miscoverage bound holds uniformly over all probability measures $\mu$.
    Moreover, this bound is sharp for measures with more than one median (\eg the uniform distribution on two points) and the factor $2$ cannot be improved.
\end{remark}

The regions are nested in the sense that $t < t'$ implies $\hC_{t} \subset \hC_{t'}$.
Moreover, $\hC_t$ depends on $t$ only through the integer $s=\ceil{n(1+t)/2}$.
To construct a $1-\eta$ confidence region for the whole set of medians $M_1(\mu)$,
we pick the smallest value of $s$ such that $2 \P{\Sigma\ge s} \le \eta$.
This is possible provided that $\eta \ge 1/{2^{n-1}}$.

\begin{theorem}
    \label{thm:confidence-region}
    Let $n\ge 1$ and $\eta\in (0,1]$ be such that $1/{2^{n-1}} \le \eta$. 
    Define $\us = \min\set[\big]{s\in \set{\floor{n/2}+1,\ldots,n}: 2\P{\Sigma\ge s} \le \eta}$
    as well as $\ut = \frac 2n \us -1$,
    where $\Sigma \sim \Bin(n, 1/2)$.
    Then $\hC_\ut$ is a confidence region for $M_1(\mu)$ with confidence level $1-\eta$, in the sense that 
    $$\P[\big]{M_1(\mu) \subset \hC_\ut} \ge 1-\eta.$$
\end{theorem}

\begin{remark}
    The values of $\us$, $\ut$ and the construction of $\hC_\ut$ do not depend on the data-generating $\mu$.
    The set $\hC_\ut$ is therefore a distribution-free non-asymptotic $1-\eta$ confidence region.
\end{remark}

\begin{remark}
    When $T$ is the double ray (\ie $T$ is isometric to the real line),
    we denote by $X_{(1)}\le \ldots \le X_{(n)}$ the order statistics obtained 
    after orienting $T$. It is straightforward to show that $\hC_\ut = [X_{(n-\us+1)}, X_{(\us)}]$,
    which coincides with the classical distribution-free $1-\eta$ confidence interval for the median 
    (see, \eg \cite[Section~7.1]{david2003order}).
\end{remark}

Consider the setting of the real line, where $X_{(1)}\le \ldots \le X_{(n)}$ denote the order statistics of an \iid $n$-sample from $\mu$.
We let $(s_n)_{n\geq 1}$ be a sequence of integers such that $s_n / n \to 1/2$.
It is well-known that if $\mu$ has a unique median $\alphas$, then $X_{(s_n)} \xrightarrow[]{a.s.} \alphas$.
Consequently, for a fixed $\eta$, the confidence interval $[X_{(n-\us+1)}, X_{(\us)}]$ shrinks to $\alphas$ with probability $1$.
We next provide a generalization of this fact to metric trees.

\begin{proposition}
    \label{prop:region-shrink}
    Assume that $\mu$ has a unique Fr\'echet median: $M_1(\mu) = \set{\alphas}$.
    \begin{enumerate}
        \item Let $(t_n)_{n\ge 1}$ be a sequence of real numbers in $(0,1]$ such that $t_n\to 0$ as $n\to \infty$.
        Then $\sup_{\alpha \in \hC_{t_n}} d(\alpha, \alphas) \xrightarrow[]{a.s.} 0$.
        \item Fix $\eta \in (0,1]$ and for each $n\ge \log_2(1/\eta) + 1$, let $\ut$ be as in \Cref{thm:confidence-region}.
        Then $\sup_{\alpha \in \hC_{\ut}} d(\alpha, \alphas) \xrightarrow[]{a.s.} 0$.
    \end{enumerate}
\end{proposition}

\newpage 
\bibliographystyle{imsart-number}
\bibliography{biblio.bib}      

\appendix
\section{Proofs}

\subsection{Proofs for Section~\ref{sec:frechet}}

\begin{lemma}
    \label{lemma:countable}
    A metric tree $T$ as in \Cref{def:tree} has countably many vertices.
\end{lemma}

\begin{proof}[Proof of \Cref{lemma:countable}]
    Consider an arbitrary vertex $v\in T$ and for each $n\ge 1$ denote by $V_n$ the set of those vertices $w$ such that the path from $v$ to $w$
    goes through exactly $n$ edges. For instance, $V_1$ is the set of neighboring vertices of $v$.
    By the local finiteness assumption and a straightforward induction on $n$, each $V_n$ is finite.
    The set of all vertices of $T$ is $\set{v} \cup\bigcup_{n\ge 1} V_n$, which is countable.
\end{proof}

\begin{lemma}
    \label{lemma:ball}
    Let $T$ denote a metric tree that satisfies Assumption~\ref{assum:lower-bound}.
    \begin{enumerate}
        \item If $B$ is a bounded subset of $T$, then $B$ contains finitely many vertices of $T$ and intersects finitely many edges of $T$.
        \item If $\S \subset T$ is a sphere, then $\S$ is finite.
    \end{enumerate}
\end{lemma}

\begin{proof}[Proof of \Cref{lemma:ball}]
    1. Since $B$ is contained in a ball centered at a vertex, we assume without loss of generality that $B$ is the closed ball $\overline B(v,r)$ where $v$ is a vertex and $r>0$.
    By Assumption~\ref{assum:lower-bound}, we can find a lower bound $\rho>0$ on the lengths of all edges.
    Any vertex $w$ contained in $B$ is then separated from $v$ by at most $\floor{r/\rho}$ consecutive edges.
    Therefore, the graph-theoretic tree $T'$ rooted at $v$ with all vertices and edges contained in $B$ has depth at most $\floor{r/\rho}$.
    By condition b) in \Cref{def:tree}, there are finitely many vertices at each depth level in $T'$.
    Thus $T'$ has finitely many vertices, and $B$ contains finitely many vertices of $T$.

    We denote by $\mathcal E'$ the set of edges of $T$ that intersect $B$.
    Then any edge $e$ that intersects $B$ must have at least one endpoint in $B$:
    indeed, if $x\in e\cap B$ then the geodesic from $v$ to $x$ contains at least one endpoint of $e$, say $w$, and $d(v,w)\le d(v,x)\le r$,
    hence $w\in B$.
    By the previous paragraph, $B$ contains finitely many vertices of $T$, and combined with the local finiteness assumption we obtain that $\mathcal E'$ is finite.

    2. By the first item, $\S$ intersects finitely many edges of $T$.
    Moreover, the intersection of $\S$ with an edge has cardinality at most $2$, hence the claim. 
\end{proof}

\begin{proof}[Proof of \Cref{prop:topology}]
    1. 
    Let $K$ denote a closed and bounded subset of $T$.
    Denoting by $\mathcal E'$ the set of edges of $T$ that intersect $K$, we have $K=\bigcup_{e\in \mathcal E'}(K\cap e)$.
    By \Cref{lemma:ball}, $\mathcal E'$ is finite.
    Each $K\cap e$ is isometric to a closed and bounded subset of $\R$, and is therefore compact.
    Thus $K$ is compact as a finite union of compact subsets.

    Since closed balls are compact, $T$ is locally compact. That a proper metric space is complete and separable is shown in \cite[Proposition~A.15]{dull2022spaces}.
    
    2. A metric tree is CAT($0$) \cite[Example 1.15(5) p.~167]{bridson1999metric}. Since $T$ is also complete, it is a Hadamard space. 
\end{proof}

The following lemma exhibits basic regularity properties of the loss.

\begin{lemma}
    \label{lemma:integral-rep}
    Let $\ell:[0,+\infty)\to \R$ be a loss function.
    \begin{enumerate}
        \item The left derivative $\ell_-':(0,+\infty)\to (0,+\infty)$ and right derivative $\ell_+':[0,+\infty)\to [0,+\infty)$
        of $\ell$ exist and are increasing.
        \item $\ell_-'$ is left-continuous, $\ell_+'$ is right-continuous, $\ell_-'(z+) = \ell_+'(z)$ for $z\ge 0$ and $\ell_+'(z-) = \ell_-'(z)$ for $z>0$.
        \item For any bounded interval $I\subset [0,+\infty)$, the restriction $\ell_{\mid I}$ is Lipschitz continuous.
        \item For every $z\in [0,+\infty)$, $\ell(z) = \ell(0) + \int_0^z \ell_-'(t) \diff t = \ell(0) + \int_0^z \ell_+'(t) \diff t$.
    \end{enumerate}
\end{lemma}

\begin{proof}[Proof of \Cref{lemma:integral-rep}]
    1. Since $\ell$ is convex, it has finite left and right derivative at each $z>0$, with $\ell_-'$ and $\ell_+'$ being increasing 
    \cite[Proposition~6.2.7]{papadopoulos2014metric}.
    Since $\ell$ is increasing, $\ell_-'$ and $\ell_+'$ are nonnegative.
    If $\ell_-'(z)=0$ for some $z>0$, then by \cite[Proposition~6.2.7]{papadopoulos2014metric} $\ell$ is constant on some open interval,
    which contradicts the strict monotonicity of $\ell$.
    
    2. For the left-continuity of $\ell_-'$ and the right-continuity of $\ell_+'$, see \cite[Remark 4.2.2 p.\ 26]{hiriart1993convex}.
    Fix some $z\ge 0$. Then for every $t\in (z,\infty)$, $\ell_+'(z) \le \ell_-'(t) \le \ell_+'(t)$ hence 
    by letting $t\downarrow z$ and by the right-continuity of $\ell_+'$ 
    we obtain that $ \ell_-'(z+) = \ell_+'(z)$. The proof of the other equality is similar.

    3. For a fixed $Z>0$ and $0\le z_1<z_2\le Z$, we have the estimate
    $$0\le \frac{\ell(z_2)-\ell(z_1)}{z_2-z_1} \le \ell_-'(z_2) \le \ell_-'(Z),$$
    hence the claim.

    4. By the last point, $\ell$ is absolutely continuous on compact intervals. Let $Z>0$ be fixed.
    By the fundamental theorem of calculus \cite[Theorem~7.18]{rudin1987real}, $\ell$ is differentiable \aev on $[0,Z]$ and for every 
    $z\in [0,Z]$, $\ell(z) - \ell(0) =  \int_0^z f_Z(t) \diff t$, 
    where $f_Z$ denotes a derivative of $\ell$.
    Since $f_Z$ and $\ell_+'$ coincide \aev on $[0,Z]$, we obtain $\ell(z) - \ell(0) =  \int_0^z \ell_+'(t) \diff t$.
    We proceed similarly with $\ell_-'$.
\end{proof}

\begin{proof}[Proof of \Cref{prop:phi-trees}]
    1. By the convexity of the loss and the triangle inequality, we obtain the following estimate:
    \begin{equation}
        \label{eq:integrand-estimate}
        \forall (\alpha, x)\in T^2, \quad -\ell_+'(d(o,x)) d(\alpha, o) \le \ell(d(\alpha,x)) - \ell(d(o,x)) \le \ell_+'(d(\alpha,x)) d(\alpha, o).
    \end{equation}  
    Since $o$ is in $F$, the leftmost term in \eqref{eq:integrand-estimate} is integrable, 
    thus $\phi(\alpha)$ is well-defined for each $\alpha\in T$ as a quantity in $\R\cup\set{+\infty}$.
    
    Since $T$ is Hadamard, by \cite[Example 8.4.7(i)]{papadopoulos2014metric} the map 
    $\alpha\mapsto d(\alpha, x)$ is convex for each $x\in T$, and by the convexity and monotonicity of $\ell$  
    we obtain convexity of $\alpha\mapsto \ell(d(\alpha, x))$.
    That $\phi$ is convex follows by integration.
    Since $\phi(o)=0< +\infty$, we note that $\phi$ is proper.

    Finally, if $\int_T \ell_+'(d(\alpha,x)) \diff\mu(x)< +\infty$ then by the rightmost inequality in \eqref{eq:integrand-estimate}, $\phi(\alpha)$ is finite.

    2. To prove lower semicontinuity, we fix $t\in \R$ and it remains to show that the sublevel set $\set{\alpha \in T: \phi(\alpha)\le t}$ is a closed subset of $T$.
    This follows from the leftmost inequality in \eqref{eq:integrand-estimate} and Fatou's lemma for functions with an integrable lower bound.

    We fix $\alpha\in \operatorname{int}(\dom \phi)$ and we show that $\phi$ is continuous at $\alpha$. 
    For convenience, we make use of the concepts and notations introduced in \Cref{sec:directions}:
    we write $\dir_\alpha = \set{\Delta_1,\ldots,\Delta_m}$.
    There exists $r>0$ such that 
    $B(\alpha, r) \subset \dom \phi$,
    and we let $r'$ denote the distance between $\alpha$ and its closest neighboring vertex.
    Let $\rho = \min(r,r')/2$ and for each $i\in [m]$
    denote by $\alpha_i$ the point in $\Delta_i$ that satisfies $d(\alpha, \alpha_i) = \rho$.
    Let us also put $C = \max_{1\le i \le m} |\phi(\alpha_i) - \phi(\alpha)|$.
    Fix $\beta$ such that $0<d(\beta,\alpha)<\rho$ and find $i$ such that $\beta\in (\alpha, \alpha_i)$.
    The convexity of $\phi$ with $t=d(\beta, \alpha)/\rho$ yields
    $\phi(\beta) \le (1-t)\phi(\alpha) + t \phi(\alpha_i) \le \phi(\alpha) + tC = \phi(\alpha) + \frac{C}{\rho} d(\beta, \alpha)$,
    from which we obtain that $\limsup_{\beta \to \alpha} \phi(\beta)\le \phi(\alpha)$.
    Combined with lower semicontinuity, it follows that $\lim_{\beta \to \alpha} \phi(\beta)=\phi(\alpha)$.

    3. Suppose that $\ell_+'$ is bounded from above, so that $\ell_+'(z)$ converges to some $L\in (0,+\infty)$ as $z$ goes to $\infty$.
    By \Cref{eq:integrand-estimate} we obtain the estimate 
    \begin{equation}
        \label{eq:integrand-estimate2}
        \forall (\alpha, x)\in T^2, \quad |\ell(d(\alpha,x)) - \ell(d(o,x))| \le L d(\alpha, o),
    \end{equation}
    which implies that $\phi$ is finite-valued.
    The convexity of $\ell$ yields $\ell(z)/z \ge \ell_+'(1)(1-1/z)+\ell(1)/z$ for each $z>0$.
    We set $M=2/\ell_+'(1)$, so that we can find some $A\ge 0$ such that  $z\ge A \implies z \le M \ell(z)$.
    We note that $\frac{\phi(\alpha)}{\ell(d(\alpha,o))} = \int_T \frac{\ell(d(\alpha,x))-\ell(d(o,x))}{\ell(d(\alpha,o))}\diff\mu(x)$
    and we apply the dominated convergence theorem to this last integral.
    We assume that $d(\alpha,o)\ge A$ and we fix some $x\in T$.
    The same convexity estimate shows that $\ell(z)\to \infty$ as $z\to\infty$, hence the ratio $\ell(d(o,x)) / \ell(d(\alpha,o))$
    goes to $0$ as $d(\alpha,o)\to \infty$.
    Regarding the other part of the integrand, note that by exchanging the roles of $\alpha$ and $x$ in \eqref{eq:integrand-estimate2}, we have
    $\left|\frac{\ell(d(\alpha,x))}{\ell(d(\alpha,o))}  - 1\right| \le \frac{L d(o,x)}{\ell(d(\alpha,o))}$,
    hence the ratio $\ell(d(\alpha,x)) / \ell(d(\alpha,o))$
    goes to $1$ as $d(\alpha,o)\to \infty$,
    and so does the integrand $\frac{\ell(d(\alpha,x))-\ell(d(o,x))}{\ell(d(\alpha,o))}$.
    Regarding domination, we observe that by \eqref{eq:integrand-estimate2},
    $\left|\frac{\ell(d(\alpha,x))-\ell(d(o,x))}{\ell(d(\alpha,o))}\right| \le \frac{L d(\alpha,o)}{\ell(d(\alpha,o))} \le LM$.
    By dominated convergence, $\frac{\phi(\alpha)}{\ell(d(\alpha,o))} \to 1$ as $d(\alpha,o)\to \infty$.

    Suppose now that $\ell_+'(z)\xrightarrow[z\to \infty]{} +\infty$.
    For convenience, let $r_\alpha = d(\alpha,o)/4$ for each $\alpha\in T$.
    For any $(\alpha,x)\in T^2$, observe that
    \begin{align*}
        &\ell(d(\alpha,x))-\ell(d(o,x)) 
        \\&= [\ell(d(\alpha,x))-\ell(d(o,x))]\indic{d(o,x)\le r_\alpha} 
            +[\ell(d(\alpha,x))-\ell(d(o,x))]\indic{d(o,x)> r_\alpha} 
        \\&\ge [\ell(d(\alpha,o)-r_\alpha)-\ell(r_\alpha)]\indic{d(o,x)\le r_\alpha} 
        -[\ell_+'(d(o,x))d(\alpha,o)]\indic{d(o,x)> r_\alpha} 
        \\&\ge \ell_+'(r_\alpha)(d(\alpha,o)-2r_\alpha)\indic{d(o,x)\le r_\alpha} 
        -d(\alpha,o) \ell_+'(d(o,x))\indic{\ell_+'(d(o,x))\ge \ell_+'(r_\alpha)} 
        \\&= \frac{\ell_+'(r_\alpha)}2 d(\alpha,o)\indic{d(o,x)\le r_\alpha} 
        -d(\alpha,o) \ell_+'(d(o,x))\indic{\ell_+'(d(o,x))\ge \ell_+'(r_\alpha)},
    \end{align*}
    from which we obtain the following lower bound:
    \begin{align*}
        \phi(\alpha) 
        \ge \frac{\ell_+'(r_\alpha)}2 d(\alpha,o)\int_T \indic{d(o,x)\le r_\alpha} \diff\mu(x)
        -d(\alpha,o)\int_T \ell_+'(d(o,x))\indic{\ell_+'(d(o,x))\ge \ell_+'(r_\alpha)} \diff\mu(x).
    \end{align*}
    Combining the fact that $o\in F$ with the limit $\ell_+'(r_\alpha)\to \infty$ as $d(\alpha,o)\to \infty$, 
    we find that the integral
    $\int_T \ell_+'(d(o,x))\indic{\ell_+'(d(o,x))\ge \ell_+'(r_\alpha)} \diff\mu(x)$
    converges to $0$.
    Moreover, $\int_T \indic{d(o,x)\le r_\alpha} \diff\mu(x)$ goes to $1$,
    hence the quantity $\frac{\ell_+'(r_\alpha)}2 \int_T \indic{d(o,x)\le r_\alpha} \diff\mu(x)$
    goes to $+\infty$,
    which concludes the proof.

    4. Since $\phi$ is convex, proper, lower semicontinuous and coercive, the subset $M(\mu)$ is nonempty, closed, convex and bounded (see \cite[Lemma~2.2.19 and Example 2.1.3]{bacak2014convex}).
    By \Cref{prop:topology}, $T$ is proper, hence $M(\mu)$ is compact.
    
    5. 
    If $\ell_+'$ is strictly increasing, by \Cref{lemma:integral-rep} the following inequality between slopes holds:
    for any $0\le z_1<z_2<z_3$, 
    \begin{equation*}
       \frac{\ell(z_2)-\ell(z_1)}{z_2-z_1} = \frac{1}{z_2-z_1} \int_{z_1}^{z_2} \ell_+'(t) \diff t
       < \ell_+'(z_2) 
       < \frac{1}{z_3-z_2} \int_{z_2}^{z_3} \ell_+'(t) \diff t = \frac{\ell(z_3)-\ell(z_2)}{z_3-z_2},
    \end{equation*} 
    and $\ell$ is strictly convex by \cite[Proposition~6.2.1]{papadopoulos2014metric}. 
    If $\ell_+'$ is not strictly increasing, then $\ell_+'$ is constant on some open interval, hence by \cite[Proposition~6.2.7]{papadopoulos2014metric} $\ell$ is affine on this interval,
    thus $\ell$ is not strictly convex.

    We suppose now that $\ell$ is strictly convex. 
    Let us consider a geodesic $\gamma:[0,1]\to T$ such that 
    $\gamma_0\neq \gamma_1$, $\phi(\gamma_0)<+\infty$ and $\phi(\gamma_1)<+\infty$ and let us fix $t\in (0,1)$.
    We prove the following pointwise estimate, which concludes the proof:
    \begin{equation}
        \label{eq:pointwise-bound}
        \forall x\in T, \quad  \ell(d(\gamma_t, x)) < (1-t) \ell(d(\gamma_0, x)) + t \ell(d(\gamma_1, x)).
    \end{equation}
    Let us pick some $x\in T$. 
    First, assume that $d(\gamma_0, x) \neq d(\gamma_1, x)$. 
    As argued earlier, $\alpha \mapsto d(\alpha, x)$ is convex, hence
    $d(\gamma_t, x) \le (1-t) d(\gamma_0, x) + t d(\gamma_1, x)$. The monotonicity and strict convexity of $\ell$ then yield \eqref{eq:pointwise-bound}.
    Second, assume that $d(\gamma_0, x) = d(\gamma_1, x)$.
    Since $T$ is CAT$(0)$, the metric $d$ satisfies the following inequality (see, \eg \cite[Theorem~1.3.3(iii)]{bacak2014convex})
    $$d^2(\gamma_t, x) \le (1-t) d^2(\gamma_0, x) + t d^2(\gamma_1, x) - t(1-t)d^2(\gamma_0, \gamma_1) < d^2(\gamma_0, x).$$
    The loss is strictly increasing, hence $\ell(d(\gamma_t, x)) < \ell(d(\gamma_0, x))$ and \eqref{eq:pointwise-bound} holds.
\end{proof}

\begin{proof}[Proof of \Cref{prop:growth}]
    1. The proof is straightforward when \ref{assum:bounded} is in force.
    Assume \ref{assum:growth}. We fix $\alpha\in T$ and we show that $\alpha\in F$.
    Let $p\ge 1$ denote an integer that satisfies $p \ge d(\alpha,o)$.
    Note that 
    $$\int_T \indic{d(o,x)< R} \ell_+'(d(\alpha,x))\diff \mu(x)
    \le \ell_+'(R+d(\alpha,o)),$$
    and
    \begin{align*}
        \int_T \indic{d(o,x)\ge  R} \ell_+'(d(\alpha,x))\diff \mu(x)
        &\le \int_T \indic{d(o,x)\ge  R} \ell_+'(d(o,x)+d(\alpha,o))\diff \mu(x) 
        \\&\le \int_T \indic{d(o,x)\ge  R} \ell_+'(d(o,x)+ p)\diff \mu(x)
        \\&\le C\int_T \indic{d(o,x)\ge R} \ell_+'(d(o,x)+(p-1))\diff \mu(x)
        \\&\le C^p \int_T  \ell_+'(d(o,x))\diff \mu(x),
    \end{align*}
    which yields $\int_T \ell_+'(d(\alpha,x))\diff \mu(x)< +\infty$.
    Since $\alpha$ was arbitrary, it follows that $F=T$, hence by \Cref{prop:phi-trees} $F=\dom \phi = T$.

    2. Assume first that $\ell_+'$ is concave. 
    It follows that for each $z\ge 0$
    $$\ell_+'(z) 
    = \ell_+'\Big(\frac{z}{z+1}(z+1) + \frac{1}{z+1}0\Big)
    \ge \frac{z}{z+1} \ell_+'(z+1) + \frac{1}{z+1} \ell_+'(0)
    \ge \frac{z}{z+1} \ell_+'(z+1),
    $$
    and therefore whenever $z\ge 1$ we have 
    $\ell_+'(z+1)\le 2\ell_+'(z)$.

    3. By the hypothesis we have for $z\ge R$, $$\frac{\ell_+'(z+1)}{\ell_+'(z)}\le \frac{M e^{C_3}}{m}\exp\big(C_1[\log^{C_2}(z+1) - \log^{C_2}(z)]\big)$$
    and it suffices to observe that the function $z\mapsto \log^{C_2}(z+1) - \log^{C_2}(z)$ is bounded on $[R,+\infty)$.
\end{proof}

\begin{lemma}
    \label{lemma:directions-equivalence}
    Let $\alpha\in T$. 
    \begin{enumerate}
        \item The relation $\sim_\alpha$ on $T\setd{\alpha}$ is an equivalence relation.
        \item For every $v,w\in T\setd{\alpha}$, $v\sim_\alpha w \iff \alpha\notin[v,w]$.
    \end{enumerate}
\end{lemma}

\begin{proof}[Proof of \Cref{lemma:directions-equivalence}]
    1. Only transitivity requires justification.
    Assume that $u,v,w$ are such that $u \sim_\alpha v$ and $v \sim_\alpha w$.
    Let us pick $x\in ([\alpha,u]\cap[\alpha,v])\setd{\alpha}$ and $y\in ([\alpha,v]\cap[\alpha,w])\setd{\alpha}$.
    Since $x$ and $y$ both lie in $[\alpha,v]$ and since a geodesic segment is isometric to an interval, it follows that 
    $x\in [\alpha, y]$ or $y\in [\alpha, x]$.
    By symmetry, we may only consider the case where $x\in [\alpha, y]$.
    We note that $[\alpha, y]\subset [\alpha,w]$, hence $x\in ([\alpha,u]\cap[\alpha,w])\setd{\alpha}$
    and $u \sim_\alpha w$.

    2. Let us show that $[\alpha,v]\cap[\alpha,w] = \set{\alpha} \iff \alpha\in[v,w]$.
    We suppose first that $[\alpha,v]\cap[\alpha,w] = \set{\alpha}$.
    We denote by $\gamma$ the unit-speed geodesic from $v$ to $\alpha$ defined on $[-d(\alpha,v), 0]$
    and $\widetilde \gamma$ the unit-speed geodesic from $\alpha$ to $w$ defined on $[0, d(\alpha,w)]$.
    We then define
   \begin{align*}
        \Gamma \colon [-d(\alpha,v), d(\alpha,w)] &\to T \\
        t &\mapsto \begin{cases}
            \gamma(t) \text{ if } t \in [-d(\alpha,v), 0], \\
            \widetilde \gamma(t) \text{ if } t \in [0, d(\alpha,w)].
        \end{cases}
    \end{align*}
    Since $[\alpha,v]\cap[\alpha,w] = \set{\alpha}$, we derive that $\Gamma$ is a unit-speed geodesic from $v$ to $w$.
    Moreover, $\alpha = \Gamma_0 \in [v,w]$.
    Conversely, we assume now that $\alpha\in[v,w]$.
    It follows that $d(v,w) = d(\alpha,v) + d(\alpha, w)$ and 
    we may introduce $\Gamma$ the unit-speed geodesic from $v$ to $w$ defined on $[-d(\alpha,v), d(\alpha,w)]$.
    Since $\Gamma$ is injective, we have $\Gamma([-d(\alpha,v), 0]) \cap \Gamma([0, d(\alpha,w)]) = \set{\alpha}$,
    hence $[\alpha,v]\cap[\alpha,w] = \set{\alpha}$.
\end{proof}

\begin{proof}[Proof of \Cref{prop:directions}]
    1. Let $\Delta\in\dir_\alpha$. Since $\Delta$ is an equivalence class, it is nonempty.
    Let $x,y\in\Delta$, hence $\alpha \notin [x,y]$. For each $z\in [x,y]$, $[x,z]$ is a subset of $[x,y]$, thus $\alpha \notin [x,z]$,
    hence $x \sim_\alpha z$, $z\in \Delta$ and $\Delta$ is convex.

    2. Let us consider $\Delta \in \dir_\alpha$ and $x\in \Delta$.
    There is a unique $i\in [m]$ such that $[\alpha, x]\cap [\alpha, v_i]\neq \set{\alpha}$, 
    hence $x\sim_\alpha v_i$ and $\Delta = \Delta_{\alpha \to v_i}$.
    We have proved that $\dir_\alpha \subset \set{\Delta_{\alpha\to v_i}:1\le i\le m}$; 
    the reverse inclusion is trivial.

    3. Consider a direction $\Delta\in \dir_\alpha$.
    It is straightforward to show that $\Delta$ is an open subset of $T\setd{\alpha}$.
    Its complement $(T\setd{\alpha})\setminus \Delta$ is the finite union of the other directions in $\dir_\alpha$,
    hence it is open.
    Consequently, $\Delta$ is a clopen subset of $T\setd{\alpha}$.
    It follows that a path in $T\setd{\alpha}$ is contained in a single direction,
    hence the directions at $\alpha$ are exactly the path components of $T\setd{\alpha}$.
\end{proof}

\begin{proof}[Proof of \Cref{thm:shape}]
    Assume for the sake of contradiction that $M(\mu)$ is not a geodesic segment. 
    Then it contains a $3$-spider with center $c$ (a vertex of $T$) and outer endpoints $v_1,v_2,v_3$ 
    (they may not be vertices of $T$, and each segment $[c,v_i]$ contains no other vertices of $T$ except for $c$ and possibly $v_i$).
    Let $r = \min_{1\le i \le 3} d(c,v_i)$ (note that $r>0$) and for each $i\in \set{1,2,3}$ define $\beta_i\in [c,v_i]$ such that $d(c,\beta_i)=r$.
    Fix an arbitrary $x\in T$.
    In the case where $x\in \Delta_{c\to \beta_1}$, we observe that $d(\beta_1, x) = |d(c,x) - r|$,
    $d(\beta_2,x) = d(\beta_3,x) = d(c,x) + r$ and therefore
    $$\frac 13 \sum_{i=1}^{3} d(\beta_i, x) = \frac 13\big(|d(c,x) - r| + 2(d(c,x) + r)\big) \ge d(c,x) + \frac r3 > d(c,x).$$
    A similar result holds in the case where $x\in \Delta_{c\to \beta_2}\cup \Delta_{c\to \beta_3}$.
    In the case where $x\notin \bigcup_{i=1}^3 \Delta_{c\to \beta_i}$, we have
    $d(\beta_1,x) = d(\beta_2,x) = d(\beta_3,x) = d(c,x) + r$
    and $\frac 13 \sum_{i=1}^{3} d(\beta_i, x) > d(c,x)$.
    Consequently, for any $x\in T$, since $\ell$ is convex and strictly increasing, 
    $$\frac 13 \sum_{i=1}^{3} \ell(d(\beta_i, x)) \ge \ell\Big(\frac 13 \sum_{i=1}^{3} d(\beta_i, x)\Big) > \ell(d(c,x)).$$
    Let us define the function $\psi:x\mapsto - [\ell(d(c, x)) - \ell(d(o, x))] + \frac 13 \sum_{i=1}^{3} [\ell(d(\beta_i, x)) - \ell(d(o, x))]$,
    so that $\psi$ takes only positive values.
    Since $\set{c,\beta_1,\beta_2,\beta_3}\subset M(\mu)$, we have $0<\int_T \psi(x) \diff\mu(x) = -\phi(c) + \frac 13 \sum_{i=1}^{3} \phi(\beta_i) = 0$, which is a contradiction.
    Hence $M(\mu)$ is a geodesic segment.
\end{proof}

\begin{proof}[Proof of \Cref{prop:derivative}]
    1. Consider $v\in \Delta$ and
    let $\gamma:[0,1]\to [\alpha,v]$ be the geodesic from $\alpha$ to $v$.
    Define additionally the convex function $\Phi:[0,1]\to (-\infty, +\infty], t\mapsto \phi(\gamma_t)$, 
    and note that $\Phi(0)=\phi(\alpha)< +\infty$.
    Because $\Phi$ is convex and $0\in \dom \Phi$, it is easily seen that the mapping $t\mapsto (\Phi(t)-\Phi(0))/t$ is increasing on $(0,1]$,
    thus $\lim_{t\downarrow 0} (\Phi(t)-\Phi(0))/t$ exists in $[-\infty, +\infty]$,
    and so does $\lim\limits_{\beta \to \alpha, \beta \in \Delta} Q(\beta)$.

    2. We pick $\beta \in \Delta \cap F$ and we denote by $v$ the neighboring vertex of $\alpha$ contained in $\Delta$.
    Let us consider a sequence $(\beta_n)_{n\ge 1}$ in $\Delta$ such that $d(\beta_n,\alpha)\to 0$.
    By discarding finitely many terms if necessary, we may assume \wlo that 
    $d(\beta_n,\alpha)< \min(d(v,\alpha), d(\beta, \alpha))$ for each $n\ge 1$.
    In particular, $\beta_n$ belongs to $(\alpha,v)\cap (\alpha, \beta)$.
    We can then decompose the distance $d(\beta_n,x)$ with respect to the location of $x$ in the tree, and we obtain
    \begin{align}
        \label{eq:decomposition}
        \forall n\ge 1, \quad \frac{\phi(\beta_n) - \phi(\alpha)}{d(\beta_n, \alpha)}
        &= \begin{aligned}[t]
            &\phantom{+}\int_{T\setminus \Delta} \frac{\ell(d(\alpha,x) + d(\alpha,\beta_n))- \ell(d(\alpha,x))}{d(\beta_n, \alpha)} \diff\mu(x)
            \\&+ \int_{\Delta_{\beta_n \to v}} \frac{\ell(d(\alpha,x) - d(\alpha,\beta_n))- \ell(d(\alpha,x))}{d(\beta_n, \alpha)} \diff\mu(x)
            \\&+ \int_{(\alpha,\beta_n]}  \frac{\ell(d(\beta_n,x))- \ell(d(\alpha,x))}{d(\beta_n, \alpha)} \diff\mu(x).
        \end{aligned}
    \end{align}
    
    To deal with the first integral, let 
    $f_n:x\mapsto \indic{T\setminus \Delta}(x)\big[\ell\big(d(\alpha,x) + d(\alpha,\beta_n)\big)- \ell(d(\alpha,x))\big]/d(\beta_n,\alpha)$,
    so that $(f_n)_{n\ge 1}$ converges pointwise to $x\mapsto \indic{T\setminus \Delta}(x) \ell_+'(d(\alpha,x))$.
    By \eqref{eq:integrand-estimate}, we have the estimate
    \begin{align*}
        \forall x\in T,\quad |f_n(x)|
        &\le \indic{T\setminus \Delta}(x)
        \ell_+'\big(\max[d(\beta_n,x),d(\alpha,x)]\big)
        \\
        &= \indic{T\setminus \Delta}(x)
        \ell_+'(d(\beta_n,x))
        \le \ell_+'(d(\beta,x)).
    \end{align*}
    Indeed, if $x\in T\setminus\Delta$, then the geodesic
    $[\beta,x]$ goes through $\beta_n$ and $\alpha$, hence $d(\beta_n,x)\le d(\beta,x)$.
    Since $\beta\in F$, the dominated convergence theorem applies, so
    that 
    $$\int_{T\setminus \Delta} \frac{\ell(d(\alpha,x) + d(\alpha,\beta_n))- \ell(d(\alpha,x))}{d(\beta_n, \alpha)}\diff\mu(x)\xrightarrow[n\to\infty]{}
    \int_{T\setminus \Delta} \ell_+'(d(\alpha,x))\diff\mu(x).$$

    For the second integral, 
    let $f_n:x\mapsto \indic{\Delta_{\beta_n \to v}}(x)\big[\ell(d(\alpha,x) - d(\alpha,\beta_n))- \ell(d(\alpha,x))\big]/d(\beta_n,\alpha)$,
    so that $(f_n)_{n\ge 1}$ converges pointwise to $x\mapsto -\indic{\Delta}(x) \ell_-'(d(\alpha,x))$.
    We have the estimate
    \begin{align*}
        \forall x\in T,\quad |f_n(x)|
        &\le \indic{\Delta_{\beta_n \to v}}(x) \ell_+'\big(\max[d(\beta_n,x), d(\alpha,x)]\big)
        \\&= \indic{\Delta_{\beta_n \to v}}(x) \ell_+'(d(\alpha,x)) 
        \le \ell_+'(d(\alpha,x)).
    \end{align*}
    Since $\alpha \in F$, the dominated convergence theorem yields
    $$\int_{\Delta_{\beta_n \to v}} \frac{\ell(d(\alpha,x) - d(\alpha,\beta_n))- \ell(d(\alpha,x))}{d(\beta_n, \alpha)}\diff\mu(x)\xrightarrow[n\to\infty]{}
    -\int_{\Delta} \ell_-'(d(\alpha,x))\diff\mu(x).$$

    Finally, for the third integral,
    let $f_n:x\mapsto \indic{(\alpha,\beta_n]}(x)\big[\ell(d(\beta_n,x))- \ell(d(\alpha,x))\big]/d(\beta_n,\alpha)$,
    so that $(f_n)_{n\ge 1}$ converges pointwise to $0$.
    We have the estimate 
    $$\forall x\in T,\quad |f_n(x)|\le \indic{(\alpha,v]}(x) \ell_+'\big(\max[d(\beta_n,x), d(\alpha,x)]\big)\le \ell_+'(d(\alpha,v)),$$
    thus by the dominated convergence theorem we obtain 
    $$\int_{(\alpha,\beta_n]}  \frac{\ell(d(\beta_n,x))- \ell(d(\alpha,x))}{d(\beta_n,\alpha)}\diff\mu(x) \xrightarrow[n\to\infty]{} 0.$$

    3. For the first claim, write 
    $$\frac{\phi(\beta) - \phi(\alpha)}{d(\beta,\alpha)} 
    = \int_T \frac{\ell(d(\beta, x)) - \ell(d(\alpha, x))}{d(\beta, \alpha)} \diff \mu(x)$$
    and apply \eqref{eq:integrand-estimate} to obtain the lower bound
    $\frac{\phi(\beta) - \phi(\alpha)}{d(\beta,\alpha)} \ge -\int_T \ell_+'(d(\alpha,x)) \diff \mu(x) > -\infty$.

    Let us prove the second claim. We pick $v\in \Delta$ and we consider $\gamma:[0,1]\to [\alpha,v]$ the geodesic from $\alpha$ to $v$,
    as well as an arbitrary sequence of positive real numbers $(t_n)_{n\geq 1}$ such that $t_n \downarrow 0$. 
    We let 
    \begin{align*}
        f_n&:x \mapsto \frac{\ell(d(\gamma(t_n), x)) - \ell(d(\alpha, x))}{t_n d(\alpha, v)}, \\
        f&:x \mapsto \indic{T\setminus \Delta} \ell_+'(d(\alpha,x)) - \indic{\Delta} \ell_-'(d(\alpha,x)), 
    \end{align*}
    and by arguing as in the proof of the second item 
    we find that $(f_n)_{n\geq 1}$ converges pointwise to $f$.
    For each $x\in T$, the function $\beta \mapsto \ell(d(\beta, x))$ is convex,
    hence $f_n(x) \downarrow f(x)$ as $n\to \infty$.
    We observe that $(\phi(\gamma(t_n)) - \phi(\alpha)) / d(\gamma(t_n), \alpha) \to \phi_\Delta'(\alpha)$,
    and since this limit is finite, we may assume by discarding finitely many terms that $\phi(\gamma(t_1)) < +\infty$, 
    which entails that $f_1$ is integrable.
    Moreover, $\alpha$ being in $F$ yields integrability of $f$.
    By the monotonicity of $(f_n)_{n\geq 1}$,
    $0 \le f_n - f \le f_1 - f$, hence $|f_n| \le f_1 + |f| - f$ and the dominated convergence theorem finishes the proof.

\end{proof}

\begin{proof}[Proof of \Cref{prop:sign-derivative}]
    1. Since $\phi'_\Delta(\alpha)<0$, there exists 
    $\beta \in \Delta$ such that $\phi(\beta)<\phi(\alpha)$, hence $\alpha\notin M(\mu)$.
    For the sake of contradiction, assume the existence of $\alphas\in M(\mu)\setminus \Delta$.
    Let $\gamma:[0,1]\to [\beta,\alphas]$ be the geodesic from $\beta$ to $\alphas$.
    Since $\alphas \notin \Delta$, we have $[\alpha,\alphas]\cap [\alpha,\beta]=\set{\alpha}$.
    Consequently, the geodesic $\gamma$ goes through $\alpha$:
    the equality $\gamma_t=\alpha$ holds for some $t\in (0,1)$, thus $$\phi(\alpha)=\phi(\gamma_t)\le (1-t)\phi(\beta)+t\phi(\alphas)<(1-t)\phi(\alpha)+t\phi(\alphas),$$
    hence $\phi(\alpha)<\phi(\alphas)$, a contradiction.

    2. Since $\phi'_\Delta(\alpha)>0$, there exists $\beta\in \Delta$ such that 
    $\forall \beta' \in T$, $\beta'\in (\alpha,\beta] \implies \phi(\beta')>\phi(\alpha)$.
    For the sake of contradiction, assume the existence of $\alphas\in M(\mu)\cap \Delta$.
    Since $\alphas \in \Delta$ and $\beta\in \Delta$, we can pick
    $\beta' \in (\alpha, \alphas] \cap (\alpha, \beta]$.
    Let $\gamma:[0,1]\to [\alpha,\alphas]$ be the geodesic from $\alpha$ to $\alphas$.
    For some positive $t$, $\gamma_t=\beta'$, thus $$\phi(\alpha)<\phi(\beta') = \phi(\gamma_t)\le (1-t)\phi(\alpha)+t\phi(\alphas),$$
    hence $\phi(\alpha)<\phi(\alphas)$, a contradiction.

    3. We let $\alpha'\in T\setd{\alpha}$ be arbitrary and we show that $\phi(\alpha)\le \phi(\alpha')$.

    Consider first the case where $\alpha' \in \Delta$.
    Letting $\gamma:[0,1]\to [\alpha,\alpha']$ be the geodesic from $\alpha$ to $\alpha'$,
    $\phi(\gamma_t)\le (1-t)\phi(\alpha) + t\phi(\alpha')$, thus for each $t\in (0,1)$, 
    \begin{equation}
        \label{eq:ineq-derivative}
        \frac{\phi(\gamma_t)-\phi(\alpha)}{t}\le \phi(\alpha')-\phi(\alpha).
    \end{equation}
    Letting $t\downarrow 0$ we obtain $0 = \phi'_{\Delta}(\alpha) d(\alpha,\alpha')\le \phi(\alpha')-\phi(\alpha)$, 
    hence 
    $\phi(\alpha)\le \phi(\alpha')$.
    
    Next, we consider the case where $\alpha' \notin \Delta\cup \set{\alpha}$.
    Let $v$ denote the neighboring vertex of $\alpha$ that satisfies $v\in \Delta$.
    Let $\gamma:[0,1]\to [\alpha',v]$ be the geodesic from $\alpha'$ to $v$. It passes through $\alpha$:
    set $t_0=d(\alpha',\alpha)/d(\alpha',v)$ so that $\gamma_{t_0} = \alpha$.
    Define additionally the convex function $\Phi:[0,1]\to (-\infty,+\infty], t\mapsto \phi(\gamma_t)$, and note that its right derivative at $t_0$ satisfies
    $\Phi_+'(t_0) = \phi'_v(\alpha) d(\alpha',v) = 0$.
    Consequently, $t_0$ is a minimizer of $\Phi$, which entails $\phi(\alpha) = \Phi(t_0)\le \Phi(0)=\phi(\alpha')$.
\end{proof}

\begin{proof}[Proof of \Cref{prop:optimality}]
    1. Consider $\alpha\in M(\mu)$. Since $\phi$ is proper, $\alpha$ is in $\dom \phi$.
    For the sake of contradiction, if there exists a direction $\Delta$ with $\phi_\Delta'(\alpha) < 0$, then by
    \Cref{prop:sign-derivative} we have the inclusion $M(\mu)\subset \Delta$. Since $\Delta$ does not contain $\alpha$,
    we have reached a contradiction. 
    
    Conversely, if $\phi_\Delta'(\alpha)=0$ for some direction $\Delta$, then by \Cref{prop:sign-derivative} we have $\alpha \in M(\mu)$.
    Otherwise, $\phi_\Delta'(\alpha)>0$ along every direction, so that
    $M(\mu)\subset \bigcap_{\Delta\in \dir_\alpha} T\setminus \Delta$.
    Since this intersection equals $\set{\alpha}$ and since $M(\mu)$ is nonempty, we must have $M(\mu)=\set{\alpha}$.

    2. This is part of the proof of the preceding item.
    
\end{proof}

\begin{proof}[Proof of \Cref{prop:localization-support}]
    1. We show that $\conv(A)$ is closed.
    Consider a sequence $(\alpha_n)_{n\ge 1}$ of points in $\conv(A)$ that converges to some $\alpha\in T$.

    First, we deal with the case where $\conv(A)$ intersects two distinct directions $\Delta$ and $\Delta'$ in $\dir_\alpha$.
    Let us pick $\beta\in \Delta \cap \conv(A)$ and $\beta'\in \Delta' \cap \conv(A)$.
    Since $\beta' \notin \Delta$ we have $\alpha \in [\beta,\beta']$, and since $\conv(A)$ is convex we obtain $\alpha\in \conv(A)$.

    Next, we assume that $\conv(A)$ is contained in $\Delta_{\alpha\to v} \cup \set{\alpha}$, where $v$ is a neighboring vertex of $\alpha$.
    There exists $N\ge 1$ such that $\alpha_n$ lies in $[\alpha,v]$ for each $n\ge N$.
    Assume for the sake of contradiction that $\alpha$ is not in $\conv(A)$.
    Then $\alpha$ is not in $A$, which is closed. We can then find $\beta \in [\alpha_n,\alpha_{n+1}]$ for some $n\ge N$, 
    such that $[\alpha,\beta]\subset (T\setminus A)$.
    Since $\alpha_n$ and $\alpha_{n+1}$ are in $\conv(A)$, it must be that $\beta\in \conv(A)$.
    Then $\conv(A)\setminus [\alpha,\beta]$ is strictly smaller than $\conv(A)$, 
    while containing $A$.
    By our assumption that $\conv(A) \subset \Delta_{\alpha\to v} \cup \set{\alpha}$,
    the subset $\conv(A)\setminus [\alpha,\beta]$ is also convex, which contradicts the minimality of $\conv(A)$.
    We have reached a contradiction, hence $\alpha \in \conv(A)$.

    2. Since $\supp \mu$ is closed, it follows that $\conv(\supp \mu)$ is closed.
    Next, let us prove that $M(\mu)\subset \conv(\supp \mu)$.
    We fix $\alpha\in T\setminus \conv(\supp \mu)$ and we show that $\alpha \in T\setminus M(\mu)$.
    If $\alpha \notin \dom \phi$, then the proof is finished, hence we assume next that $\phi(\alpha) < +\infty$.
    Since $T$ is a Hadamard space and $\conv(\supp \mu)$ is closed and convex, the metric projection on $\conv(\supp \mu)$ is well-defined \cite[Theorem~2.1.12]{bacak2014convex}, 
    and we denote it by $\pi$. Let $x\in \conv(\supp \mu)$.
    By the Pythagorean inequality \cite[Theorem~2.1.12(ii)]{bacak2014convex} and the strict monotonicity of $\ell$,
    \begin{equation}
        \label{eq:ineq-proj}
        \ell(d(\alpha,x)) > \ell(d(\pi(\alpha),x)).
    \end{equation}
    In particular, $\pi(\alpha) \in \dom \phi$.
    Since $\supp \mu\subset \conv(\supp \mu)$, we have further
    $$\phi(\alpha) - \phi(\pi(\alpha)) = \int_{\conv(\supp \mu)} \ell(d(\alpha,x)) - \ell(d(\pi(\alpha),x)) \diff \mu(x) > 0,$$
    where the last inequality follows from \eqref{eq:ineq-proj} and $\mu(\conv(\supp \mu))=1$.
    As a result, $\alpha \notin M(\mu)$.

    3. If $\supp \mu \subset G$ and $G$ is convex, then $\conv(\supp \mu) \subset G$ and therefore $M(\mu)\subset G$.

\end{proof}

\subsection{Proofs for Section~\ref{sec:stickiness}}

\begin{proof}[Proof of \Cref{thm:localization-sticky}]
    1. This follows from the second item of \Cref{prop:optimality}.

    2. a) \Cref{prop:optimality} yields $\alpha\in M(\mu)$.
    The implication $i\notin I_0 \implies \phi_{\Delta_i}'(\alpha) > 0$ and the equality 
    $$\bigcap_{i\notin I_0} (T\setminus \Delta_i) = \set{\alpha} \cup \bigcup_{i\in I_0} \Delta_i$$
    justify the inclusion $M(\mu) \subset \set{\alpha} \cup \bigcup_{i\in I_0} \Delta_i$.

    b) Suppose for the sake of contradiction that $I_0$ has more than two elements.
    Then \wlo $\phi'_{\Delta_i}(\alpha)=0$ for each $i\in [3]$.
    Fix $i\in [3]$ for now.
    By the definition of $\phi'_{\Delta_i}(\alpha)$, there exists a sequence $(\beta_n)_{n\ge 1}$ of points in $\Delta_i$
    such that $\beta_n\to \alpha$, $\phi(\beta_n)<+\infty$ and $Q(\beta_n)\to 0$.
    For each $n\ge 1$, let $f_n: T\to \R, x\mapsto \frac{\ell(d(\beta_n, x)) - \ell(d(\alpha, x))}{d(\beta_n, \alpha)}$.
    Since $\set{\beta_n,\alpha}\subset \dom \phi$, 
    we obtain that $f_n$ is integrable and $\int_T f_n(x)\diff\mu(x) = Q(\beta_n)$.
    We also define $g:T\to [0,\infty)$, $x\mapsto \indic{x\neq \alpha} \ell_-'(d(\alpha,x))$.
    For any $x\in T\setminus(\Delta_i \cup \set{\alpha})$, the convexity of $\ell$ yields
    $$\ell(d(\beta_n, x)) - \ell(d(\alpha, x)) 
    \ge \ell_-'(d(\alpha, x))\big(d(\beta_n, x)-d(\alpha, x)\big) 
    = \ell_-'(d(\alpha, x))d(\beta_n, \alpha),$$
    hence $f_n(x)\ge g(x)$.
    For any $x\in \Delta_i$, by the triangle inequality,
    $$\ell(d(\beta_n, x)) - \ell(d(\alpha, x)) 
    \ge \ell_-'(d(\alpha, x))\big(d(\beta_n, x)-d(\alpha, x)\big) 
    \ge -\ell_-'(d(\alpha, x))d(\beta_n, \alpha),$$
    hence $f_n(x)\ge -g(x)$.
    Since $f_n(\alpha)\ge 0$, we have shown that for every $x\in T$,
    $f_n(x) \ge g(x) \indic{T\setminus(\Delta_i \cup \set{\alpha})}(x) - g(x)\indic{\Delta_i}(x)$,
    or equivalently $g(x) \indic{\Delta_i}(x) \ge \frac 12 g(x) \indic{T\setd{\alpha}}(x) - \frac 12 f_n(x)$.
    Since $g\ge 0$ and $f_n$ is integrable, we may integrate this inequality and let $n\to \infty$ to obtain that
    $\int_{\Delta_i} g(x) \diff\mu(x) \ge \frac 12 \int_{T\setd{\alpha}} g(x) \diff\mu(x)$.

    Next, let us show that $\int_{T} g(x) \diff\mu(x) < +\infty$. Since $f_n$ is integrable we already have
    $$\int_{T\setminus(\Delta_i \cup \set{\alpha})} g(x) \diff\mu(x) < +\infty.$$
    Pick a $j\in \set{1,2,3}\setd{i}$ and note that for the same reasons 
    $\int_{T\setminus(\Delta_j \cup \set{\alpha})} g(x) \diff\mu(x) < +\infty$.
    The inclusion $\Delta_i \subset T\setminus(\Delta_j \cup \set{\alpha})$ yields
    $\int_{\Delta_i} g(x) \diff\mu(x) < +\infty$, hence the claim.

    Summing over $i\in \set{1,2,3}$ yields 
    $$\int_{T\setd{\alpha}} g(x) \diff\mu(x) 
    \ge \sum_{i=1}^{3}\int_{\Delta_i} g(x) \diff\mu(x) 
    \ge \frac 32 \int_{T\setd{\alpha}} g(x) \diff\mu(x).$$
    Combined with $\int_{T} g(x) \diff\mu(x) < +\infty$, we obtain that
    $\int_{T\setd{\alpha}} g(x) \diff\mu(x) = 0$.
    However, since $\ell_-'$ takes only positive values, $g$ is positive on $T\setd{\alpha}$,
    hence $\mu(\set{\alpha})=1$,
    which is a contradiction.

    c) Suppose \wlo that $\phi_{\Delta_1}'(\alpha) = \phi_{\Delta_2}'(\alpha) = 0$.
    With the same notation as in the proof of b), we obtain that 
    $\int_{T\setd{\alpha}} g(x) \diff\mu(x) 
    \ge \sum_{i=1}^{2}\int_{\Delta_i} g(x) \diff\mu(x) 
    \ge  \int_{T\setd{\alpha}} g(x) \diff\mu(x)$,
    therefore $\int_{\Delta_1} g(x) \diff\mu(x) = \int_{\Delta_2} g(x) \diff\mu(x) = \frac 12 \int_{T\setd{\alpha}} g(x)\diff\mu(x)$.
    Combined with $\int_{T} g(x) \diff\mu(x) < +\infty$, we obtain that $$\int_{T\setminus(\Delta_1\cup\Delta_2\cup\set{\alpha})} g(x) \diff\mu(x) = 0.$$
    The positivity of $g$ on $T\setminus(\Delta_1\cup\Delta_2\cup\set{\alpha})$
    then forces
    $\mu(\Delta_{1}\cup \Delta_{2}\cup \set{\alpha})=1.$

    3. This is part of \Cref{prop:sign-derivative}.
\end{proof}

\begin{proof}[Proof of \Cref{prop:sticky-geometry}]
    1. For each $i \in [m]$, we define $A_i= \int_{\Delta_i} \ell_+'(d(\alpha,x)) \diff \mu(x)$,
    as well as $A = \sum_{i=1}^m A_i = \int_{T\setminus\set{\alpha}} \ell_+'(d(\alpha,x)) \diff \mu(x)$.

    Fix $i\in[m]$. 
    By the differentiability assumption, note that $A_i$ coincides with $\int_{\Delta_i} \ell_-'(d(\alpha,x)) \diff \mu(x)$.
    Since $\alpha\in\operatorname{int}(F)$, we have $\alpha\in F$ and $\Delta_i \cap F \neq \emptyset$.
    Consequently, the closed form \eqref{eq:derivative} is available in every direction.
    Since $T\setminus\Delta_i = \set{\alpha} \cup \bigcup_{j\neq i} \Delta_j$,
    it follows that 
    $$\phi'_{\Delta_i}(\alpha) = \ell_+'(0)\mu(\set{\alpha}) + \sum_{j\neq i} A_j - A_i = \ell_+'(0)\mu(\set{\alpha}) + A - 2A_i.$$
    Summing over $i \in [m]$ proves the claimed equality.

    2. Suppose that $\alpha$ is sticky and $\ell_+'(0)\mu(\set{\alpha})=0$. 
    Then $0 < \sum_{i=1}^{m} \phi_{\Delta_i}'(\alpha) = (m-2)A$.
    Since $A\ge 0$, it follows that $m-2>0$, hence $m\ge 3$ and $\alpha$ is a branching vertex.

    3. Suppose that $m=2$ and $\alpha\in M(\mu)$.
    By \Cref{prop:optimality},
    $\phi'_{\Delta_1}(\alpha)$ and $\phi'_{\Delta_2}(\alpha)$ are
    nonnegative, while the first item of \Cref{prop:sticky-geometry} gives
    $\phi'_{\Delta_1}(\alpha)+\phi'_{\Delta_2}(\alpha) = 2\ell_+'(0)\mu(\set{\alpha})$.
    The equivalence follows immediately.
\end{proof}

\begin{lemma}
    \label{lemma:interior-F}
    Let $\alpha\in T$.
    \begin{enumerate}
        \item If $\alpha$ is a leaf, the following implication holds: $\alpha \in F  \implies \alpha \in \operatorname{int}(F)$.
        \item If $\alpha$ is not a leaf, $\alpha \in \operatorname{int}(F) \iff \exists \rho  \in (0,+\infty), \E{\ell_+'(d(\alpha, X) + \rho)} < +\infty$.
        \item $\alpha \in \operatorname{int}(F)$ if and only if $\alpha \in F$ and $\forall \Delta \in \dir_\alpha$, $\Delta \cap F \neq \emptyset$.
    \end{enumerate}
\end{lemma}

\begin{proof}[Proof of \Cref{lemma:interior-F}]
    1. Suppose that $\alpha$ is in $F$. We denote by $v$ the neighboring vertex of $\alpha$ and we show the inclusion $B(\alpha, d(\alpha, v)) \subset F$.
    Fix an arbitrary $\beta \in B(\alpha, d(\alpha, v))$ and consider $x\in T$.
    If $x\in [\alpha, \beta]$, then $d(\beta, x) = d(\alpha, \beta) - d(\alpha,x) \le d(\alpha, \beta) < d(\alpha, v)$.
    Otherwise, $d(\beta, x) = d(\alpha,x) - d(\alpha, \beta) \le d(\alpha, x)$.
    Since $\ell_+'$ is nonnegative and increasing, it follows that 
    $\ell_+'(d(\beta, x)) \le \ell_+'(d(\alpha, v)) + \ell_+'(d(\alpha, x))$ for every $x\in T$, hence $\beta \in F$.

    2. If there exists $\rho>0$ such that $\E{\ell_+'(d(\alpha, X) + \rho)} < +\infty$, then by the triangle inequality and the monotonicity of $\ell_+'$,
    it follows easily that $\alpha \in \operatorname{int}(F)$.
    Let us prove the converse. We assume that $B(\alpha, \varepsilon)\subset F$ for some $\varepsilon>0$.
    Since $\alpha$ is not a leaf, we may find two directions $\Delta \neq \Delta' \in \dir_\alpha$, as well as $\beta \in \Delta\cap B(\alpha, \varepsilon)$
    and $\beta' \in \Delta' \cap B(\alpha, \varepsilon)$.
    We let $\rho = \min(d(\alpha, \beta), d(\alpha, \beta'))$ and we consider an arbitrary $x\in T$.
    In the case where $x\in \Delta$, we have $d(\beta', x) = d(\beta', \alpha) + d(\alpha, x) \ge \rho + d(\alpha, x)$.
    If $x\notin \Delta$, $d(\beta, x) = d(\beta, \alpha) + d(\alpha, x) \ge \rho + d(\alpha, x)$.
    Consequently, $\ell_+'(d(\alpha, x) + \rho) \le \ell_+'(d(\beta, x)) + \ell_+'(d(\beta', x))$ for every $x\in T$.
    Since $\beta, \beta'\in F$, the conclusion follows.

    3. A ball centered at $\alpha$ intersects all the directions at $\alpha$, 
    and the forward implication follows immediately.
    Conversely, we suppose that $\alpha \in F$ and $\forall \Delta \in \dir_\alpha$, $\Delta \cap F \neq \emptyset$.
    If $\alpha$ is a leaf, then $\alpha\in \operatorname{int}(F)$ by the first item.
    We consider the case where $\alpha$ is not a leaf.
    We may find two directions $\Delta \neq \Delta' \in \dir_\alpha$, 
    as well as $\beta \in \Delta\cap F$ and $\beta' \in \Delta' \cap F$.
    As in the proof of the second item, we let $\rho = \min(d(\alpha, \beta), d(\alpha, \beta'))$
    and we obtain that $\ell_+'(d(\alpha, x) + \rho) \le \ell_+'(d(\beta, x)) + \ell_+'(d(\beta', x))$ for every $x\in T$,
    Integration over $x$ yields $\E{\ell_+'(d(\alpha, X) + \rho)} < +\infty$, hence $\alpha\in \operatorname{int}(F)$.
\end{proof}

\begin{proof}[Proof of \Cref{prop:stickiness-R}]
    1. Let us verify that \Cref{prop:sticky-geometry} is applicable.
    The equality \eqref{eq:differentiability} clearly holds under condition (b).
    We now justify why it also holds under condition (a).
    We let $D = \set{z\in (0,+\infty): \ell_-'(z) \neq \ell_+'(z)}$. 
    By \Cref{lemma:integral-rep}, $D$ coincides with $\set{z\in (0,+\infty): \ell_+'(z-) \neq \ell_+'(z+)}$
    \ie $D$ is the set of discontinuity points of the increasing function $\ell_+'$,
    hence $D$ is countable.
    We observe that
    \begin{equation}
        \label{eq:union-spheres}
        \set{x\in T\setminus\set{\alpha}: \ell_-'(d(\alpha,x)) \neq \ell_+'(d(\alpha,x))}
        = \bigcup_{z\in D} \set{x\in T: d(\alpha, x) = z},
    \end{equation}
    which is a countable union of spheres.
    By \Cref{lemma:ball}, each sphere is finite, hence the left-hand side of \eqref{eq:union-spheres} is countable.
    Since $\mu$ is atomless, this set has $\mu$-measure $0$.

    Under our assumptions, we have $\ell_+'(0)\mu(\set{\alpha}) = 0$.
    Let us consider $\alpha\in M(\mu)$, with $\dir_\alpha = \set{\Delta_1,\ldots,\Delta_m}$.
    By \Cref{prop:optimality}, all the directional derivatives at $\alpha$ are nonnegative.
    Combined with \Cref{prop:sticky-geometry}, we obtain that 
    $0 \le \sum_{i=1}^{m} \phi_{\Delta_i}'(\alpha) = (m-2) \int_{T\setminus\set{\alpha}} \ell_+'(d(\alpha,x)) \diff \mu(x)$.
    By our assumption on $T$, we have $m\in \set{1,2}$, hence $(m-2) \int_{T\setminus\set{\alpha}} \ell_+'(d(\alpha,x)) \diff \mu(x) \le 0$,
    $\sum_{i=1}^{m} \phi_{\Delta_i}'(\alpha) = 0$
    and finally $\phi_{\Delta_i}'(\alpha) = 0$ for each $i\in [m]$.

    2. Since $\alpha$ is a leaf, $m=1$ and $0 \le \phi_{\Delta_1}'(\alpha) = -\int_{T\setminus\set{\alpha}} \ell_+'(d(\alpha,x)) \diff \mu(x)$.
    Since $\ell_+'$ is nonnegative, we must have $\mu = \delta_\alpha$, which invalidates condition (a) and forces (b), thus $\ell_+'(0)=0$
    and $\alpha$ is one-sided partly sticky.

    3. Since $\alpha$ is not a leaf and $m\in \set{1,2}$, it holds that $m=2$ and the claim follows immediately from the third item of \Cref{prop:sticky-geometry}.
\end{proof}

\begin{proof}[Proof of \Cref{thm:robust-sticky}]
    In this proof we will say for convenience that \textit{$\alpha$ is $\mu$-sticky} if $\alpha$ is sticky when the measure under study is $\mu$,
    and we say similarly that \textit{$\alpha$ is $\mu$-partly sticky}.
    We denote by $v_1,\ldots,v_m$ the neighboring vertices of $\alpha$.

    1. Suppose $\alpha$ is $\mu$-sticky and define  $L = \ell_+'(\diam(T))$,
    with the abuse of notation that $L = \lim_{z\to\infty}\ell_+'(z) < +\infty$ when $T$ is unbounded. 
    Let also $\varepsilon = \min_j \phi_{v_j}'(\alpha)/(4L)$.
    Let $\nu$ be such that $\tv(\nu, \mu)\le \varepsilon$, let $\tphi$ denote the objective function associated to $\nu$ and fix some $i \in [m]$. 
    Note that \Cref{assum:bounded} holds, hence $\tphi$ is well-defined with $\dom \tphi = T$.
    Moreover, \ref{assum:bounded} also makes the second item of \Cref{prop:derivative} applicable.
    Classically,
    $$\tv(\nu_1,\nu_2) = \frac 12 \sup\set[\Big]{\int_T f(x) \diff\nu_1(x)-\int_T f(x) \diff\nu_2(x)\mid  f:T\to [-1,1] \text{ is measurable}}$$
    (see, \eg \cite[Lemma~1 p.~432]{shiryaev2016proba}).
    Since the function 
    \begin{equation}
        \label{eq:varphi}
        \varphi: x\mapsto \indic{T\setminus \Delta_{\alpha \to v_i}}(x) \ell_+'(d(\alpha,x))
        - \indic{\Delta_{\alpha \to v_i}}(x) \ell_-'(d(\alpha,x))
    \end{equation}
    is bounded by $L$,
    we have the estimate
    $$|\tphi_{v_i}'(\alpha) - \phi_{v_i}'(\alpha)| = \left|\int_T \varphi(x) \diff\nu(x)-\int_T \varphi(x) \diff\mu(x)  \right| 
    \le 2 L \tv(\nu, \mu) \le \frac{\phi_{v_i}'(\alpha)}{2},$$
    which implies $\tphi_{v_i}'(\alpha)>0$, thus $\alpha$ is $\nu$-sticky and $M(\nu) = \set{\alpha}$.

    Conversely, for some $\varepsilon>0$ we suppose that $\forall \nu, \tv(\nu, \mu)\le \varepsilon \implies M(\nu) = \set \alpha$
    and we assume for the sake of contradiction that $\alpha$ is not $\mu$-sticky.
    We may assume \wlo that $\varepsilon \in (0,1]$.
    Since $\alpha\in M(\mu)$, $\alpha$ is $\mu$-partly sticky and there exists $i$ with $\phi_{v_i}'(\alpha)=0$, 
    \ie $\int_T \varphi(x) \diff\mu(x) = 0$ with $\varphi$ defined in \eqref{eq:varphi}.
    Since $\ell$ is strictly increasing, $\ell_-'(d(\alpha,v_i))>0$ hence $\varphi(v_i)<0$.
    Next, define the mixture measure $\nu = (1-\varepsilon)\mu + \varepsilon \delta_{v_i}$ so that 
    \begin{equation}
        \label{eq:proof-robust}
        \tphi_{v_i}'(\alpha) = (1-\varepsilon)\phi_{v_i}'(\alpha) +  \varepsilon \varphi(v_i) = \varepsilon \varphi(v_i)<0
    \end{equation}
    and $\tv(\nu, \mu)\le \varepsilon$.
    By our initial assumption, the closeness of the measures implies $M(\nu) = \set \alpha$, which contradicts the inequality \eqref{eq:proof-robust}.

    2. Suppose $\alpha$ is $\mu$-sticky. We reuse the notations from the proof of the first item.
    Say for concreteness that $\ell'$ is $L$-Lipschitz on $I$.
    Let us show that for any $\nu\in \mathcal P_1$, the objective $\tphi$ is well-defined.
    Note that $\ell'(z) = |\ell'(z) - \ell'(0)|\le Lz$ for each $z\in I$, hence
    $\int_T \ell'(d(\beta, x))\diff\nu(x) \le L \int_T d(\beta, x)\diff\nu(x) < +\infty$ for every $\beta\in T$
    and consequently $\tphi$ is well-defined with $\dom \tphi = T$.

    Next, we prove that $\varphi$ is $L$-Lipschitz.
    In the case where $\set{x,y}\subset \Delta_{\alpha\to v_i}$ or  $\set{x,y}\subset T\setminus \Delta_{\alpha\to v_i}$,
    we have $|\varphi(x) - \varphi(y)| = |\ell'(d(\alpha,x)) - \ell'(d(\alpha,y))| \le L d(x,y)$.
    In the case where $x\in \Delta_{\alpha\to v_i}$ and $y\notin \Delta_{\alpha\to v_i}$, 
    $|\varphi(x) - \varphi(y)| = \ell'(d(\alpha,x)) + \ell'(d(\alpha,y)) \le L(d(\alpha,x) + d(\alpha,y)) = Ld(x,y)$.

    Let $\varepsilon = \min_j \phi_{v_j}'(\alpha)/(2L)$ and let $\nu\in \mathcal P_1$ be such that $W_1(\nu, \mu)\le \varepsilon$.
    By the duality formula for $W_1$,
    $$|\tphi_{v_i}'(\alpha) - \phi_{v_i}'(\alpha)| = \left|\int_T \varphi(x) \diff\nu(x)-\int_T \varphi(x) \diff\mu(x)  \right| 
    \le L W_1(\nu, \mu) \le \frac{\phi_{v_i}'(\alpha)}{2},$$
    hence the claim.

    For the converse, we leverage a similar argument as in the setting of total variation.
    Since $M(\mu) = \set{\alpha}$ and $v_i \neq \alpha$, we have $\mu \neq \delta_{v_i}$ hence $\int_T d(v_i, x)\diff\mu(x) > 0$.
    We let $\varepsilon' = \min(\varepsilon /\int_T d(v_i, x)\diff\mu(x),1)$
    and the mixture measure is $\nu = (1-\varepsilon')\mu + \varepsilon' \delta_{v_i}$.
    We show next that $W_1(\nu, \mu) \le \varepsilon$.
    Let us denote by $X$ a random element having distribution $\mu$, by $\pi_1$ (resp., $\pi_2$) the distribution of $(X,X)$ (resp., $(v_i, X)$)
    and we define the coupling $\pi = (1-\varepsilon')\pi_1 + \varepsilon' \pi_2$.
    By the definition of $W_1$ it holds that $W_1(\nu, \mu) \le \int_{T^2} d(x,y) \diff \pi(x,y) 
    = \varepsilon' \int_T d(v_i, x)\diff \mu(x)\le \varepsilon$.
    The rest of the proof is identical to the previous one.
\end{proof}

\subsection{Proofs for Section~\ref{sec:stats}}

\begin{lemma}
    \label{lemma:Borel-selection}
    Let $n\ge 1$. 
    There exists a Borel measurable $f_n:T^n \to T$
    such that for each $(x_1,\ldots,x_n) \in T^n$,
    $f_n(x_1,\ldots,x_n) \in \argmin_{\alpha\in T} \sum_{k=1}^{n} \ell(d(\alpha, x_k))$.
    Consequently, the random element $f_n(X_1,\ldots,X_n)$ is a measurable selection from the subset $M(\hmu_n)$.
\end{lemma}

\begin{proof}
    We define the function 
    \begin{align*}
        g \colon \quad \quad \quad T^n \times T &\to \R \\
        (x_1,\ldots, x_n, \alpha) &\mapsto \sum_{k=1}^{n} \ell(d(\alpha, x_k))
    \end{align*}
    and we show that Corollary~1(ii) in \cite{brown1973measurable} applies.
    Let us fix $(x_1,\ldots,x_n) \in T^n$.
    By \Cref{prop:phi-trees} applied to the measure $\mu = \frac 1n \sum_{k=1}^{n} \delta_{x_k}$, 
    the set denoted by $I$ in \cite{brown1973measurable} equals $T^n$.
    Moreover, $T$ is $\sigma$-compact (\ie it is the countable union of compact subsets)
    and the mapping $\alpha \mapsto g(x_1,\ldots,x_n,\alpha)$ is continuous on $T$.
    Consequently, \cite[Corollary~1(ii)]{brown1973measurable} applies and the claim follows.
\end{proof}

\begin{proof}[Proof of Theorem \ref{thm:lln-sets}] 
    We apply \cite[Theorem~2.3]{artstein1995consistency} to derive that $\hphi_n$ epi-converges to $\phi$ with probability $1$.
    Let us verify that the assumptions needed in \cite{artstein1995consistency} hold in our setting.
    The loss $\ell$ is continuous by \Cref{lemma:integral-rep}, hence their Assumption 2.1 is verified.
    Let us check their Assumption 2.2. 
    Pick $\alpha_0\in T$ and note that 
    for every $x\in T$ and every $\alpha$ in $B(\alpha_0, 1)$,
    we have by \eqref{eq:integrand-estimate} that $\ell(d(\alpha,x)) - \ell(d(o, x)) \ge -(1+d(\alpha_0, o))\ell_+'(d(o,x))$,
    and since $o\in F$ we validate their Assumption 2.2.
    Applying \cite[Theorem~2.3]{schotz2022strong} (with a minor adjustment to account for the fact that $\phi$ may take the value $+\infty$) 
    then yields the strong consistency in outer limit of $(M(\hmu_n))_{n\ge 1}$.

    Let us exhibit some $R>0$ such that $\P[\big]{\exists N\ge 1, \forall n\ge N, M(\hmu_n) \subset \overline B(o, R)} = 1$.
    By \Cref{prop:phi-trees}, $\phi$ is coercive, hence there exists $R>0$ such that for every $\alpha\in T$, $d(\alpha, o)\ge R \implies \phi(\alpha) >0$.
    Let $\S = \set{\alpha\in T: d(\alpha, o)= R}$.
    \Cref{lemma:ball} implies that $\S$ is finite.
    Define the event $\Omega' = \bigcap_{\beta \in \S} \set{\hphi_n(\beta)\to \phi(\beta)}$.
    The strong law of large numbers (a variant thereof in the case where $\phi(\beta)=+\infty$, see \cite[Remark 2 p.~18]{shiryaev2019proba}),
    yields $\P{\Omega'}=1$.
    Define $\Omega'' =
    \bigcup_{N\ge 1}\bigcap_{n\ge N} \bigcap_{\beta \in \S} \set{\hphi_n(\beta)> 0}$
    and note the inclusion $\Omega'\subset \Omega''$.
    Pick any $\alpha$ such that $d(\alpha, o)> R$, let $\gamma:[0,1]\to [o,\alpha]$ denote the geodesic from $o$ to $\alpha$ and set $t = R/d(\alpha,o)\in (0,1)$ so that $\gamma_t\in \S$.
    Since $\hphi_n$ is convex, it holds on $\Omega''$ that 
    for all sufficiently large $n$,
    $$\hphi_n(o) = 0 < \hphi_n(\gamma_t) \le (1-t)\hphi_n(o) + t\hphi_n(\alpha) = t\hphi_n(\alpha) \le\hphi_n(\alpha).$$
    We have shown $\Omega''\subset \set*{\exists N\ge 1, \forall n\ge N, M(\hmu_n) \subset \overline B(o, R)}$, hence the claim.

    It follows that with probability $1$, 
    $(M(\hmu_n))_{n\ge 1}$ is eventually precompact, where the terminology is from \cite{schotz2022strong}.
    By \cite[Theorem~2.4(ii)]{schotz2022strong}, $(M(\hmu_n))_{n\ge 1}$ is strongly consistent in one-sided Hausdorff distance.
\end{proof}

\begin{proof}[Proof of \Cref{thm:sticky-lln}]
    1. Let $i\in [m]$ be fixed. 
    For each $k\in [n]$ we define the random variable $$Y_k = \indic{T\setminus \Delta_{i}}(X_k) \ell_+'(d(\alpha,X_k))
    - \indic{\Delta_{i}}(X_k) \ell_-'(d(\alpha,X_k)),$$
    so that $\hphi_{\Delta_i}'(\alpha) = \frac 1n \sum_{k=1}^{n} Y_k$. 
    The $Y_k$ are \iid and integrable, thus by the strong law of large numbers $\hphi_{\Delta_i}'(\alpha)$ converges almost surely to $\E{Y_1}$. 
    Our assumption that $\alpha \in \operatorname{int}(F)$ entails both $\alpha\in F$ and $\Delta_i \cap F \neq \emptyset$, hence 
    $\E{Y_1} = \phi_{\Delta_i}'(\alpha)$.
    Besides, since each $\phi_{\Delta_i}'(\alpha)$ is positive, the following inclusion holds:
    $$\bigcap_{i=1}^m \set[\Big]{\hphi_{\Delta_i}'(\alpha) \underset{n}{\to} \phi_{\Delta_i}'(\alpha)} \, \subset\,  \liminf_{n\ge 1} \Big(\bigcap_{i=1}^m \set*{\hphi_{\Delta_i}'(\alpha) > 0}\Big).$$
    Consequently, $\P[\big]{\exists N\ge 1, \forall n\ge N, M(\hmu_n)=\set{\alpha}}=1.$

    2.a) The proof is similar to that of the first item.

    2.b) By \Cref{thm:localization-sticky} and since $\set{\alpha}\cup \bigcup_{i\in I_0} \Delta_{i}$ is closed, 
    we obtain the inclusion $\supp \mu \subset \set{\alpha}\cup \bigcup_{i\in I_0} \Delta_{i}$.
    Moreover, with probability $1$ we have $\supp \hmu_n \subset \supp \mu$.
    Since $\set{\alpha}\cup \bigcup_{i\in I_0} \Delta_{i}$ is also convex, 
    applying \Cref{prop:localization-support} with the measure $\hmu_n$ yields $\P{M(\hmu_n)\subset \set{\alpha}\cup \bigcup_{i\in I_0} \Delta_{i}}=1$.
    
    3. The proof is similar to that of the first item.
\end{proof}

\begin{proof}[Proof of \Cref{prop:converse-empirical-stickiness}]
    Let us consider the case where $\mu\neq \delta_\alpha$ and let us assume for contradiction that $\phi_\Delta'(\alpha)\le 0$ for some $\Delta\in \dir_\alpha$.
    For each $k,n\ge 1$ we define the random variables $$Y_k = \indic{T\setminus \Delta}(X_k) \ell_+'(d(\alpha,X_k))
    - \indic{\Delta}(X_k) \ell_-'(d(\alpha,X_k)) \text{ and } S_n = \sum_{k=1}^{n} Y_k,$$ 
    so that $\hphi_{\Delta}'(\alpha) = \frac 1n S_n$. 
    Note that $Y_1$ is integrable since $|Y_1|\le \ell_+'(d(\alpha,X_1))$ and $\alpha\in F$.
    
    Let us show that $\E{Y_1}\le 0$.
    We recall that $\phi_\Delta'(\alpha)=\lim\limits_{\beta \to \alpha, \beta \in \Delta} \int_T \frac{\ell(d(\beta, x))-\ell(d(\alpha, x))}{d(\beta, \alpha)} \diff\mu(x)$.
    Let $(\beta_p)_{p\ge 1}$ be a sequence in $\Delta$ such that $d(\beta_p,\alpha)\to 0$
    and define $f_p:x\mapsto \big(\ell(d(\beta_p, x))-\ell(d(\alpha, x))\big) / d(\beta_p, \alpha)$.
    It follows from \eqref{eq:integrand-estimate} that $f_p(x) \ge -\ell_+'(d(\alpha, x))$ for each $x\in T$.
    By considering cases according to the location of $x$ in the tree and decomposing distances as in \eqref{eq:decomposition},
    we obtain that $(f_p)_{p\ge 1}$ converges pointwise to $x\mapsto \indic{T\setminus \Delta}(x) \ell_+'(d(\alpha,x)) - \indic{\Delta}(x) \ell_-'(d(\alpha,x))$.
    Since $\alpha\in F$, we may apply Fatou's lemma for functions with an integrable lower bound and obtain that 
    $\E{Y_1}=\int_T \liminf_p f_p(x)\diff \mu(x) \le \liminf_p \int_T f_p(x)\diff \mu(x) = \phi_\Delta'(\alpha) \le 0$. 
    
    Whenever $\alpha$ is a minimizer of $\hphi_n$ we have $\hphi_{\Delta}'(\alpha)\ge 0$ hence $S_n\ge 0$, and it follows that 
    \begin{equation}
        \label{eq:RW-nonnegative}
        \P[\big]{\exists N\ge 1, \forall n\ge N, S_n\ge 0}=1.
    \end{equation}
    In the case where $\E{Y_1}< 0$, the strong law of large numbers guarantees that $\frac 1n S_n \xrightarrow[]{a.s.} \E{Y_1}$,
    hence $\P{\exists N\ge 1, \forall n\ge N, S_n< 0}=1$, which contradicts \eqref{eq:RW-nonnegative}.
    We assume next that $\E{Y_1}= 0$.
    Since $\ell_-'$ takes only positive values and $\ell_+'(z)=0$ implies $z=0$,
    the following inclusion of events holds: $\set{Y_1=0} \subset \set{X_1 = \alpha}$.
    By our assumption that $\mu\neq \delta_\alpha$, it follows that $\P{Y_1=0}<1$.
    A classical result on random walks with mean-zero increments (see, \eg \cite[Corollary~2 p.\ 154]{chow2003probability})
    states that $\P{\liminf_n S_n = -\infty}=1$, which contradicts \eqref{eq:RW-nonnegative}.
\end{proof}

\begin{proof}[Proof of \Cref{thm:sticky-nonasymp-subgaussian}]
    We only consider the sticky case as the other cases are dealt with in a similar fashion.
    Let $i\in [m]$ be fixed. 
    For each $k\in [n]$ we define the random variable $Y_k$ as 
    $$Y_k = \indic{\Delta_{i}}(X_k) \ell_-'(d(\alpha,X_k)) 
        - \indic{T\setminus \Delta_{i}}(X_k) \ell_+'(d(\alpha,X_k)).$$
    From the estimate $|Y_k|\le \ell_+'(d(\alpha,X_k))$, it follows that $Y_k$ is sub-Gaussian.
    
    Since $\alpha$ is sticky, $\phi_{\Delta_i}'(\alpha)>0$ and by Hoeffding's inequality, %
    \begin{align*}
        \P{\hphi_{\Delta_i}'(\alpha) \le 0} 
        =\P*{\sum_{k=1}^{n} (Y_k-\E{Y_k}) \ge n \phi_{\Delta_i}'(\alpha)}
        \le \exp\left(-\frac{n\phi_{\Delta_i}'(\alpha)^2}{2\sigma_i^2} \right).
    \end{align*}
    \Cref{prop:sign-derivative} yields 
    the equality $M(\mu)=\set{\alpha}$ as well as the inclusion $\bigcap_{i=1}^m \set{\hphi_{\Delta_i}'(\alpha) > 0} \subset \set{M(\hmu_n)=\set{\alpha}}$.
    By the union bound,
    $$\P{M(\hmu_n)=\set{\alpha}} \ge 1 - \sum_{i=1}^m \exp\left(-\frac{n\phi_{\Delta_i}'(\alpha)^2}{2\sigma_i^2}\right).$$
\end{proof}

\begin{proof}[Proof of \Cref{prop:Gamma-S}]
    1. First, we consider the one-sided case and we fix $z\in [0,\rho_0]$.
    Since $\rho_0 < \rhos$, it follows that $\Gamma_z$ lies on the segment that connects $\alphas$ to its neighboring vertex in $\Delta_+$.
    If $x\in (\alphas, \Gamma_z]$, then $d(\Gamma_z,x) = d(\alphas, \Gamma_z) - d(\alphas, x) = z-S(x)$.
    If $x\in \Delta_+\setminus (\alphas, \Gamma_z]$, then $d(\Gamma_z,x) = d(\alphas, x) - d(\alphas, \Gamma_z) = S(x)-z$.
    If $x\in T\setminus \Delta_+$, $d(\Gamma_z,x) = d(\alphas, \Gamma_z) + d(\alphas, x) = z-S(x)$,
    hence in every subcase $d(\Gamma_z,x) = |z - S(x)|$ for every $x\in T$.
    A similar reasoning proves the equality in the two-sided case,
    with the added constraint that $x\in \Delta_- \cup \set{\alphas} \cup \Delta_+$.

    2. For every $z>0$, we have $\Gamma_z\in \Delta_+$ hence $S(\Gamma_z) = d(\alphas, \Gamma_z) = d(\alphas, \gamma_+(z))=z$.
    In the two-sided case, if $z<0$ then $\Gamma_z\in \Delta_-$ and $S(\Gamma_z) = -d(\alphas, \Gamma_z) = -d(\alphas, \gamma_-(|z|))= -|z|=z$.
    The equality also clearly holds at $z=0$, hence $S \circ \Gamma = \Id$.

    Let us show that $\Gamma$ is an isometry.
    In the one-sided case, we may apply the first item and obtain that $d(\Gamma_z, \Gamma_{z'}) = |z - S(\Gamma_{z'})| = |z - d(\alphas, \gamma_+(z'))| = |z-z'|$
    for every $z,z'\in [0,\rho_0]$.
    The proof in the two-sided case is similar.
\end{proof}

\begin{proof}[Proof of \Cref{prop:Psi}]
    1. We fix $z\in C\cap [-\rho_0,\rho_0]$ and we show that $\E[\big]{|\psi(z-S(X)) - \psi(S(X))|} < +\infty$.
    When $z>0$, the convexity of $\psi$ yields the estimate 
    $|\psi(z-S(X)) - \psi(0-S(X))|\le z\max\big(|\psi_-'(z-S(X))|, |\psi_+'(-S(X))|\big)$.
    When $z<0$ (which only happens in the two-sided case), we have instead
    $|\psi(z-S(X)) - \psi(0-S(X))|\le |z|\max\big(|\psi_+'(z-S(X))|, |\psi_-'(-S(X))|\big)$.
    Since 
    \begin{equation}
        \label{eq:psi-derivative}
        \forall u\in \R, \, \psi_+'(u) = \ell_+'(u)\indic{u\ge 0} - \ell_-'(|u|)\indic{u<0}
        \, \text{ and } \,  \psi_-'(u) = \ell_-'(u)\indic{u> 0} - \ell_+'(|u|)\indic{u\le 0},
    \end{equation}
    we obtain that $$|\psi(z-S(X)) - \psi(S(X))| \le \rho_0 \ell_+'\big(\max(|z-S(X)|, |S(X)|)\big) \le \rho_0 \ell_+'(d(\alphas, X) + \rho_0),$$
    and the integrability follows from \Cref{assum:clt-1}, hence $\Psi$ is well-defined.
    The convexity of $\Psi$ follows from that of $\psi$.
    In the one-sided case, by \Cref{prop:Gamma-S}, $\psi(z-S(X)) - \psi(S(X))$ coincides with $\ell(d(\Gamma_z, X)) - \ell(d(\alphas, X))$.
    In the two-sided case, we recall from \Cref{thm:localization-sticky} that $\P{X\in \Delta_-\cup\set{\alphas}\cup\Delta_+}=1$,
    hence $\psi(z-S(X)) - \psi(S(X)) = \ell(d(\Gamma_z, X)) - \ell(d(\alphas, X))$ almost surely and 
    we obtain in both cases that $\Psi(z) = \phi(\Gamma_z) - \phi(\alphas)$.
    Since $\alphas = \Gamma_0$ is the unique global minimizer of $\phi$, this last equality implies that $0$ is the unique minimizer of $\Psi$.

    2. We fix $z\in \operatorname{int}(C\cap [-\rho_0,\rho_0])$,
    we consider a sequence $(z_n)_{n\ge 1}$ in $C\cap [-\rho_0,\rho_0]$ which decreases strictly to $z$
    and we assess the ratio 
    $$\frac{\Psi(z_n) - \Psi(z)}{z_n - z} = \E[\Big]{\frac{\psi(z_n - S(X)) - \psi(z - S(X))}{z_n - z}}.$$
    Since $\psi$ is convex, 
    $$\left|\frac{\psi(z_n - S(X)) - \psi(z - S(X))}{z_n - z}\right| \le \max\big(|\psi_-'(z_n-S(X))|, |\psi_+'(z-S(X))|\big)
    \le \ell_+'(d(\alphas, X) + \rho_0)
    $$
    and by dominated convergence we obtain that $\psi_+'(z-S(X))$ is integrable and $\Psi_+'(z) = \E{\psi_+'(z-S(X))}$.
    A similar argument shows that $\Psi_-'(z) = \E{\psi_-'(z-S(X))}$ and $\Psi_+'(0)=\E{\psi_+'(-S(X))}$.
    Besides, it follows from the first item that $\Psi_+'(0) = \phi_{\Delta_+}'(\alphas) = 0$.
    In the two-sided case, $\Psi_-'(0) = -\phi_{\Delta_-}'(\alphas) = 0$.
    As a consequence, 
    $0 = \Psi_+'(0) - \Psi_-'(0) = \E{\psi_+'(-S(X)) - \psi_-'(-S(X))}$.
    Since $\psi$ is convex, the random variable $\psi_+'(-S(X)) - \psi_-'(-S(X))$ is nonnegative,
    and we obtain that in the two-sided case,
    \begin{equation}
        \label{eq:psi-derivatives-equal}
        \psi_+'(-S(X)) = \psi_-'(-S(X)) \text{ a.s.}
    \end{equation}
    
    3. From \eqref{eq:psi-derivative} it follows that
    $|\psi_+'(-S(X))| \le \ell_+'(|S(X)|) = \ell_+'(d(\alphas, X))$
    and combined with the integrability assumption we obtain
    that $\sigma^2 < +\infty$.
    Let us assume for contradiction that $\sigma^2=0$,
    which forces $\psi_+'(-S(X)) = 0$ \as 
    Since $\psi_+'$ has the closed form \eqref{eq:psi-derivative},
    since $\ell_-'$ takes positive values and $\ell_+'(z)=0$ implies $z=0$,
    it follows that $-S(X) = 0$ \as, hence $X=\alphas$ \as, which contradicts our initial assumption that $\mu\neq \delta_\alphas$.

    4. In the one-sided case, for each $n\ge 1$ and $z\in [0,\rho_0]$, $\hphi_n(\Gamma_z) - \hphi_n(\alphas)$ coincides everywhere with $\hPsi_n(z)$ by \Cref{prop:Gamma-S}.
    In the two-sided case, \Cref{prop:Gamma-S} guarantees the inclusion of events 
    $$\bigcap_{k\ge 1} \set[\big]{X_k \in \Delta_-\cup\set{\alphas}\cup\Delta_+} \subset \set[\big]{\forall n\ge 1, \forall z\in [-\rho_0,\rho_0],\,\hphi_n(\Gamma_z) - \hphi_n(\alphas) = \hPsi_n(z)}.$$
    By \Cref{thm:localization-sticky}, the left-hand side has probability $1$, hence the claim.

    5. Since $\alphas \in \operatorname{int}(F)$, 
    \Cref{thm:sticky-lln} applies: with probability $1$ there exists $N\ge 1$ such that for every $n\ge N$,
    $M(\hmu_n) \subset \set{\alphas} \cup \Delta_+$ in the one-sided case 
    (resp., $M(\hmu_n) \subset \Delta_- \cup \set{\alphas} \cup \Delta_+$ in the two-sided case).
    Moreover, by \Cref{thm:lln-sets}, $(M(\hmu_n))_{n\ge 1}$ is strongly consistent in one-sided Hausdorff distance.
    By intersecting these two almost-sure events, it follows that the event 
    \begin{equation}
        \label{eq:event}
        \set[\big]{\exists N\ge 1, \forall n\ge N, \, M(\hmu_n) \subset \Gamma(C\cap [-\rho_0,\rho_0])}
    \end{equation}
    has probability $1$.
    On the event \eqref{eq:event}, for large enough $n$ we have $\halpha_n = \Gamma_{z_0}$ where ${z_0}\in C\cap [-\rho_0, \rho_0]$, hence 
    $\hm_n = S(\halpha_n) = S(\Gamma_{z_0}) = z_0$, thus 
    $\Gamma_{\hm_n} = \Gamma_{z_0} = \halpha_n$.
    Finally, with probability 1, for large enough $n$,
    \begin{align*}
        \hPsi_n(\hm_n) &= \hphi_n(\Gamma_{\hm_n}) - \hphi_n(\alphas)
        = \hphi_n(\halpha_n) - \hphi_n(\alphas)
        \le \hphi_n(\Gamma_z) - \hphi_n(\alphas) = \hPsi_n(z),
    \end{align*}
    where $z\in C\cap [-\rho_0,\rho_0]$ is arbitrary, thus 
    $\hm_n\in \argmin_{z\in C\cap [-\rho_0,\rho_0]} \hPsi_n(z)$.
\end{proof}

\begin{lemma}
    \label{lemma:clt}
    Consider the setting of \Cref{sec:CLT}.
    \begin{enumerate}
        \item Assume that $\E{\ell_+'(d(\alphas, X))^2} < +\infty$. The probabilities $\P{\hphi_{\Delta_+}'(\alphas)>0}$ and $\P{\hphi_{\Delta_+}'(\alphas)\ge 0}$
        go to $1/2$ as $n\to \infty$.
        Additionally, in the two-sided case, $\P{\hphi_{\Delta_-}'(\alphas)>0}$ and $\P{\hphi_{\Delta_-}'(\alphas)\ge 0}$
        converge to $1/2$.
        \item Assume \ref{assum:clt-2}. Let $(z_n)_{n\ge 1}$ be a sequence of real numbers in $\operatorname{int}(C\cap [-\rho_0,\rho_0])$ such that 
        $\lim_n z_n = 0$. %
        Then \begin{equation}
            \label{eq:convergence-hPsi}
            \sqrt n \big((\hPsi_n)_+'(z_n) -  \Psi_+'(z_n)\big) \xrightarrow[]{(d)} N(0, \sigma^2).
        \end{equation}
    \end{enumerate}
\end{lemma}

\begin{proof}[Proof of \Cref{lemma:clt}]
    1. For each $k\in [n]$, let $$Y_k = \indic{T\setminus \Delta_+}(X_k) \ell_+'(d(\alphas,X_k))
    - \indic{\Delta_+}(X_k) \ell_-'(d(\alphas,X_k)),$$
    so that $\hphi_{\Delta_+}'(\alphas) = \frac 1n \sum_{k=1}^{n} Y_k$. 
    Our assumption that $\E{\ell_+'(d(\alphas, X))^2} < +\infty$ implies that 
    $Y_k$ is square-integrable, as well as  
    $\alphas\in F$.
    We know that $\phi_{\Delta_+}'(\alphas) = 0 < +\infty$, hence
    by the third item of \Cref{prop:derivative},
    the $Y_k$ have mean zero.  
    Moreover, we observe that $Y_k = \psi_+'(-S(X_k))$.
    The classical central limit theorem then yields
    convergence in distribution of $(\sqrt n \hphi_{\Delta_+}'(\alphas))_{n\ge 1}$ 
    to $N(0, \sigma^2)$.
    It follows that $\P{\hphi_{\Delta_+}'(\alphas)>0} \to_n 1/2$ and $\P{\hphi_{\Delta_+}'(\alphas)\ge 0} \to_n 1/2$.
    In the two-sided case, the convergence of $\P{\hphi_{\Delta_-}'(\alphas)>0}$ and $\P{\hphi_{\Delta_-}'(\alphas)\ge 0}$ 
    is obtained in a similar fashion.

    2. For each $k,n \ge 1$, we let $Y_{n,k} = \psi_+'(z_n - S(X_k))$ and we note by \Cref{prop:Psi} that $\E{Y_{n,k}} = \Psi_+'(z_n)$.
    
    Let us show that $Y_{n,1} \xrightarrow[]{L^2} \psi_+'(-S(X_1))$.
    Since $\psi$ is convex, the second item of \Cref{lemma:integral-rep} carries over to $\psi$.
    Consequently, if $z_n \downarrow 0$, then $Y_{n,1} \xrightarrow[]{a.s.} \psi_+'(-S(X_1))$,
    while if $z_n \uparrow 0$ (only in the two-sided case), then $Y_{n,1} \xrightarrow[]{a.s.} \psi_-'(-S(X_1))$.
    Combining these facts with \eqref{eq:psi-derivatives-equal}, we obtain in the general case that $Y_{n,1} \xrightarrow[]{a.s.} \psi_+'(-S(X_1))$.
    From there, we finish the proof via dominated convergence.
    Regarding domination, 
    let us note the following estimate:  
        \begin{equation}
            \label{eq:domination}
            \forall z\in C\cap [-\rho_0, \rho_0], \, |\psi_+'(z-S(X))| \le \ell_+'(d(\alphas, X) + \rho_0).
        \end{equation}
    Indeed, it follows from \eqref{eq:psi-derivative} that
    $|\psi_+'(z-S(X))| \le \ell_+'(|z-S(X)|) \le \ell_+'(d(\alphas, X) + \rho_0)$.
    We apply \eqref{eq:domination} and \Cref{assum:clt-2} twice to derive that 
    $Y_{n,1}$ and $\psi_+'(-S(X_1))$ are in $L^2$, and
    $|Y_{n,1} - \psi_+'(-S(X_1))|^2 \le 4\ell_+'(d(\alphas, X_1) + \rho_0)^2$.
    The left-hand side converges almost surely to $0$, while the right-hand side is integrable, 
    hence it follows from the dominated convergence theorem that $Y_{n,1} \xrightarrow[]{L^2} \psi_+'(-S(X_1))$.

    We observe that $\sqrt n \big((\hPsi_n)_+'(z_n) -  \Psi_+'(z_n)\big) = \frac1{\sqrt n} \sum_{k=1}^n (Y_{n,k} - \Psi_+'(z_n))$
    and we define the random variable 
    $$U_n\coloneqq\frac1{\sqrt n} \sum_{k=1}^n (Y_{n,k} - \Psi_+'(z_n)) - \frac1{\sqrt n} \sum_{k=1}^n \psi_+'(-S(X_k)),$$
    which is in $L^2$, with mean $0$ and variance 
    $$\V{U_n} = \V{Y_{n,1} - \psi_+'(-S(X_1))} \le \norm[\big]{Y_{n,1} - \psi_+'(-S(X_1))}_{L^2}^2,$$
    which goes to $0$.
    Consequently, $U_n = o_\PP(1)$ and 
    $$\sqrt n \big((\hPsi_n)_+'(z_n) -  \Psi_+'(z_n)\big) = \frac1{\sqrt n} \sum_{k=1}^n \psi_+'(-S(X_k)) + o_\PP(1).$$
    \Cref{eq:convergence-hPsi} then follows from the classical central limit theorem and Slutsky's theorem.
\end{proof}

\begin{proof}[Proof of \Cref{prop:proba-localization-two-sided}]
    Note the inclusions of events
    $$\set{\hphi_{\Delta_+}'(\alphas) <0} 
    \subset \set{M(\hmu_n) \subset \Delta_+}
    \subset \set{\halpha_n \in \Delta_+}
    \subset \set{\hphi_{\Delta_+}'(\alphas) \le 0},
    $$
    and apply the first item of \Cref{lemma:clt} to obtain that 
    $\lim_n\P{M(\hmu_n) \subset \Delta_+} = \lim_n\P{\halpha_n \in \Delta_+} = \frac 12$.
    The direction $\Delta_-$ is treated similarly.
\end{proof}

\begin{proof}[Proof of \Cref{thm:two-sided-clt}]
    1. We fix $t>0$ and we let $z_n = t n^{-1/(2a_+)}$, as well as $\Omega_n = \set{\hm_n\in \argmin_{z\in [-\rho_0,\rho_0]} \hPsi_n(z)}$.
    Since $\hPsi_n$ is convex on $\R$, the following inclusion of events holds:
    $$\Big(\set[\Big]{(\hPsi_n)_+'(z_n) >0} \cap \Omega_n\Big) 
    \subset 
    \Big(\set[\Big]{ \hm_n \le z_n}\cap \Omega_n\Big) 
    \subset 
    \Big(\set[\Big]{(\hPsi_n)_+'(z_n)\ge 0}\cap \Omega_n\Big).$$
    By the fifth item of \Cref{prop:Psi}, $\P{\Omega_n}\to 1$.
    Therefore, it suffices to establish the convergence of $\P[\big]{(\hPsi_n)_+'(z_n)>0}$ and $\P[\big]{(\hPsi_n)_+'(z_n)\ge 0}$ as $n\to \infty$.
    By discarding finitely many terms if necessary, we assume \wlo that $z_n\in (0, \rho_0)$ for each $n\ge 1$.
    Condition (b) of \Cref{thm:two-sided-clt} guarantees that $\sqrt n \Psi_+'(z_n) \to K_+ t^{a_+}$.
    Then, it follows from \eqref{eq:convergence-hPsi} in \Cref{lemma:clt} that $\sqrt n (\hPsi_n)_+'(z_n) \xrightarrow[]{(d)} N(K_+ t^{a_+}, \sigma^2)$.
    Consequently, $\lim_n\P[\big]{(\hPsi_n)_+'(z_n)>0} = \lim_n\P[\big]{(\hPsi_n)_+'(z_n)\ge 0} = \P{Z\le \frac{K_+ }{\sigma}t^{a_+}}$.

    2. The proof is similar to that of the first item.

    3. Let $Z'=\frac{\sgn(Z)+1}{2} \Big(\frac{\sigma}{K_+}|Z|\Big)^{1/a}  
        + \frac{\sgn(Z)-1}{2} \Big(\frac{\sigma}{K_-}|Z|\Big)^{1/a}$.
    For each $t\neq 0$, the CDF of $Z'$ satisfies
    $$F_{Z'}(t) =  \P[\Big]{Z\le \frac{K_+}{\sigma}t^{a}} \indic{t> 0}
    + \P[\Big]{Z\le -\frac{K_-}{\sigma}|t|^{a}} \indic{t< 0}.$$
    For convenience, let us denote the CDF of $n^{1/(2a)}\hm_n$ by $F_n$, so that
    $F_n(t) \to_n F_{Z'}(t)$ for each $t\neq 0$.
    Since $\R\setminus\{0\}$ is dense in $\R$,
    \cite[Lemma~3(i) p.~275]{chow2003probability} yields
    the first distributional convergence.
    The second convergence follows from 
    the equality $|\hm_n| = d(\halpha_n, \alphas)$ and the continuous mapping theorem.
\end{proof}

\begin{lemma}
    \label{lemma:assum-growth-second-moment}
    Assume \ref{assum:growth} holds and $\E{\ell_+'(d(o,X))^2} < +\infty$.
    Then for every $\alpha\in T$ and every $\rho > 0$, it holds that 
    $\E{\ell_+'(d(\alpha, X) + \rho)^2} < +\infty$.
\end{lemma}

\begin{proof}[Proof of \Cref{lemma:assum-growth-second-moment}]
    Let us denote by $q\ge 1$ an integer such that $q\ge d(o,\alpha) + \rho$.
    We note that $\indic{d(o,X)<R} \ell_+'(d(\alpha, X) + \rho)^2 \le \ell_+'(R + d(o,\alpha) + \rho)^2$
    and
    \begin{align*}
        \indic{d(o,X)\ge R} \ell_+'(d(\alpha, X) + \rho)^2
        &\le \indic{d(o,X)\ge R} \ell_+'(d(o, X) + q)^2
        \\&\le C^2 \indic{d(o,X)\ge R} \ell_+'(d(o, X) + q-1)^2
        \\&\le C^{2q} \ell_+'(d(o, X))^2,
    \end{align*}
    hence the claim.
\end{proof}

\begin{lemma}
    \label{lemma:equiv-assum}
    Consider the two-sided partly sticky setting of \Cref{sec:two-sided-clt}.
    Assume \ref{assum:clt-1}: $\E{\ell_+'(d(\alphas, X) + \rho_0)} < +\infty$ for some $\rho_0\in (0, \rhos)$.
    For $\varepsilon \in \set{+,-}$, define $\Phi^\varepsilon: [0,\rho_0] \to \R$, $t\mapsto \phi(\gamma_\varepsilon(t))$.
    The following assertions 1. and 2. are equivalent:
    \begin{enumerate}
        \item $\frac{\Psi_+'(z)}{z^{a_+}} \to K_+$ as $z\downarrow 0$ and $\frac{\Psi_+'(z)}{|z|^{a_-}} \to -K_-$ as $z\uparrow 0$,
        \item $\frac{(\Phi^+)_+'(t)}{t^{a_+}} \to K_+$ as $t\downarrow 0$ and  $\frac{(\Phi^-)_+'(t)}{t^{a_-}} \to K_-$ as $t\downarrow 0$,
    \end{enumerate}
    where $a_-, a_+,K_-,K_+$ are positive constants.
\end{lemma}

\begin{proof}[Proof of \Cref{lemma:equiv-assum}]
    By the first item of \Cref{prop:Psi}, for every $z\in (0,\rho_0)$, 
    $\Psi(z) 
    = \phi(\gamma_+(z)) - \phi(\gamma_+(0)) 
    = \Phi^{+}(z) - \Phi^{+}(0) $
    and consequently $\Psi_+'(z) = (\Phi^{+})_+'(z)$.
    It follows immediately that 
    $$\frac{\Psi_+'(z)}{z^{a_+}} \xrightarrow[z\downarrow 0]{} K_+ \iff \frac{(\Phi^+)_+'(t)}{t^{a_+}} \xrightarrow[t\downarrow 0]{} K_+.$$

    Similarly, for every $z\in (-\rho_0,0)$, 
    $\Psi(z) 
    = \phi(\gamma_-(-z)) - \phi(\gamma_-(0)) 
    = \Phi^{-}(-z) - \Phi^{-}(0) $
    and therefore $\Psi_+'(z) = -(\Phi^-)_-'(|z|)$.
    Let us suppose that $\frac{(\Phi^-)_+'(t)}{t^{a_-}} \to K_-$ as $t\downarrow 0$.
    We fix some arbitrary $\lambda \in (0,1)$. Since $\Phi^{-}$ is convex, 
    $$\frac{(\Phi^{-})_+'((1-\lambda)|z|)}{|z|^{a_-}} 
    \le \frac{(\Phi^{-})_-'(|z|)}{|z|^{a_-}}
    \le \frac{(\Phi^{-})_+'(|z|)}{|z|^{a_-}}.
    $$
    Letting $z\uparrow 0$, it follows that 
    $$(1-\lambda)^{a_-} K_- \le \liminf_{z\uparrow 0}\frac{(\Phi^{-})_-'(|z|)}{|z|^{a_-}}
    \le \limsup_{z\uparrow 0}\frac{(\Phi^{-})_-'(|z|)}{|z|^{a_-}}
    \le K_-,
    $$
    and letting $\lambda \to 0$ we obtain that $\lim_{z\uparrow 0}(\Phi^{-})_-'(|z|) / |z|^{a_-} = K_-$,
    hence $\Psi_+'(z) / |z|^{a_-} \to -K_-$ as $z\uparrow 0$. The proof of the converse is similar.
\end{proof}

\begin{proof}[Proof of \Cref{corol:clt-regular}]
    For any $\rho_0 > 0$, the mean value theorem yields the estimate
    $\ell'(d(\alphas, X) + \rho_0) \le \ell'(d(\alphas, X)) + \rho_0 L$,
    hence $$\ell'(d(\alphas, X) + \rho_0)^2 \le 2\ell'(d(\alphas, X))^2 + 2\rho_0^2 L^2$$
    and \Cref{assum:clt-2} is satisfied.

    Let us define the functions $s:\R\to \R$, $z\mapsto \indic{z\ge 0} - \indic{z < 0}$ and $g: \R\to \R$, $z\mapsto (\ell'(|z|) - \ell'(0))(\indic{z\ge 0} - \indic{z < 0})$,
    so that $\psi_+'(z) 
    = \ell_+'(z)\indic{z \ge 0} - \ell_-'(|z|)\indic{z<0} 
    = \ell'(|z|)s(z)
    = g(z) + \ell'(0)s(z)
    $.
    Recall from the second item of \Cref{prop:Psi} that $\Psi_+'(z) = \E{\psi_+'(z-S(X))}$, hence
    $\Psi_+'(z) = \E{g(z-S(X))} + \ell'(0)\E{s(z-S(X))}$ for each $z\in (-\rho_0,\rho_0)$.

    The function $g$ is differentiable with $g'(z) = \ell''(|z|)$, hence $0\le g'\le L$ and $g$ is $L$-Lipschitz.
    It follows by the dominated convergence theorem that $z\mapsto \E{g(z-S(X))}$ is differentiable at $0$ with derivative $\E{\ell''(d(\alphas, X))}$.

    Next, we consider the case where $\ell'(0)>0$ and we turn to the differentiability of $z\mapsto \E{s(z-S(X))}$.
    We denote by $\F$ the CDF of the random variable $S(X)$,
    so that $\E{s(z-S(X))} = 2\F(z)-1$.
    We assume \wlo that $\varepsilon\le \rho_0$.
    For any Borel subset $B$ of $(-\varepsilon,\varepsilon)$, we observe that 
    $\P{S(X)\in B} = \P{X\in \Gamma(B)} = \int_B f(\Gamma(z)) \diff z$
    and we define $$h:z\mapsto \begin{cases}
        f(\gamma_+(z)) &\text{ if } z>0, \\
        f(\gamma_-(-z)) &\text{ if } z<0, \\
        \eta &\text{ if } z=0.
    \end{cases}$$
    We denote by $S\push \mu$ the pushforward measure and by $(S\push \mu)_{\mid (-\varepsilon, \varepsilon)}$ its restriction to $(-\varepsilon, \varepsilon)$.
    Then $(S\push \mu)_{\mid (-\varepsilon, \varepsilon)}$ has density $h$ \wrt the Lebesgue measure on $(-\varepsilon, \varepsilon)$ and $h$ is continuous at $0$ by condition (c).
    It follows that $\F$ is differentiable at $0$ with $\F'(0) = \eta$,
    hence $\Psi_+'$ is differentiable at $0$ with derivative $\E{\ell''(d(\alphas, X))} + 2\ell'(0) \eta$.

    Since $K>0$ by assumption, \Cref{thm:two-sided-clt} applies with 
    $a_+=a_-=1$, $K_+=K_-=K$ and $\sigma^2 = \E{\ell'(d(\alphas, X))^2}$. 
    Note that $\sigma^2=\E{\psi_+'(-S(X))^2}=\E{\ell'(|-S(X)|)^2} = \E{\ell'(d(\alphas, X))^2}$.
\end{proof}

\begin{proof}[Proof of \Cref{prop:one-sided-collapse}]
    Consider $\Delta\in \dir_{\alphas}\setd{\Delta_+}$. Then $\phi_\Delta'(\alphas)>0$
    and
    $$\hphi_\Delta'(\alphas) = 
    \frac1n \sum_{k=1}^{n}[\indic{T\setminus \Delta}(X_k) \ell_+'(d(\alphas,X_k))- \indic{\Delta}(X_k) \ell_-'(d(\alphas,X_k))].$$
    Since $\alphas \in \operatorname{int}(F)$, we may apply the second item of \Cref{prop:derivative} and the strong law of large numbers to derive 
    $\hphi_\Delta'(\alphas) \xrightarrow[]{a.s.} \phi_\Delta'(\alphas)$.
    It follows that $\indic{\hphi_\Delta'(\alphas) > 0} \xrightarrow{a.s.} 1$,
    hence by dominated convergence
    $\P{\hphi_\Delta'(\alphas)>0} \to 1$ as $n\to \infty$.

    Combining these limits with the first item of \Cref{lemma:clt} and with the inclusions of events
    $$\bigcap_{\Delta \in \dir_{\alphas}} \set{\hphi_\Delta'(\alphas)>0} 
    \subset \set{M(\hmu_n) = \set{\alphas}}
    \subset \set{\halpha_n = \alphas}
    \subset \set{\hphi_{\Delta_+}'(\alphas) \ge 0},
    $$
    we obtain that $\lim_n \P{M(\hmu_n) = \set{\alphas}} 
    =\lim_n \P{\halpha_n = \alphas}  = \frac 12$.
    Besides,
    $$\set{\hphi_{\Delta_+}'(\alphas) <0} 
    \subset \set{M(\hmu_n) \subset \Delta_+}
    \subset \set{\halpha_n \in \Delta_+}
    \subset \set{\hphi_{\Delta_+}'(\alphas) \le 0},
    $$
    hence $\lim_n\P{M(\hmu_n) \subset \Delta_+} = \lim_n\P{\halpha_n \in \Delta_+} = \frac 12$.
\end{proof}

\begin{proof}[Proof of \Cref{thm:one-sided-clt}]
    1. We repeat the same proof as for the first item in \Cref{thm:two-sided-clt},
    with the only difference that 
    $\Omega_n = \set{\hm_n\in \argmin_{z\in [0,\rho_0]} \hPsi_n(z)}$.

    2. By the fifth item of \Cref{prop:Psi},
    $\P{\hm_n \in [0,\rho_0]}\to 1$, hence 
    $\P*{n^{1/(2a_+)} \hm_n \le t} \to 0$ for each $t<0$.

    3. First, we show that $\P{n^{1/(2a_+)} \hm_n \le 0} \to 1/2$.
    Note that $\P{\hm_n < 0} \le \P{\hm_n \notin [0,\rho_0]}\to 0$.
    Moreover, by \Cref{prop:one-sided-collapse}, $\P{\hm_n = 0} \to 1/2$,
    hence the claim.

    Next, we let $Z'=(\frac{\sigma}{K_+}\max(Z,0))^{1/a_+}$
    and we have shown that its CDF satisfies
    $\P{n^{1/(2a_+)}\hm_n \le t} \to_n F_{Z'}(t)$ for each $t\in \R$.
    The convergence in distribution of $(n^{1/(2a_+)}\hm_n)_{n\ge 1}$ follows,
    and we obtain that of $(n^{1/(2a_+)}d(\halpha_n, \alphas))_{n\ge 1}$ by the continuous mapping theorem.
\end{proof}

\subsection{Proofs for Section~\ref{sec:collisions}}

Before turning to the proof of \Cref{thm:estimator-stickiness}, we need two technical lemmas.

\begin{lemma}
    \label{lemma:bound-proba-equality}
    Let $Y,Y':\Omega \to \R$ denote two \iid random variables.
    Then $\P{Y=Y'}\le \sup_{y\in \R} \P{Y=y}.$
\end{lemma}

\begin{proof}[Proof of \Cref{lemma:bound-proba-equality}]
    We denote by $\nu$ the probability measure on $\R$ induced by $Y$.
    By the law of the unconscious statistician and Tonelli's theorem,
    $$\P{Y=Y'} = \int_{\R^2} \indic{y=y'} \diff (\nu \otimes \nu)(y,y')
    = \int_\R \nu(\set{y}) \diff \nu(y) \le \sup_{y\in \R} \nu(\set{y}) = \sup_{y\in \R} \P{Y=y}.
    $$
\end{proof}

\begin{lemma}
    \label{lemma:bound-supremum}
    Let $Y, Y_1, Y_2,\ldots$ be random variables such that $Y_n \xrightarrow[]{(d)} Y$.
    If $Y$ is atomless (\ie $\P{Y=y} = 0$ for every $y\in \R$), then $\sup_{y\in \R} \P{Y_n=y} \to 0$ as $n\to \infty$.
\end{lemma}

\begin{proof}[Proof of \Cref{lemma:bound-supremum}]
    We denote by $F_n$ (resp., $F$) the CDF of $Y_n$ (resp., $Y$).
    Since $Y$ is atomless, $F$ is continuous.
    By P\'olya's theorem \cite[Satz~I]{polya1920uber},
    $(F_n)_{n\ge 1}$ converges uniformly on $\R$ to $F$.
    Moreover, for each $y\in \R$,
    $$\P{Y_n=y} = F_n(y) - F_n(y-) = F_n(y) - F(y) - [F_n(y-) - F(y-)] \le 2 \norm{F_n - F}_\infty,$$
    hence the claim.
\end{proof}

\begin{proof}[Proof of \Cref{thm:estimator-stickiness}]
    We observe that $\halpha_{n,1},\ldots, \halpha_{n,j_n}$ are \iid
    and we let $p_n = \P{\halpha_{n,1} = \halpha_{n,2}}$,
    so that $\E{\hp_n} = p_n$.
    We view 
    $$\hp_n = \frac{2}{j_n(j_n-1)}\sum_{1\le j < j' \le j_n} \indic{\halpha_{n,j} = \halpha_{n,j'}}$$ 
    as a function of the random elements $\halpha_{n,1},\ldots, \halpha_{n,j_n}$.
    Changing the value of $\halpha_{n,1}$ changes the value of at most $j_n - 1$ summands, 
    and each summand changes by at most $1$ in absolute value.
    Consequently, the value of $\hp_n$ changes by at most $2/j_n$ in absolute value.
    It follows by the bounded differences inequality
    \cite[Theorem 6.2]{boucheron2013concentration}
    that 
    \begin{equation}
        \label{eq:deviation-pn}
        \forall t\ge 0, \quad \P{|\hp_n - p_n| \ge t} \le 2e^{-j_nt^2/2}
    \end{equation}
    Throughout the proof we let $q_n = \P{\halpha_{n,1} = \alphas}$.

    1. a) We define the random variable $J_n = \sum_{j=1}^{j_n} \indic{\halpha_{n,j}\neq \alphas}$, so that 
    $J_n \sim \Bin(j_n, 1-q_n)$.
    Exactly $j_n - J_n$ blocks have Fr\'echet $\ell$-mean equal to $\alphas$,
    hence $\sum_{1\le j < j' \le j_n} \indic{\halpha_{n,j} = \halpha_{n,j'}} \ge \binom{j_n - J_n}{2}$
    and $\hp_n \ge \frac{(j_n - J_n)(j_n - J_n - 1)}{j_n(j_n - 1)}$.
    This last ratio is greater than $5/8$ as soon as the number of bad blocks $J_n$ is less than $j_n/5$.
    Indeed, if $J_n < j_n/5$, then $5(j_n - J_n) > 4j_n$, \ie $5(j_n - J_n) \ge 4j_n + 1$,
    hence both $j_n - J_n \ge (4j_n + 1)/5$ and $j_n - J_n - 1 \ge (4j_n -4)/5$,
    therefore $\frac{(j_n - J_n)(j_n - J_n - 1)}{j_n(j_n - 1)} \ge \frac{(4j_n + 1)(4j_n -4)}{25 j_n (j_n-1)} > 5/8$.
    The inclusion of events $\set{\hR_n \neq 0} \subset \set{ J_n \ge \ceil{j_n / 5}}$ follows.
    By Markov's inequality,
    \begin{equation}
        \label{eq:proba-not-zero-markov}
        \P{\hR_n \neq 0} \le \P{ J_n \ge \ceil{j_n / 5}} = \P{ J_n \ge j_n / 5} \le 5(1-q_n).
    \end{equation}
    Since $\alphas\in \operatorname{int}(F)$ and $b_n \to +\infty$, \Cref{thm:sticky-lln} applies and yields $q_n \to 1$,
    hence the claim.

    1. b) By the Borel--Cantelli lemma, it suffices to show that $\sum_{n=2}^{\infty} \P{\hR_n \neq 0} < +\infty$.
    We partition the integers into the subsets $\mN_1 = \set{n\geq 2: j_n \le 5}$ and $\mN_2 = \set{n\geq 2: j_n \ge 6}$
    and we show separately that $\sum_{n\in \mN_1} \P{\hR_n \neq 0}$ and $\sum_{n\in \mN_2} \P{\hR_n \neq 0}$ are finite.

    Let us leverage \eqref{eq:proba-not-zero-markov} and prove that $\sum_{n\in \mN_1} (1-q_n) < +\infty$.
    For $i\in [m]$ and $k\geq 1$, we define $Y_{i,k} = \indic{T\setminus \Delta_{i}}(X_k) \ell_+'(d(\alphas,X_k))
    - \indic{\Delta_{i}}(X_k) \ell_-'(d(\alphas,X_k))$.
    Since $\E{\ell_+'(d(\alphas, X))^2} < +\infty$, $Y_{i,k}$ is square-integrable and $\E{Y_{i,k}} = \phi_{\Delta_i}'(\alphas)>0$.
    For each $n\ge 2$,
    the second item of \Cref{prop:optimality} applied to the empirical measure $\frac{1}{b_n} \sum_{k\in B_{n,1}} \delta_{X_k}$ yields the inclusion of events
    $\bigcap_{i=1}^m \set{\frac{1}{b_n} \sum_{k\in B_{n,1}} Y_{i,k} > 0} \subset \set{\halpha_{n,1} = \alphas}$
    and consequently
    \begin{equation}
        \label{eq:bound-q_n-1}
        1-q_n \le \sum_{i=1}^{m} \P[\Big]{\frac{1}{b_n} \sum_{k\in B_{n,1}} Y_{i,k}\le 0},
    \end{equation}
    hence
    \begin{equation}
        \label{eq:bound-q_n-2}
        \sum_{n\in \mN_1} (1-q_n) \le \sum_{i=1}^{m} \sum_{j=2}^5 \sum_{\substack{n \ge 2,\\ j_n=j}}
        \P[\Big]{\frac{1}{\floor{n/j}} \sum_{k\in B_{n,1}} Y_{i,k}\le 0}
        .
    \end{equation}
    Fix $j\in \set{2,\ldots,5}$ and note that for each integer $b\ge 1$, $|\set{n\geq 2: \floor{n/j} = b}| \le j$.
    Therefore,
    \begin{align}
        \label{eq:bound-sum-b}
        \sum_{\substack{n \ge 2,\\ j_n=j}}\P[\Big]{\frac{1}{\floor{\frac nj}} \sum_{k\in B_{n,1}} Y_{i,k}\le 0}
        &\le 
        \sum_{b=1}^\infty j \P[\Big]{\frac{1}{b} \sum_{k=1}^b Y_{i,k}\le 0} \nonumber
        \\&\le
        j\sum_{b=1}^\infty \P[\Big]{\Big|\frac{1}{b} \sum_{k=1}^b Y_{i,k} - \phi_{\Delta_i}'(\alphas)\Big| \ge \phi_{\Delta_i}'(\alphas)}.
    \end{align}
    By the Hsu--Robbins theorem \cite[Corollary~2 p.~393]{chow2003probability} applied to the sequence $(Y_{i,k})_{k\ge 1}$,
    the sum in the right-hand side of \eqref{eq:bound-sum-b} is finite, hence it follows from \eqref{eq:bound-q_n-2} that $\sum_{n\in \mN_1} (1-q_n) < +\infty$
    and $\sum_{n\in \mN_1} \P{\hR_n \neq 0} < +\infty$.
    
    We turn to the finiteness of $\sum_{n\in \mN_2} \P{\hR_n \neq 0}$.
    For convenience, we let $h_n = \ceil{j_n / 5}$
    and we recall the estimate $\P{\hR_n \neq 0} \le \P{J_n \ge h_n}$ from \eqref{eq:proba-not-zero-markov}.
    The inclusion of events 
    $$\set{J_n \ge h_n} \subset \bigcup_{\substack{\mJ \subset [j_n]\\ |\mJ| = h_n}} \bigcap_{j\in \mJ} \set{\halpha_{n,j}\neq \alphas}$$
    implies the bound
    $\P{J_n \ge h_n} \le \binom{j_n}{h_n} (1-q_n)^{h_n} \le 2^{j_n} (1-q_n)^{h_n}$,
    hence $\P{\hR_n \neq 0} \le 2^{j_n} (1-q_n)^{h_n} \le (32(1-q_n))^{h_n}$.
    From \eqref{eq:bound-q_n-1} and Chebyshev's inequality, we obtain that
    \begin{align}
        1-q_n 
        \le \sum_{i=1}^{m} \P[\Big]{\Big|\frac{1}{b_n} \sum_{k\in B_{n,1}} Y_{i,k} - \phi_{\Delta_i}'(\alphas) \Big|\ge \phi_{\Delta_i}'(\alphas)} 
        \le \frac{A}{b_n}, 
    \end{align}
    where $A = \sum_{i=1}^{m} \frac{\V{Y_{i,1}}}{\phi_{\Delta_i}'(\alphas)^2}$,
    thus $\P{\hR_n \neq 0} \le (32A/b_n)^{h_n}$.
    By the definition of $b_n$, it holds that $n/j_n < b_n +1$, hence $n < 2j_n b_n \le 10 h_n b_n$.
    Consequently,
    \begin{align*}
        \P{\hR_n \neq 0} &\le n^{-2} 10^2 h_n^2 b_n^2 \Big(\frac{32A}{b_n}\Big)^{h_n}
         = n^{-2} (320 A)^2 h_n^2 \Big(\frac{32A}{b_n}\Big)^{h_n-2}
    \end{align*}
    Since $b_n\to +\infty$, by discarding finitely many terms, we may assume \wlo that $32A/b_n \le 1/2$ for each $n\in \mN_2$.
    If $n\in \mN_2$, then $h_n\ge 2$ and we obtain the bound
    $\P{\hR_n \neq 0} \le n^{-2} (320 A)^2 h_n^2 /2^{h_n-2}$.
    The sequence $(h^2 /2^{h-2})_{h\geq 2}$ is bounded,
    which entails $\P{\hR_n \neq 0} \le A'/n^2$ for some $A'>0$ and every $n\in \mN_2$,
    hence finally $\sum_{n\in \mN_2} \P{\hR_n \neq 0} < +\infty$.

    2. a) We observe that
    \begin{align*}
        p_n &= \P{\halpha_{n,1} = \halpha_{n,2}, \halpha_{n,1} = \alphas} + \P{\halpha_{n,1} =\halpha_{n,2}, \halpha_{n,1} \neq \alphas}
        \\&= q_n^2 + \P{\halpha_{n,1} = \halpha_{n,2},  \halpha_{n,1} \neq \alphas}
        \\&\le q_n^2 + \P{\halpha_{n,1} \neq \alphas, \halpha_{n,2} \neq \alphas}
        = q_n^2 + (1-q_n)^2,
    \end{align*}
    hence $q_n^2 \le p_n \le q_n^2 + (1-q_n)^2$.
    Since we assumed 
    $\alphas \in \operatorname{int}(F)$, $\E{\ell_+'(d(\alphas, X))^2} < +\infty$ and $b_n \to +\infty$,
    \Cref{prop:one-sided-collapse} applies and yields the convergence $q_n \to 1/2$.
    By discarding finitely many terms, we may therefore assume \wlo that $3/16 \le p_n \le 9/16$.
    Let us note that if $\hp_n\notin (1/8, 5/8]$, then $|\hp_n - p_n| \ge 1/16$,
    hence $\set{\hR_n \neq 1} \subset \set{|\hp_n - p_n| \ge 1/16}$.
    Combined with \eqref{eq:deviation-pn}, we obtain the bound 
    $\P{\hR_n \neq 1} \le 2e^{-j_n / 512}$.
    Since $j_n\to +\infty$, $\P{\hR_n \neq 1} \to 0$. 

    2. b) Continuing from the last estimate,
    $\sum_{n=2}^\infty \P{\hR_n \neq 1} \le \sum_{n=2}^\infty 2e^{-j_n / 512} < +\infty$ and 
    the claim follows from the Borel--Cantelli lemma.

    3. By the definition of $\hR_n$ and Markov's inequality,
    $\P{\hR_n \neq 2} = \P{\hp_n > 1/8} \le 8 p_n$ and it remains to establish that $p_n \to 0$.
    We write $I_0=\set{i_-, i_+}$, and to improve legibility we let $\Delta_+, a_-, a_+, K_-, K_+$ stand for $\Delta_{i_+}, a_{i_-}, a_{i_+}, K_{i_-}, K_{i_+}$.
    We define 
    \begin{align*}
        Y_n &= b_n^{1/(2a_+)} S(\halpha_{n,1}) \indic{\halpha_{n,1} \in \Delta_+} 
                + b_n^{1/(2a_-)} S(\halpha_{n,1}) \indic{\halpha_{n,1} \notin \Delta_+}, \\
        Y_n' &= b_n^{1/(2a_+)} S(\halpha_{n,2}) \indic{\halpha_{n,2} \in \Delta_+} 
                + b_n^{1/(2a_-)} S(\halpha_{n,2}) \indic{\halpha_{n,2} \notin \Delta_+},
    \end{align*}
    and we note that by \Cref{lemma:bound-proba-equality}, 
    \begin{equation}
        \label{bound-p_n}
        p_n\le \P{Y_n = Y_n'} \le \sup_{y\in \R} \P{Y_n = y}.
    \end{equation}
    By \Cref{lemma:equiv-assum}, condition (b) of \Cref{thm:two-sided-clt} is fulfilled.
    We let $Y=\frac{\sgn(Z)+1}{2} (\frac{\sigma}{K_+}|Z|)^{1/a_+}  
        + \frac{\sgn(Z)-1}{2} (\frac{\sigma}{K_-}|Z|)^{1/a_-}$.
    Then for each $t\neq 0$, the CDF of $Y$ satisfies
    $$F_{Y}(t) =  \P[\Big]{Z\le \frac{K_+}{\sigma}t^{a_+}} \indic{t> 0}
    + \P[\Big]{Z\le -\frac{K_-}{\sigma}|t|^{a_-}} \indic{t< 0}.$$
    Moreover, $b_n\to +\infty$, hence we may apply \Cref{thm:two-sided-clt}
    and we obtain that $(F_{Y_n})_{n\geq 1}$ converges pointwise on $\R\setd{0}$ to $F_Y$.
    Since $\R\setminus\{0\}$ is dense in $\R$,
    it follows from \cite[Lemma~3(i) p.~275]{chow2003probability} that $Y_n$ converges to $Y$ in distribution.
    Since $Y$ is atomless, \Cref{lemma:bound-supremum} yields
    $\sup_{y\in \R} \P{Y_n = y} \to 0$, and combined with \eqref{bound-p_n} we obtain that $p_n \to 0$.
\end{proof}

\begin{proof}[Proof of \Cref{corol:estimator-stickiness}]
    \Cref{assum:clt-2} clearly implies that $\alphas \in \operatorname{int}(F)$ and $\E{\ell_+'(d(\alphas, X))^2} < +\infty$.
    By considering cases according to the value of $|I_0|$ and applying \Cref{thm:estimator-stickiness},
    we obtain that $\P{\hR_n = |I_0|} \xrightarrow[n\to\infty]{} 1$.

    Since $\hR_n$ and $|I_0|$ have values in $\set{0,1,2}$,
    $$\E[\big]{\big|\hR_n - |I_0|\big|} = \E[\big]{\indic{\hR_n \neq |I_0|}\big|\hR_n - |I_0|\big|}
    \le 2 \P{\hR_n \neq |I_0|}\xrightarrow[n\to \infty]{}0.$$
\end{proof}

\subsection{Proofs for Section~\ref{sec:empirical-difference-quotients}}

\begin{lemma}
    \label{lemma:estimator-measurable}
    Let $n\ge 2$. The mapping $\hR_n':\Omega \to \set{0,1,2}$ is measurable.
\end{lemma}

\begin{proof}
    We define the mapping 
    \begin{align*}
        g_n \colon \Omega \times T &\to [0,\infty) \\
        (\omega, \beta) &\mapsto 
        |d(\beta, \halpha_{A_n}^\omega) - r_n| 
        + \big(\hphi_{B_n}^\omega(\beta) - \hphi_{B_n}^\omega(\halpha_{A_n}^\omega)- r_n\tau_n\big)_+
    \end{align*}
    and we observe that $g_n$ is Carathéodory in the sense of \cite[Definition~4.50]{aliprantis2006infinite}.
    Let us also introduce the random set $\omega \mapsto C_n^\omega \coloneqq \set{\beta \in T: g_n(\omega, \beta) = 0}$.
    For each $\omega\in \Omega$, $C_n^\omega$ is closed.
    By \cite[Corollary~18.8]{aliprantis2006infinite}, $C_n$ has measurable graph. 
    Since $T$ is Polish and the probability space is complete,
    \cite[Theorem~1.3.3(iii)]{molchanov2017theory} implies that $C_n$ is a measurable closed set. 
    Moreover, by \cite[Proposition~1.1.15(i)]{molchanov2017theory}, its cardinality $|C_n|$ is measurable,
    hence so is $\hR_n'$.
\end{proof}

The proof of \Cref{thm:estimator-stickiness-2} is preceded by a series of lemmas.
Throughout, the assumptions of \Cref{thm:estimator-stickiness-2} are silently in force.

\begin{lemma}
    \label{lemma:good-sphere}
    The following holds:
    \begin{enumerate}
        \item $d(\halpha_{A_n},\alphas) = o_\PP(r_n)$.
        \item For each $n\ge 2$, let $\Omega_n = \set{d(\halpha_{A_n},\alphas) \le r_n/2}$, so that $\P{\Omega_n}\to 1$.
        There exists $N\ge 2$ such that for each $n\ge N$ and on the event $\Omega_n$, the sphere $\S_n$ coincides with $\set{\beta_i: 1\le i \le m}$,
        where each $\beta_i$ is in $\Delta_i$.
    \end{enumerate}
\end{lemma}

\begin{proof}
    1. In the case where $\alphas$ is sticky, it follows from $|A_n|\to \infty$ and \Cref{thm:sticky-lln} that $d(\halpha_{A_n},\alphas) = o_\PP(r_n)$.
    In the one-sided partly sticky case, where $I_0 = \set{i}$,
    the assumption that $(\Phi_i)_+'(t) / t^{a_i} \to K_i$ implies condition (b) of \Cref{thm:one-sided-clt}.
    Combined with the assumption that $r_n |A_n|^{1/(2a_i)} \to +\infty$, 
    we obtain from \Cref{thm:one-sided-clt} that $d(\halpha_{A_n},\alphas) = o_\PP(r_n)$.
    Let us turn to the two-sided partly sticky case, where $I_0=\set{i_-, i_+}$.
    For legibility, we write $a_-, a_+, K_-, K_+$ instead of $a_{i_-}, a_{i_+}, K_{i_-}, K_{i_+}$.
    We adapt the notation from \Cref{sec:two-sided-clt} by letting $\hm_{A_n} = S(\halpha_{A_n})$, so that
    $d(\halpha_{A_n},\alphas) = |\hm_{A_n}|$. 
    We fix $\varepsilon>0$ and we show that $\P{|\hm_{A_n}| \ge \varepsilon r_n}$ goes to $0$.
    This probability equals
    $\P{|A_n|^{1/(2a_+)}\hm_{A_n} \ge \varepsilon |A_n|^{1/(2a_+)} r_n}
    + \P{|A_n|^{1/(2a_-)}\hm_{A_n} \le -\varepsilon |A_n|^{1/(2a_-)} r_n}$.
    By assumption, $|A_n|^{1/(2a_+)} r_n$ and $|A_n|^{1/(2a_-)} r_n$ both go to $+\infty$.
    Let us consider an arbitrary $M>0$ and note that for $n$ large enough we have 
    $\varepsilon r_n \min(|A_n|^{1/(2a_+)}, |A_n|^{1/(2a_-)}) > M$.
    In that case, 
    \begin{equation}
        \label{eq:deviation-hm}
        \P{|\hm_{A_n}| \ge \varepsilon r_n} \le \P{|A_n|^{1/(2a_+)}\hm_{A_n} > M} + \P{|A_n|^{1/(2a_-)}\hm_{A_n} \le -M}.
    \end{equation}
    By \Cref{lemma:equiv-assum}, condition (b) of \Cref{thm:two-sided-clt} is fulfilled.
    Letting $n\to \infty$ in \Cref{eq:deviation-hm} and applying \Cref{thm:two-sided-clt} yields
    $\limsup_n \P{|\hm_{A_n}| \ge \varepsilon r_n} \le 1-\P{Z\le \frac{K_+ }{\sigma}M^{a_+}} +  \P{Z\le -\frac{K_- }{\sigma}M^{a_-}}$.
    Letting $M\to \infty$, it follows that $\limsup_n \P{|\hm_{A_n}| \ge \varepsilon r_n} = 0$, hence the conclusion.

    2. 
    Since $r_n\to 0$, we may pick $N\ge 2$ such that $n\ge N \implies r_n \le \frac{2\rho_0}{3}$.
    Let $n\ge N$ be arbitrary.
    In the case where $\halpha_{A_n} = \alphas$, we let $\beta_i \coloneqq \gamma_i(r_n)$ for each $i\in [m]$.
    In the case where $\halpha_{A_n} \neq \alphas$, we denote by $j_n \in [m]$ the index such that $\halpha_{A_n} \in \Delta_{j_n}$
    and on $\Omega_n$ 
    we have $r_n + d(\halpha_{A_n},\alphas) \le 3r_n/2 \le \rho_0$ and 
    we define $\beta_{j_n} \coloneqq \gamma_{j_n}(r_n + d(\halpha_{A_n},\alphas))$ as well as 
    $\beta_{i} \coloneqq \gamma_{i}(r_n - d(\halpha_{A_n},\alphas))$ for each $i\in [m]\setd{j_n}$.
    In both cases, it holds on $\Omega_n$ that $\set{\beta_i: 1\le i \le m} \subset \S_n$, with $\beta_i \in \Delta_i$.
    Conversely, on $\Omega_n$, $\S_n \subset \overline B(\alphas, \rho_0)$,
    hence $\S_n$ has cardinality $m$ and $\S_n = \set{\beta_i: 1\le i \le m}$.
\end{proof}

\begin{lemma}
    \label{lemma:quotient-closed-form}
    For each $i\in [m]$, we define the function 
    \begin{align*}
        s_i \colon (0,\rho_0] &\to [0,+\infty) \\
        t &\mapsto \frac{\phi(\gamma_i(t)) - \phi(\alphas)}{t}.
    \end{align*}
    \begin{enumerate}
        \item $s_i$ is increasing and $s_i(t) \xrightarrow[t\downarrow 0]{} \phi_{\Delta_i}'(\alphas)$.
        \item For each $n\ge N$ and on $\Omega_n$ we define $Q_{n,i} = \frac{\phi(\beta_i) - \phi(\halpha_{A_n})}{r_n}$.
        \begin{enumerate}
            \item In the case where $\halpha_{A_n} = \alphas$, $Q_{n,i} = s_i(r_n)$ for each $i\in [m]$.
            \item In the case where $\halpha_{A_n} \neq \alphas$, 
            \begin{align*}
                Q_{n,j_n} &= \Big(1+\frac{d(\halpha_{A_n},\alphas)}{r_n}\Big) s_{j_n}(r_n + d(\halpha_{A_n},\alphas)) 
                - \frac{d(\halpha_{A_n},\alphas)}{r_n} s_{j_n}( d(\halpha_{A_n},\alphas)),
                \\
                Q_{n,i} &= \Big(1-\frac{d(\halpha_{A_n},\alphas)}{r_n}\Big) s_{i}(r_n - d(\halpha_{A_n},\alphas)) 
                - \frac{d(\halpha_{A_n},\alphas)}{r_n} s_{j_n}( d(\halpha_{A_n},\alphas))
                \text{ for } i\neq j_n.
            \end{align*}
        \end{enumerate}
    \end{enumerate}
\end{lemma}

\begin{proof}
    The first item follows from the convexity of $t\mapsto \phi(\gamma_i(t))$ and the second item is straightforward algebra.
\end{proof}

\begin{lemma}
    \label{lemma:quotient-separation}
    \begin{enumerate}
        \item Assume that $I_0 \neq \emptyset$ and let $i\in I_0$. There exists $C_i\ge 0$ such that for each $n\ge N$ and on $\Omega_n$: $$Q_{n,i} \le C_i r_n^{a_i}.$$
        \item Assume that $I_0^c \neq \emptyset$. Let $\epsilon = \min_{i\notin I_0} \phi_{\Delta_i}'(\alphas)$,
        $M = \max_{1\le i\le m} s_i(\rho_0/3)$ and
        $\Omega_n' = \set[\big]{d(\halpha_{A_n},\alphas) \le \frac{\epsilon}{2(\epsilon + M)} r_n}$.
        For each $n\ge N$ and on $\Omega_n \cap \Omega_n'$: $$\forall i\in I_0^c, \,Q_{n,i} \ge \frac \epsilon 2.$$
        \item There is some $N'\ge 2$ that satisfies the following. For each $n\ge N'$ and $i\in I_0$, 
         $C_i r_n^{a_i} \le \frac{\tau_n}3$. Additionally, if $I_0^c \neq \emptyset$, $n\ge N' \implies \frac \epsilon 2 > \frac{4\tau_n}{3}$.
    \end{enumerate}
\end{lemma}

\begin{proof}
    1. Let $i\in I_0$.
    Since $\Phi_i$ is convex, for every $t\in (0,\rho_0)$, $s_i(t) = (\Phi_i(t) - \Phi_i(0))/t \le (\Phi_i)_+'(t)$. 
    By condition (b) of \Cref{thm:estimator-stickiness-2}, 
    there exists $A_i \ge 0$ and $\delta_i \in (0,\rho_0)$ such that $\forall t\in (0, \delta_i]$, $s_i(t) \le A_i t^{a_i}$.
    For each $t\in [\delta_i, \rho_0]$, we have the estimate $s_i(t) / t^{a_i} \le s_i(\rho_0) / {\delta_i}^{a_i}$.
    We let $B_i = \max(A_i, s_i(\rho_0) / {\delta_i}^{a_i})$, so that 
    $s_i(t) \le B_i t^{a_i}$ holds for each $t\in (0,\rho_0]$.

    Fix some $n\ge N$.
    In the case where $\halpha_{A_n} = \alphas$, we have on $\Omega_n$ that $Q_{n,i} = s_i(r_n) \le B_i r_n^{a_i}$ for each $i\in I_0$.
    Assume now that $\halpha_{A_n} \neq \alphas$. 
    If $j_n \in I_0$, by \Cref{lemma:quotient-closed-form} it holds on $\Omega_n$ that 
    $$Q_{n,j_n} \le \Big(1+\frac{d(\halpha_{A_n},\alphas)}{r_n}\Big) s_{j_n}(r_n + d(\halpha_{A_n},\alphas)) \le \frac 32 s_{j_n}(3r_n/2) 
    \le \Big(\frac{3}{2}\Big)^{a_{j_n}+1} B_{j_n} r_n^{a_{j_n}}.$$
    Similarly, for each $i\in I_0\setd{j_n}$, $$Q_{n,i} \le \Big(1-\frac{d(\halpha_{A_n},\alphas)}{r_n}\Big) s_{i}(r_n - d(\halpha_{A_n},\alphas)) \le s_{i}(r_n) \le B_i r_n^{a_i}.$$
    Finally, we define $C_i = (3/2)^{a_{i}+1} B_i$ for each $i\in I_0$ and we have shown that $Q_{n,i} \le C_i r_n^{a_i}$.

    2. Let $n\ge N$ be arbitrary. 
    In the case where $\halpha_{A_n} = \alphas$, it holds on $\Omega_n$ that $Q_{n,i} = s_i(r_n) \ge \phi_{\Delta_i}'(\alphas) \ge \epsilon$ for each $i\in I_0^c$.
    Assume next that $\halpha_{A_n} \neq \alphas$. If $j_n \in I_0^c$, by \Cref{lemma:quotient-closed-form} we have on $\Omega_n$ that 
    \begin{align*}
        Q_{n,j_n} 
    &= s_{j_n}(r_n + d(\halpha_{A_n},\alphas)) + \tfrac{d(\halpha_{A_n},\alphas)}{r_n} \big[s_{j_n}(r_n + d(\halpha_{A_n},\alphas)) - s_{j_n}(d(\halpha_{A_n},\alphas))\big]
    \\&\ge s_{j_n}(r_n + d(\halpha_{A_n},\alphas)) \ge \epsilon.
    \end{align*}
    For each $i\in I_0^c\setd{j_n}$ and on $\Omega_n \cap \Omega_n'$, 
    $d(\halpha_{A_n},\alphas) \le r_n/2 \le \rho_0/3$, hence 
    $$Q_{n,i} \ge \Big(1-\frac{d(\halpha_{A_n},\alphas)}{r_n}\Big) \epsilon
    - \frac{d(\halpha_{A_n},\alphas)}{r_n} M \ge \frac \epsilon 2.$$

    3. This follows from $r_n^{a_i}/\tau_n \to 0$ for each $i\in I_0$, as well as from $\tau_n \to 0$.
\end{proof}

\begin{lemma}
    \label{lemma:bound-quotient-convexity}
    Let $n\ge N$. Consider $\beta, \talpha \in T$ that satisfy $d(\beta, \talpha) = r_n$ and $d(\talpha, \alphas)\le r_n/2$.
    Define the random variable $U = \frac{\ell(d(\beta, X)) - \ell(d(\talpha,X))}{r_n}$.
    Then $|U| \le \ell_+'(d(\alphas, X) + \rho_0)$.
\end{lemma}

\begin{proof}
    By the convexity estimate \eqref{eq:integrand-estimate}, we have $|U| \le \ell_+'\big(\max[d(\beta, X), d(\talpha, X)]\big)$.
    We observe that $$d(\beta, X) \le  d(\beta,\talpha) + d(\talpha, \alphas) + d(\alphas, X) \le d(\alphas, X) + 3r_n/2 \le d(\alphas, X) + \rho_0$$
    and
    $$d(\talpha, X) \le d(\talpha, \alphas) + d(\alphas, X) \le  d(\alphas, X) + r_n/2 \le d(\alphas, X) + \rho_0,$$
    hence the claimed inequality.
\end{proof}

\begin{lemma}
    \label{lemma:deviation-quotients}
    For each $n\ge N$ and on $\Omega_n$, let $\widehat Q_{n,i}=\frac{\hphi_{B_n}(\beta_i) - \hphi_{B_n}(\halpha_{A_n})}{r_n}$ for every $i\in [m]$.
    Define also the subset $\Omega_n''$ of $\Omega_n$ by $\Omega_n'' = \set[\big]{\max_{1\le i \le m} |\widehat Q_{n,i} - Q_{n,i}| \le \tau_n/3}$
    and let $\mathtt v = \E[\big]{\ell_+'(d(\alphas, X) + \rho_0)^2}$.
    Then $\P{\Omega_n''} \ge 1-\frac{9 \mathtt v m}{|B_n| \tau_n^2} - \P{\Omega_n^c}$ for every $n\ge N$.
\end{lemma}

\begin{proof}
    Fix some $n\ge N$. For short, we denote by $\mathbf X_{A_n}$ the $|A_n|$-tuple of observations indexed by $A_n$, and similarly we use the notation $\mathbf X_{B_n}$.
    By \Cref{lemma:Borel-selection}, we may write $\Omega_n = \set{\mathbf X_{A_n} \in D_n}$ for some Borel subset $D_n$,
    thus
    $$\indic{\Omega_n''^c} = \indic{\mathbf X_{A_n} \in D_n, \max_{1\le i \le m} |\widehat Q_{n,i} - Q_{n,i}| > \tau_n/3} + \indic{\mathbf X_{A_n} \notin D_n} = g_n(\mathbf X_{A_n}, \mathbf X_{B_n})$$
    for some Borel function $g_n$.
    Let us define the function $h_n: T^{|A_n|} \to \R, \mathbf x \mapsto \E{g_n(\mathbf x, \mathbf X_{B_n})}$. 
    In the case where $\mathbf x\in D_n$, 
    we let $\talpha = f_{|A_n|}(\mathbf x)$ (using the notation from \Cref{lemma:Borel-selection}),
    and for each $i\in [m]$
    we apply \Cref{lemma:bound-quotient-convexity} and Chebyshev's inequality 
    to the random variable 
    $\frac{1}{|B_n|}\sum_{k\in B_n} \frac{\ell(d(\beta_i, X_k)) - \ell(d(\talpha,X_k))}{r_n}$,
    where each summand has finite variance by \Cref{assum:clt-2}.
    The union bound then yields
    $h_n(\mathbf x) \le \frac{9 \mathtt v m}{|B_n| \tau_n^2}$.
    In the case where $\mathbf x\notin D_n$, we have $h_n(\mathbf x)=1$.
    Since $\mathbf X_{A_n}$ and $\mathbf X_{B_n}$ are independent, it follows that 
    $$\P{\Omega_n''^c} = \E{\E{g_n(\mathbf X_{A_n}, \mathbf X_{B_n}) \mid \mathbf X_{A_n}}} = \E{h_n(\mathbf X_{A_n})} \le \frac{9 \mathtt v m}{|B_n| \tau_n^2} + \P{\Omega_n^c}.$$
\end{proof}

\begin{proof}[Proof of \Cref{thm:estimator-stickiness-2}]
    Let us fix some $n\ge \max(N,N')$.

    First, we consider the case where $I_0 \neq \emptyset$.
    On $\Omega_n''$ we have by \Cref{lemma:quotient-separation} that for each $i\in I_0$, $\widehat Q_{n,i} \le Q_{n,i} + \tau_n/3 \le 2\tau_n/3 \le \tau_n$.
    Assume now that $I_0^c \neq \emptyset$. 
    On $\Omega_n' \cap \Omega_n''$ it holds by \Cref{lemma:quotient-separation} that for every $i\in I_0^c$,
    $\widehat Q_{n,i} \ge Q_{n,i} - \tau_n/3 \ge \epsilon/2 - \tau_n/3 > \tau_n$.
    It follows that when $I_0 \neq \emptyset$ and $I_0^c \neq \emptyset$, the inclusion $\Omega_n' \cap \Omega_n'' \subset \set{\hR_n' = |I_0|}$ holds.
    Assume on the other hand that $I_0^c = \emptyset$, \ie $I_0 = [m]$, which entails $1\le m \le 2$. 
    Then $\Omega_n'' \subset \set{\hR_n' = |I_0|}$ is true as well.

    Finally, we deal with the case $I_0 = \emptyset$.
    On $\Omega_n' \cap \Omega_n''$ we already observed that for every $i\in I_0^c$, $\widehat Q_{n,i} > \tau_n$,
    thus $\Omega_n' \cap \Omega_n'' \subset \set{\hR_n' = |I_0|}$.

    We proved that $\Omega_n'' \subset \set{\hR_n' = |I_0|}$ or $\Omega_n' \cap \Omega_n'' \subset \set{\hR_n' = |I_0|}$.
    Since $d(\halpha_{A_n},\alphas) = o_\PP(r_n)$ we have the convergence $\P{\Omega_n'}\to 1$, and 
    since $|B_n| \tau_n^2 \to \infty$ we also have $\P{\Omega_n''}\to 1$, hence finally $\P{\hR_n' = |I_0|}\to 1$.
    The $L^1$ convergence is obtained as in \Cref{thm:estimator-stickiness}.
\end{proof}

\subsection{Proofs for Section \ref{sec:medians}}

\begin{proof}[Proof of \Cref{prop:endpoints-median}]
    Assume for the sake of contradiction that $\alpha_+$ is neither in $D$ nor in the support. 
    We consider two cases according to the cardinality of $\dir_{\alpha_+}$.

    First, $|\dir_{\alpha_+}|=1$, meaning that $\alpha_+$ is a leaf and $\dir_{\alpha_+} = \set[\big]{T\setd{\alpha_+}}$.
    The first-order optimality condition at $\alpha_+$ is $1-2\mu(T\setd{\alpha_+}) \ge 0$, hence $\mu(\set{\alpha_+})\ge 1/2$, which contradicts $\alpha_+ \notin \supp \mu$.

    Next, $|\dir_{\alpha_+}|=2$. We denote by $v,w$ the neighboring vertices of $\alpha_+$ and we assume \wlo that 
    \begin{equation}
        \label{eq:disjoint}
        [v,\alpha_+)\cap M_1(\mu)=\emptyset.
    \end{equation}
    Indeed, if $\alpha_+ \neq \alpha_-$, we relabel $v$ and $w$ so that $v\notin \Delta_{\alpha_+ \to \alpha_-}$.
    Since $T\setminus \supp \mu$ is open, so is the intersection $\left(T\setminus \supp \mu\right) \cap (v,w)$, which contains $\alpha_+$.
    Consequently, we can pick some $\beta\in \left(T\setminus \supp \mu\right) \cap (v,\alpha_+)$ that satisfies 
    $\mu\left([\beta,\alpha_+]\right)=0$.
    This last equality implies $\mu(\Delta_{\beta\to v}) = \mu(\Delta_{\alpha_+\to v})$ and $\mu(\Delta_{\beta\to w}) = \mu(\Delta_{\alpha_+\to w})$,
    thus $\phi_v'(\beta) = \phi_v'(\alpha_+)$ and $\phi_w'(\beta) = \phi_w'(\alpha_+)$.
    Therefore $\beta$ satisfies the same optimality conditions as $\alpha_+$ and we must have $\beta \in M_1(\mu)$.
    However, by construction $\beta \in (v,\alpha_+)$, which contradicts \eqref{eq:disjoint}.
    We proceed identically with $\alpha_-$.
\end{proof}

\begin{proof}[Proof of \Cref{prop:unique-median}]
    1. The subsets $T\setminus \Delta_{\alpha_-\to \alpha_+}$ and $T\setminus \Delta_{\alpha_+\to \alpha_-}$ 
    are clearly closed and convex; let us show they are disjoint. 
    For the sake of contradiction, assume there exists some $\alpha$ in the intersection.
    Then $\alpha_- \in [\alpha, \alpha_+)$ and $\alpha_+ \in [\alpha, \alpha_-)$,
    hence $d(\alpha, \alpha_+) = d(\alpha, \alpha_-) + d(\alpha_-, \alpha_+) > d(\alpha, \alpha_-)$,
    and by symmetry $d(\alpha, \alpha_-) > d(\alpha, \alpha_+)$, a contradiction.
    A necessary optimality condition for $\alpha_-$ is $\phi_{\alpha_+}'(\alpha_-)\ge 0$, 
    which we rewrite as $\mu(T\setminus \Delta_{\alpha_-\to \alpha_+})\ge \frac 12$.
    Symmetrically, $\mu(T\setminus \Delta_{\alpha_+\to \alpha_-})\ge \frac 12$ 
    hence $\mu(T\setminus \Delta_{\alpha_-\to \alpha_+}) = \mu(T\setminus \Delta_{\alpha_+\to \alpha_-}) = \frac 12$.

    2. We rewrite the assumption as $\phi_{\alpha_+}'(\alpha_-) = \phi_{\alpha_-}'(\alpha_+)=0$ 
    and it follows from \Cref{prop:sign-derivative} that $\set{\alpha_-,\alpha_+}\subset M_1(\mu)$.
    Since $M_1(\mu)$ is convex, this implies $[\alpha_-,\alpha_+]\subset M_1(\mu)$.

    3. The forward implication follows trivially from the first item. 
    Conversely, assume the existence of such $C_1,C_2$ and pick $c_1\in C_1$, $c_2\in C_2$.
    The set $[c_1,c_2]\cap C_1$ is a closed convex subset of $[c_1,c_2]$, hence it is a geodesic segment of the form $[c_1, \alpha_-]$,
    where $\alpha_-\in C_1$.
    Similarly, $[c_1,c_2]\cap C_2 = [c_2, \alpha_+]$ with $\alpha_+\in C_2$.
    Since $C_1$ and $C_2$ are disjoint, we have $\alpha_- \neq \alpha_+$.
    By definition, $C_1 \cap (\alpha_-, \alpha_+] = C_2 \cap [\alpha_-, \alpha_+)= \emptyset$.
    Since $C_1$ is convex, it follows that $C_1 \subset T\setminus \Delta_{\alpha_-\to \alpha_+}$,
    and similarly $C_2 \subset T\setminus \Delta_{\alpha_+\to \alpha_-}$.
    Since $T\setminus \Delta_{\alpha_-\to \alpha_+}$ and $T\setminus \Delta_{\alpha_+\to \alpha_-}$ are disjoint, 
    we obtain that $\mu(T\setminus \Delta_{\alpha_-\to \alpha_+}) = \mu(T\setminus \Delta_{\alpha_+\to \alpha_-}) = \frac 12$,
    thus $\phi_{\alpha_+}'(\alpha_-) = \phi_{\alpha_-}'(\alpha_+)=0$ and by the second item $[\alpha_-,\alpha_+]\subset M_1(\mu)$.
\end{proof}

\begin{proof}[Proof of \Cref{lemma:binomial}]
    If $k>n$, then $f$ is identically zero. 
    If $k=1$, then $n=1$ and the claim holds trivially.
    We may therefore assume next that $1<k\le n$.
    By the integral representation for the CDF of the binomial distribution (see, \eg \cite[Exercise 21 p.~204]{shiryaev2016proba}), it holds that
    $f(p) = \P{n-\Sigma_{n,p}\le n-k} = k\binom{n}{k}\int_{0}^{p}x^{k-1}(1-x)^{n-k}\diff x$.
    The second derivative of $p\mapsto \int_{0}^{p}x^{k-1}(1-x)^{n-k}\diff x$ is $p^{k-2}(1-p)^{n-k-1}\big(k-1-p(n-1)\big)$.
    Since $p\le 1/2$ we find that 
    $p(n-1) + 1 \le (n+1)/2$, and since $k>n/2$ it follows that $(n+1)/2 \le k$, hence 
    $f''\ge 0$ and $f$ is convex on $[0,1/2]$.
    Let us write $p\in [0,1/2]$ as the convex combination $p = (2p) \cdot \frac 12 + (1-2p) \cdot 0$, 
    from which it follows that $f(p)\le 2pf(1/2)$.
\end{proof}

\begin{proof}[Proof of \Cref{prop:control-derivatives-median}]
    For $1\le i\le m$, define $N_i = n\hmu_n(\Delta_i)$ and $N_0 = n\hmu_n(\set{\alphas})$,
    so that $(N_0,\ldots,N_m)$ has a multinomial distribution with $n$ trials and $m+1$ outcomes.
    Note the equality of events 
    $\set{\min_{1\le i\le m} \hphi_{\Delta_i}'(\alphas)\le -t} = \bigcup_{i=1}^m \set{N_i\ge n(1+t)/2}$,
    where the $\set{N_i\ge n(1+t)/2}$ are pairwise disjoint since $t>0$.
    Define $k=\ceil[\big]{\frac{n(1+t)}{2}}$ and 
    observe that $\P{N_i\ge n(1+t)/2} = \P{\Sigma_{n, \mu(\Delta_i)} \ge k}$,
    with the notation from \Cref{lemma:binomial}.
    Since $\alphas\in M_1(\mu)$, we note additionally that $0\le \mu(\Delta_i) \le 1/2$ and
    $$\P*{\min_{1\le i\le m} \hphi_{\Delta_i}'(\alphas)\le -t} = \sum_{i=1}^{m} f(\mu(\Delta_i)).$$
    By \Cref{lemma:binomial}, 
    \begin{equation}
        \label{eq:final-bound-convexity}
        \P*{\min_{1\le i\le m} \hphi_{\Delta_i}'(\alphas)\le -t} \stackrel{(i)}{\le} 2f(1/2) \sum_{i=1}^{m} \mu(\Delta_i) 
        \stackrel{(ii)}{\le} 2f(1/2) = 2 \P*{\Sigma_{n, 1/2}\ge \frac{n(1+t)}{2}}.
    \end{equation}
    Hoeffding's inequality then yields the final inequality. 
\end{proof}

\begin{proof}[Proof of \Cref{prop:region}]
    1. The first equality follows from $\hphi_{\Delta}'(\alpha) = 1-2\hmu_n(\Delta)$ and $n\hmu_n(\Delta)= \sum_{k=1}^{n} \indic{X_k\in \Delta}$.
    Next, we show that $\hC_t$ is a subset of the right-hand side of \eqref{eq:region-equality}.
    Let us fix $\alpha\in \hC_t$ and $K$ a subset of $[n]$ with cardinality $s$.
    By the first reformulation of $\hC_t$, the $s$ points $(X_k)_{k\in K}$ cannot all lie in a single $\Delta \in \dir_\alpha$, hence $\alpha$ belongs to 
    $\conv(\set{X_k:k\in K})$.
    For the opposite inclusion, we proceed by contrapositive and we consider an $\alpha \in T\setminus \hC_t$.
    We can then find a direction $\Delta\in \dir_\alpha$ such that $\sum_{k=1}^{n} \indic{X_k\in \Delta} \ge s$,
    hence there exists $K\subset [n]$ with cardinality $s$ such that $X_k\in \Delta$ for each $k\in K$.
    Since $\Delta$ is convex, $\conv(\set{X_k:k\in K}) \subset \Delta \subset T\setd{\alpha}$ and therefore $\alpha\notin \conv(\set{X_k:k\in K})$.

    2. The optimality condition of \Cref{prop:optimality} applied to the measure $\hmu_n$ entails that $M_1(\hmu_n) \subset \hC_t$, hence $\hC_t$ is nonempty.
    By \Cref{prop:localization-support} and \eqref{eq:region-equality}, $\hC_t$ is an intersection of closed convex subsets, thus it is closed and convex.
    It also follows from \eqref{eq:region-equality} that $\hC_t$ is bounded. Since $T$ is proper, we deduce that $\hC_t$ is compact.

    It is straightforward to show that for a given $K\subset [n]$, $$\conv(\set{X_k:k\in K}) = \bigcup_{(i,j)\in K^2} [X_i, X_j].$$
    By \eqref{eq:region-equality} and \cite[Theorem~1.3.25(iv)]{molchanov2017theory}, the measurability of $\hC_t$ follows from the measurability of $[X_i, X_j]$
    for a fixed $(i,j)$, which we shall establish next.
    First, since $T$ is a separable metric space, the Borel $\sigma$-algebra on $T^2$ coincides with the product $\sigma$-algebra $\B(T)\otimes \B(T)$,
    hence $(X_i,X_j)$ is a measurable $T^2$-valued random element.
    We denote by $\mathcal K'$ the collection of nonempty compact subsets of $T$, by $d_H$ the Hausdorff distance on $\mathcal K'$
    and we define the mapping $g:T^2 \to \mathcal K'$, $(x,y)\mapsto [x,y]$. 
    We show the estimate $d_H([x,y], [x',y'])\le \max(d(x,x'), d(y,y'))$ for $x,x',y,y'\in T$, which will imply that $g$ is continuous.
    Let us denote by $\gamma:[0,1]\to [x,y]$ the geodesic from $x$ to $y$ and by $\gamma':[0,1]\to [x',y']$ the geodesic from $x'$ to $y'$.
    For each $t\in [0,1]$, $d(\gamma_t, [x',y']) \le d(\gamma_t, \gamma_t')$ and since the space $T$ is Hadamard it holds that 
    $d(\gamma_t, \gamma_t') \le (1-t)d(x,x') + td(y,y') \le \max(d(x,x'), d(y,y'))$ (see, \eg \cite[Proposition~1.1.5]{bacak2014convex}). 
    Consequently, $\sup_{z\in [x,y]} d(z, [x',y']) \le \max(d(x,x'), d(y,y'))$.
    Similarly, $\sup_{z'\in [x',y']} d(z', [x,y])\le \max(d(x,x'), d(y,y'))$, hence the claim.
    Since $g$ is continuous, it is Borel measurable, hence $[X_i,X_j]$ is a measurable $\mathcal K'$-valued random element.
    By \cite[Theorem~D.6(iii)]{molchanov2017theory}, the Borel $\sigma$-algebra on $\mathcal K'$ coincides with the trace of the Effros $\sigma$-algebra,
    hence $[X_i,X_j]$ is a measurable closed subset.

    3. Let us write $M_1(\mu) = [\alpha_-, \alpha_+]$, so that $\set{M_1(\mu) \subset \hC_t} = \set{\alpha_-\in \hC_t} \cap \set{\alpha_+\in \hC_t}$.
    We observe that $\set{\alpha_-\in \hC_t} = \bigcap_{\Delta\in \dir_{\alpha_-}} \set[\big]{\sum_{k=1}^{n} \indic{X_k \in \Delta} < \frac{n(1+t)}{2}}$,
    which belongs to the $\sigma$-algebra, and so does $\set{\alpha_+\in \hC_t}$. It follows that $\set{M_1(\mu) \subset \hC_t}$ is measurable.

    In the case where $M_1(\mu) = \set{\alphas}$, we may write $\set{M_1(\mu) \not\subset \hC_t} = \set{\alphas\notin \hC_t}$ 
    and apply \Cref{prop:control-derivatives-median} to obtain the claim.
    We assume next that $M_1(\mu) = [\alpha_-, \alpha_+]$ with $\alpha_-\neq \alpha_+$.
    Let us introduce the subsets
    $T_- = T\setminus \Delta_{\alpha_- \to \alpha_+}$, $T_+ = T\setminus \Delta_{\alpha_+ \to \alpha_-}$,
    the random variables
    $N_- = n\hmu_n(T_-)$, $N_+ = n\hmu_n(T_+)$ and 
    the event $\Omega_n = \bigcap_{k=1}^n \set{X_k\in T_- \cup T_+}$.
    
    Let us show the equality of events
    \begin{equation}
        \label{eq:events-median}
        \set{M_1(\mu) \subset \hC_t} \cap \Omega_n = \set{N_- < s, N_+ < s} \cap \Omega_n.
    \end{equation}
    We consider $\omega\in \Omega$ such that $N_-^{\omega} < s$, $N_+^{\omega} < s$ and $X_1^\omega,\ldots,X_n^\omega \in T_- \cup T_+$.
    Let us pick some $\beta\in M_1(\mu)$ and an arbitrary $\Delta \in \dir_\beta$.
    We argue that $\Delta$ intersects at most one of $T_-$ or $T_+$. 
    Indeed, if $x_-\in T_-$ and $x_+\in T_+$, 
    then the geodesic $[x_-,x_+]$ contains the whole segment $[\alpha_-,\alpha_+]$,
    and hence contains $\beta$. 
    Thus $x_-$ and $x_+$ cannot both belong to $\Delta$.
    Consequently, $n \hmu_n^\omega(\Delta) 
    = n \hmu_n^\omega(\Delta \cap (T_-\cup T_+)) 
    \le \max(n\hmu_n^\omega(T_-), n\hmu_n^\omega(T_+)) = \max(N_-^\omega, N_+^\omega) < s$.
    By the first item, it follows that $\beta\in \hC_t^\omega$, hence $M_1(\mu) \subset \hC_t^\omega$.
    Conversely, let us consider $\omega\in \Omega_n$ such that $M_1(\mu) \subset \hC_t^\omega$.
    Since $\alpha_- \neq \alpha_+$, we may pick some $\beta \in (\alpha_-, \alpha_+)$, and $\beta$ is in $\hC_t^\omega$.
    We define the directions $\Delta_- = \Delta_{\beta \to \alpha_-}$ and $\Delta_+ = \Delta_{\beta \to \alpha_+}$.
    Let us observe that $T_- \subset \Delta_-$, $T_+ \subset \Delta_+$ and $\Delta_- \cap \Delta_+ = \emptyset$.
    Indeed, if $x\in T_-$ then $\alpha_- \in [x, \alpha_+]$, and since $\beta \in (\alpha_-, \alpha_+)$ we have $\beta \notin [x,\alpha_-]$,
    hence $x\in \Delta_-$.
    Consequently, $\Delta_- \cap (T_-\cup T_+) = T_-$ and $\Delta_+ \cap (T_-\cup T_+) = T_+$.
    Since $\omega$ is in $\Omega_n$, we have additionally $\hmu_n^\omega(\Delta_-) = \hmu_n^\omega(\Delta_-\cap (T_-\cup T_+)) = \hmu_n^\omega(T_-)$,
    hence $n\hmu_n^\omega(\Delta_-) = N_-^\omega$ and similarly $n\hmu_n^\omega(\Delta_+) = N_+^\omega$.
    It follows from the first item that $n\hmu_n^\omega(\Delta_-) < s$ and $n\hmu_n^\omega(\Delta_+) < s$,
    hence $N_-^\omega<s$ and $N_+^\omega<s$.
    
    By \Cref{prop:unique-median}, $\Omega_n$ has probability $1$ and by combining with \eqref{eq:events-median} 
    we obtain that $\P{M_1(\mu) \not\subset \hC_t} = \P{(N_- \ge s) \cup (N_+ \ge s)}$.
    Since $N_- + N_+ \le  n$ and $s>n/2$, the events $\set{N_- \ge s}$ and $\set{N_+ \ge s}$ are disjoint.
    By \Cref{prop:unique-median}, $\mu(T_+) = 1/2$, thus $N_+$ and $N_-$ each follow $\Bin(n, 1/2)$
    and $\P{M_1(\mu) \not\subset \hC_t} = 2\P{\Sigma \ge s}$.
\end{proof}

\begin{proof}[Proof of \Cref{thm:confidence-region}]
    Note that $\us$ is well-defined because $\eta \ge 1/2^{n-1}$.
    Since $\floor{n/2} < \us \le n$ we have $\ut \in (0,1]$ and we may apply \Cref{prop:region}
    to obtain $\P{M_1(\mu) \subset \hC_\ut} \ge 1 -2 \P{\Sigma \ge \us} \ge 1-\eta$.
\end{proof}

\begin{proof}[Proof of \Cref{prop:region-shrink}]
    1. We fix an arbitrary $r\in (0,+\infty)$ and we denote by $\S_r$ the sphere centered at $\alphas$ with radius $r$,
    which is finite by \Cref{lemma:ball}.
    We also let $B_r = B(\alphas, r)$ and $s_n = \ceil*{\frac{n(1+t_n)}{2}}$.
    
    Let us prove by contraposition that 
    \begin{equation}
        \label{eq:inclusion}
        \bigcap_{\beta \in \S_r} \set{\hmu_n(\Delta_{\beta \to \alphas}) \ge s_n/n} \subset \set{\hC_{t_n} \subset B_r}.
    \end{equation}
    We suppose there exists $\alpha \in \hC_{t_n}$ such that $\alpha \notin B_r$, \ie $d(\alpha, \alphas) \ge r$.
    We therefore find $\beta \in \S_r \cap [\alpha, \alphas]$ and we observe that $\Delta_{\beta \to \alphas} \subset \Delta_{\alpha \to \alphas}$.
    By combining this inclusion with the first item of \Cref{prop:region}, we get
    $\hmu_n(\Delta_{\beta \to \alphas}) \le \hmu_n(\Delta_{\alpha \to \alphas}) < s_n/n$,
    which finishes the proof of \eqref{eq:inclusion}.

    Consider an arbitrary $\beta \in \S_r$. 
    The strong law of large numbers yields $\hmu_n(\Delta_{\beta \to \alphas}) \xrightarrow[]{a.s.} \mu(\Delta_{\beta \to \alphas})$.
    Since $\beta \neq \alphas$ and $M_1(\mu) = \set{\alphas}$, $\beta$ is not a Fr\'echet median of $\mu$, hence by \Cref{prop:optimality}
    there exists a direction $\Delta \in \dir_\beta$ such that $\phi_\Delta'(\beta)<0$, \ie $\mu(\Delta) > 1/2$.
    \Cref{prop:sign-derivative} guarantees that $\alphas \in \Delta$, hence $\Delta = \Delta_{\beta \to \alphas}$
    and finally $\mu(\Delta_{\beta \to \alphas}) > 1/2$.
    Combined with the limit $s_n/n \to 1/2$, we obtain $\P{\liminf_n\set{\hmu_n(\Delta_{\beta \to \alphas}) \ge s_n/n}} = 1$, 
    and since $\S_r$ is finite,
    $$\P[\Big]{\liminf_n \bigcap_{\beta \in \S_r} \set{\hmu_n(\Delta_{\beta \to \alphas}) \ge s_n/n}} = 1.$$
    It then follows from \eqref{eq:inclusion} that
    $\P{\bigcup_{N\ge 1} \bigcap_{n\ge N} \set{\hC_{t_n} \subset B_r}} = 1$.
    Since $r$ was arbitrary, this entails
    $\P{\bigcap_{m\ge 1}\bigcup_{N\ge 1} \bigcap_{n\ge N} \set{\hC_{t_n} \subset B_{1/m}}} = 1$
    and finally $\P{\sup_{\alpha \in \hC_{t_n}} d(\alpha, \alphas) \to 0} = 1$.

    2. It suffices to show that $\ut \to 0$, \ie that $\us / n \to 1/2$ as $n\to \infty$.
    Hoeffding's inequality yields
    $\P{\Sigma \ge \frac n2 + t} \le \exp(-2t^2/n)$ for each $t\ge 0$.
    Consequently, the integer $$s_n'\coloneqq \ceil*{\frac n2 + \sqrt{\frac n2 \log \frac 2\eta}\,}$$
    satisfies $2\P{\Sigma \ge s_n'} \le \eta$ and for sufficiently large $n$ we also have 
    $s_n' \in \set{\floor{n/2}+1,\ldots,n}$, thus $\us \le s_n'$.
    Moreover, $\us \ge \floor{n/2}+1$, hence by squeezing we obtain that $\us / n \to 1/2$.
\end{proof}

\end{document}